\documentclass[11pt,3p, authoryear]{elsarticle}
\usepackage{graphicx}
\usepackage{algorithm}  
\usepackage{algorithmic}
\usepackage{amssymb,url,amsmath}
\usepackage{bm} % bold symbols
\usepackage[dvipsnames]{xcolor}
\usepackage{booktabs}	% nice tables
\usepackage{epstopdf}
\usepackage{color}
\usepackage{pifont}
\usepackage{tikz}
\usetikzlibrary{arrows}
\usepackage{eurosym}

\usepackage{pgfplots}
\usepackage{setspace}
\usepackage{rotating}
\usepackage{circuitikz}
\usepackage{bbm}

\newtheorem{proposition}{Proposition}

\newproof{pf}{Proof}
\tikzset{>=latex}

\usepackage{bm}
\newcommand*{\QEDA}{\hfill\ensuremath{\blacksquare}}

\newcommand{\CC}{\mathcal{C}}
\newcommand{\X}{\mathbb{X}}
\newcommand{\AC}{\mathcal{A}}
\newcommand{\SC}{\mathcal{S}}
\newcommand{\R}{\mathbb{R}}

\allowdisplaybreaks
\makeatletter
\def\ps@pprintTitle{%
   \let\@oddhead\@empty
   \let\@evenhead\@empty
   \def\@oddfoot{\reset@font\hfil\thepage\hfil}
   \let\@evenfoot\@oddfoot
}
\makeatother
\begin{document}
\begin{frontmatter}

\title{Enabling inter-area reserve exchange through stable benefit allocation mechanisms}
%WH

\author[a1]{Orcun Karaca\corref{cor1}}
\ead{okaraca@ethz.ch}
\author[a3]{Stefanos Delikaraoglou}
\ead{sdelikar@mit.edu}
\author[a2]{Gabriela Hug}
\ead{hug@ethz.ch}
\author[a1]{Maryam Kamgarpour}
\ead{maryamk@ethz.ch}
\cortext[cor1]{Corresponding author}
\address[a1]{Automatic Control Laboratory, D-ITET, ETH Z{\"u}rich, Switzerland}
\address[a3]{Laboratory for Information and Decision Systems, MIT, USA}
\address[a2]{Power Systems Laboratory, D-ITET, ETH Z{\"u}rich, Switzerland}

\begin{abstract}
The establishment of a single European day-ahead market has accomplished the integration of the regional day-ahead markets. However, the reserves provision and activation remain an exclusive responsibility of regional operators. This limited spatial coordination and the sequential market structure hinder the efficient utilization of flexible generation and transmission, since their capacities have to be ex-ante allocated between energy and reserves. To promote reserve exchange, recent work has proposed a preemptive model that defines the optimal inter-area transmission capacities for energy and reserves reducing the expected system cost. This decision-support tool, formulated as a stochastic bilevel program, respects the current architecture but does not suggest area-specific costs that guarantee sufficient incentives for all areas to accept the solution. To this end, we formulate a preemptive model in a framework that allows the application of coalitional game theory methods to obtain a stable benefit allocation, i.e., an outcome immune to coalitional deviations ensuring willingness of all areas to coordinate. We show that benefit allocation mechanisms can be formulated either at the day-ahead or the real-time stages, in order to distribute the expected or the scenario-specific benefits, respectively. For both games, the proposed benefits achieve minimal stability violation, while allowing for a tractable computation with limited queries to the bilevel program. Our case studies, based on an illustrative and a more realistic test case, compare our method with well-studied benefit allocations, namely, the Shapley value and the nucleolus. We show that our method performs better in stability, tractability, and fairness, which would potentially be dictated by a criterion chosen by the regulator.

\end{abstract}
\end{frontmatter}
\newpage

\section{Introduction}

{The existing market architectures and the predominant models for system balancing were designed in a time when fully controllable generators with nonnegligible marginal costs were prevalent.}
{However, as the increasing shares of variable and partly predictable renewables displace controllable generation the dispatch flexibility decreases, while the operational uncertainty characteristics become increasingly complex~\citep{denholm2011grid}.}
{In light of this new operational paradigm, there is an imperative need to re-evaluate the current market design}, and there has been a surge of interest in proposing new market frameworks~\citep{neuhoff2005large,ahlstrom2015evolution}.
\looseness=-1

According to the European electricity market design, the bulk volume of energy trading takes place in the day-ahead market, which is typically cleared 12-36 hours before the actual delivery, based on single-valued point forecasts of the stochastic power output of renewables. In turn, a balancing market is cleared close to real-time operation in order to compensate any deviations from the day-ahead schedule.  Apart from these energy-only trading floors, a reserve capacity auction is organized, usually prior to the day-ahead market, in order to ensure that sufficient capacity is set aside for the provision of balancing services. 
Following this sequential clearing approach, the current market structure attains only limited coordination between day-ahead and balancing stages.
%between the scheduling and the balancing processes. 
Aiming to enhance this temporal coupling, recently proposed dispatch models employ stochastic programming, which allows for a probabilistic description of renewables' forecast errors, in order to co-optimize the day-ahead and the reserve capacity markets \citep{pritchard2010single,morales2012pricing}. However, these approaches cannot be directly applied to the {European} electricity market, since they would require significant restructuring of the current market framework.

In terms of geographical considerations, the European electricity market is fully coordinated only at the day-ahead stage, while reserve capacity and balancing markets are still operated on a regional level~\citep{gebrekiros2015analysis}. However, the European Commission (EC) regulation has already established a detailed guideline framework \citep{EC1} that lays out the rules also for the integration of the European balancing markets in order to improve the security of supply and the efficiency of the balancing system. 
In this ongoing process that is expected to be completed by 2023 \citep{ENTSOE1}, several questions remain open regarding the specific structure and final coordination arrangement among the various ENTSO-E (European Network of Transmission System Operators for Electricity) member countries, since there is no binding legislation that enforces transmission system operators (TSOs) to enter such collaborations. The recent study in \citep{EC2} to investigate the costs and the benefits of different organizational models for the integration of balancing markets, shows that the 10-year net present value (NPV) of coordinated balancing ranges from 1,7 to 3,8 B\euro, depending on the degree of coordination in the inter-area exchanges.
% In addition, it is shown that the implementation of different `market design initiatives' that promote the regional dimensioning and the cross-border exchange of reserves can yield a 10-year NPV up to 27 B\euro.
Nonetheless, any coordination arrangement in the procurement and activation of reserves depends upon the availability of inter-regional transmission capacity. Given that the energy and reserve capacity markets are cleared sequentially, the cross-border transmission resources have to be ex-ante allocated between energy and reserves. 
However, the exact methodology for allocating inter-area transmission capacity remains in consultation with ENTSO-E members \citep{ENTSOE2}. To this end, the goal of this paper is to propose an approach for the design of a transmission allocation mechanism that attains high efficiency both from a technical and an economic perspective.
\looseness=-1
% Even though the European system has recently established an internal day-ahead market, the market lacks spatial coordination with respect to the reserve capacity and balancing trading floors~\citep{gebrekiros2015analysis}.
% If the complete integration of the provision and activation of reserves is to be achieved, the sequential electricity market structure requires that, similar to the flexible generation, the inter-area transmission capacities have to be ex-ante allocated between energy and reserve exchanges. 

The reservation of inter-area interconnections for reserve exchange withdraws transmission resources from the day-ahead market, where the main volume of electricity is being traded. Hence, a sub-optimal reservation may lead to significant efficiency losses at the day-ahead stage. To prevent this situation, the~\citet{ACER} mandates to perform detailed analyses demonstrating that such reservation would increase overall social welfare. Up to this date, inter-area transmission capacity is typically allocated entirely to day-ahead energy exchanges. One notable exemption is the Skagerrak interconnector between Western Denmark and Norway, in which $15\%$ of the capacity is permanently set aside for reserve exchange~\citep{energinet}. Nevertheless, this allocation is static, while the true optimum varies dynamically depending on generation, load, and system uncertainties. As such, \cite{delikaraoglou2018optimal} developed a preemptive transmission allocation model that defines the optimal inter-area transmission capacity allocation to improve both spatial and temporal coordination {at the reserve procurement stage}. 
\looseness=-1
% The proposed model was formulated as a stochastic bilevel programming problem, {which is a powerful mathematical tool that enables the decision-maker, i.e., system/market operator, to anticipate the effect of its decisions on the subsequent market clearings~\citep{morales2014electricity}.}

{
The recently proposed preemptive transmission allocation model focuses on the minimization of the expected system cost, assuming implicitly full coordination among the regional operators. This assumption is in line with the current state of the day-ahead market, which is fully integrated across Europe, or even for the balancing markets in certain regions, e.g., in the Nordic system all reserve activation offers are pooled into a common merit-order list and are available to all TSOs~\citep{bondy2014operational}. However, the initial version of the preemptive model} does not suggest an area-specific cost allocation which guarantees that all areas have sufficient incentives to accept the proposed solution. To address this issue, we {integrate} the preemptive model in a {mathematical} framework that allows the application of tools from coalitional game theory in order to obtain a stable benefit allocation, i.e., an outcome immune to coalitional deviations ensuring that all areas are willing to coordinate via the preemptive model. 
The {concepts} and the tools from coalitional game theory have recently been widely used in the energy community. The Shapley value has been employed in similar problems regarding the distribution of social welfare among TSOs participating in an imbalance netting cooperation~\citep{avramiotis2018investigations} as well as the benefit allocation in transmission network expansion problems~\citep{ruiz2007effective} and in cross-border interconnection development in the North Sea offshore grid~\citep{kristiansen2018mechanism}. 
% Other applications of the Shapley value in the energy field include cooperation problems in the Eurasian gas supply system \citep{nagayama2014network} and the $\text{CO}_2$ emissions abatement in mainland China~\citep{he2018estimation}.} 
However, the Shapley value is in general not within the core, i.e., the set of stable outcomes.
\citet{baeyens2013coalitional} shared the expected profits from aggregating wind power generation using the core benefit allocations, whereas a similar concept was applied for prosumer cooperation in a combined heat and electricity system in~\citep{mitridatidesign}, and for cross-border transmission expansion in the Northeast Asia in~\citep{churkin2019can}.
%~\citep{guajardo2015common} %contreras1999coalition, %kristiansen2018generic
% in allocating the benefits in transmission network expansion problems~\citep{contreras1999coalition,ruiz2007effective}, international cooperation problems in Eurasian gas supply system \citep{nagayama2014network}, $\text{CO}_2$ emissions abatement in mainland China~\citep{he2018estimation} \sdcomm{and in} cross-border interconnection development in North Sea offshore grid~\citep{kristiansen2018mechanism,kristiansen2018generic}.

In contrast to the case studies in the aforementioned works, realistic instances of our problem can potentially exhibit an empty core. To this end, we utilize the least-core as a solution concept, since it achieves minimal stability violation, i.e., minimal incentives for coalitional deviations~\citep{maschler1979geometric}. To obtain a unique outcome, we propose the approximation of a fairness criterion, which {is at the discretion of the regulatory authorities to define.
We propose two variations of the benefit allocation mechanism that can be executed either at the day-ahead or the real-time stage to distribute the expected or the actual benefits (i.e., when the uncertainty is revealed), respectively.} % In addition, depending on when the total benefits are computed, the problem of allocating the benefits in reserve exchanges can be formulated as a coalitional game at the day-ahead or the real-time stages.
We illustrate how these two formulations establish a trade-off between allocating the risk to the regulator or to the regional operators.
% , which arises from the minimization of the expected costs in the preemptive model. 
% Moreover, both the core and the least-core requires the exhaustive enumeration of all possible coalitional deviations, which would grow exponentially in the number of areas.
% and quickly become intractable for the stochastic bilevel program under consideration. 
To overcome the exhaustive enumeration of all coalitional deviations, we show that the least-core selecting benefit allocations in this work can be computed efficiently via an iterative constraint generation algorithm. Similar algorithms were utilized to compute an outcome from the core in combinatorial auctions~\citep{day2007fair}, inventory games~\citep{drechsel2010computing} and electricity markets~\citep{karaca2019designing}. In contrast, our algorithm treats the least-core and it is formulated for a general nonconvex problem.
\looseness=-1

Our contributions are as follows. We formulate the coalition-dependent version of the preemptive transmission allocation model to consider coalitional deviations. We then study the coalitional game that treats the benefits as an ex-ante process with respect to the uncertainty realization {and} we provide a condition under which the core is nonempty. In case this condition is not found, we show that the least-core, which attains minimal stability violation, also ensures the individual rationality property. We then propose the least-core selecting mechanism as a benefit allocation that achieves minimal stability violation, while enabling the approximation of an additional fairness criterion. In order to {implement} this mechanism with only a few queries to the preemptive model, we formulate a constraint generation algorithm. {In addition, we formulate a variation of the} coalitional game that {allocates} the benefits {in} an ex-post process, which can be applied only after {the uncertainty realization is known}. For this game, we provide conditions under which the core is empty {and} we propose an ex-post version of our benefit allocation mechanism.
The ex-ante and ex-post versions of this mechanism can achieve the same fundamental properties for the regional operators either for every uncertainty realization or in expectation, respectively.
Finally, our results are verified first with an illustrative three-area nine-node system and then with a more realistic case study based on a larger IEEE test system. 
\looseness=-1
% Via this ex-post extension, we can achieve the fundamental properties associated with ex-ante benefit allocations in expectation.

The remainder of the paper is organized as follows. Section~\ref{sec:elecmarkframe} describes the organizational structure of electricity markets and introduces a set of necessary assumptions to obtain tractable models. Section~\ref{sec:transcapalloc} discusses the issues related to reserve exchanges and different coordination arrangements for transmission allocations {and motivates the formulation of} the preemptive transmission allocation model. Section~\ref{sec:3} {introduces some necessary background from} coalitional game theory, whereas Section~\ref{sec:4} focuses on the coalitional games arising from the preemptive model, {which provide the basis for the design of the} benefit allocation mechanisms that accomplish the implicit coordination requirements outlined in the previous section.
The numerical case studies are presented in Section~\ref{sec:CaseStudies} and Section~\ref{sec:conc} concludes the paper and gives suggestions for future work.

\section{Electricity market framework}\label{sec:elecmarkframe}

\subsection{Sequential electricity market design and modeling assumptions}\label{sec:seqassm}

The existing market design based on the sequential and independent clearing of reserves, day-ahead, and balancing markets suffers from two main caveats that become increasingly pronounced as we move towards larger shares of renewable energy production. On the one hand, the day-ahead schedule is optimized based on purely deterministic inputs, i.e., single-valued forecasts of renewable production.
As a result, the day-ahead market is not responsive to the uncertainty associated with the forecast errors of stochastic renewables and thus it is weakly coordinated with the real-time balancing actions. On the other hand, the decoupling of energy and upward/downward reserve capacity trading into two independent auctions ignores the substitution and complementary properties of these two services and leads to inefficient reserve procurement and energy schedules. In order to enable the inter-area exchange of reserves given this decoupling of reserve and energy markets, the operator has to prescribe a certain share of the interconnection capacities to the reserve capacity market and in turn withdraw this headroom from the day-ahead energy trading.
\looseness=-1

From a theoretical perspective, these two issues can be contained if reserve capacity procurement, day-ahead energy schedules, and real-time re-dispatch actions are jointly optimized based on a probabilistic description of the uncertainty. As several recent works show, e.g., \citet{bouffard2005market}, \citet{pritchard2010single}, and \citet{morales2012pricing} among others, the co-optimization of these services using a stochastic programming framework reduces the total expected system cost and it establishes a benchmark in terms of perfect temporal coordination. However, the adoption of a stochastic dispatch model as a market-clearing algorithm requires significant restructuring of the current European market framework that electricity sector stakeholders may not be willing to embrace.
Owing to this reason, we restrict ourselves to the status-quo market architecture and we embody its design attributes in our methodology aiming to mitigate the resulting inefficiencies.

In the following, we build the mathematical models of the different trading floors based on a set of realistic assumptions that allows us to capture the main attributes of the European market, while maintaining computational tractability. 
In line with the current practice, we consider a zonal network representation during the reserve procurement, whereas the full network topology is taken into account in the day-ahead and the balancing market-clearing models, using a DC power flow approximation. Note that our network model can be readily adapted to a zonal pricing scheme, where the inter-zonal transmission energy flows are constrained by the available transfer capacity, following the internal electricity market paradigm. However, since this work focuses on transmission allocation issues concerning primarily the operators, we believe that allowing for a more complete network representation will increase the value of the proposed models as decision-support~tools.

On the generation side, we consider that all market participants are perfectly competitive. Day-ahead energy offers are submitted in price-quantity pairs that internalize the marginal production cost as well as the unit commitment and inter-temporal constraints, e.g. ramping limits, in accordance with the portfolio bidding practice in the European market.
Moreover, we assume that the reserve capacity offer prices provide adequate incentives to the flexible generators for the provision of real-time balancing services, i.e., no price premiums are required in balancing to compensate for the opportunity cost of the capacity withdrawn from the day-ahead market.
In terms of stochastic renewable in-feed, we focus on wind power generation and we model forecast errors using a finite set of scenarios. Assuming null production costs, the corresponding offer price and the spillage cost are set equal to zero. 
% Wind power forecast errors are modeled using a finite set of scenarios, while other sources of uncertainty such as demand
% variations can be considered in a similar way. 
On the consumption side, we consider inelastic demand with a large penalty on lost load and thus the social welfare maximization becomes equivalent to the cost minimization.
\looseness=-1

Finally, we assume that the current implementation of the sequential market provides a budget balanced method to allocate the system cost to all the areas, i.e. the costs of the reserve capacity market, the day-ahead market and the balancing market are allocated to the areas without any deficit or surplus. In the numerics, for the case of no inter-area exchange of reserves, we provide and discuss one such allocation method based on the zonal and nodal prices that assigns producer and consumer surpluses to their corresponding areas, and divides the congestion rent equally between the adjacent areas, see~\citep{kristiansen2018mechanism} and Section~\ref{sec:illu_ex}.

\subsection{Mathematical formulation}
\label{sec:seq}
%Following the sequential market structure described above, we provide below the mathematical formulation of the reserve capacity, the day-ahead and the balancing market-clearing problems.
For the models below, main notation is stated in~\ref{app:not}.

\subsubsection{Reserve capacity market} Having as a fixed input the upward/downward reserve requirements $RR_a^{+}/RR_a^{-}$ in each area~$a$ and a pre-defined share $\chi_e$ of the transmission capacity of each inter-area link~$e$ allocated to reserves exchange, the reserve capacity market clearing is formulated as:
\begin{subequations}
	\label{mod:A-Res}
	\begin{flalign}
		 &\underset{\Phi_{\text R}}{\text{minimize}} \quad  \sum_{i \in \mathcal I}\left(C^{+}_{i} r^{+}_{i} + C^{-}_i r^{-}_{i}\right) \label{eq:A-Res-obj}	 % we need to assume C^{+}\leq C \leq C^{-}. Otherwise, this can easily be gamed.
	\end{flalign}\vspace{-1cm}
	\begin{flalign}\vspace{-1cm}
	    \text{subject\ to}\ & \sum_{i \in \mathcal{M}_{a}^{\mathcal I}} r_{i}^{+} + \sum_{e\in\mathcal E} \mathcal{H}(e,a) r_{e}^{+} \geq RR_a^{+}, \quad \forall a\in\AC, \label{eq:A-Res-req-up}\\
	    & \sum_{i \in \mathcal{M}_{a}^{\mathcal I}} r_{i}^{-} + \sum_{e\in\mathcal E} \mathcal{H}(e,a) r_{e}^{-} \geq RR_a^{-}, \quad \forall a\in\AC,  \label{eq:A-Res-req-dn}\\
		& 0 \le r_{i}^{+} \le  {R}^{+}_{i}, \quad \forall i\in \mathcal I, \quad 0 \le r_{i}^{-} \le  {R}^{-}_{i}, \quad \forall i\in\mathcal I, \label{eq:A-Res-cdn-max} \\% 
		& -\chi_{e} T_{e}\le r_{e}^{+} \le  \chi_{e} T_{e}, \quad \forall e \in \mathcal E, \quad
		-\chi_{e} T_{e} \le r_{e}^{-} \le \chi_{e} T_{e}, \quad \forall e \in \mathcal E, \label{eq:A-Res-X2} 
		%& \textcolor{WildStrawberry}{\chi_{e}=0,\ \forall e\ \mathrm{s.t.}\ \exists a\in \mathcal{C}\setminus \mathcal{S}\ \mathrm{with}\ \mathcal{H}(e,a)\neq0},
	\end{flalign}
\end{subequations}
where $\Phi_{\text R} = \{ r^{+}_{i},r^{-}_{i}, \forall i ; r_{e}^{+}, r_{e}^{-}, \forall e  \} $ is the set of optimization variables. The objective function \eqref{eq:A-Res-obj} to be minimized is the cost of reserve procurement. 
Constraints \eqref{eq:A-Res-req-up} and \eqref{eq:A-Res-req-dn} ensure, respectively, that the upward and downward reserve requirements of each area are satisfied either by procuring reserve capacity from intra-area generators or via inter-area reserves exchange that is modelled using the incidence matrix $\mathcal{H}(e,a)$. 
As shown in~\ref{app:not}, for each link $e$ with sending and receiving ends in areas $a_s(e)$ and $a_r(e)$, respectively,
$\mathcal{H}(e,a)$ is equal to 1 (-1) if reserve import (export) is considered from (to) area $a=a_s(e)$ ($a=a_r(e)$) and zero for any other area. With this definition availability of cross-border reserves within the neighboring areas is guaranteed by \eqref{eq:A-Res-req-up} and \eqref{eq:A-Res-req-dn}. We underline that directed links are used as a notational convention, and both $r_e^+$ and $r_e^-$ are free of sign.
Upward and downward capacity offers of dispatchable power plants are enforced by constraints \eqref{eq:A-Res-cdn-max}. %Note that for inflexible power plants these capacities are given by zero.
The set of constraints \eqref{eq:A-Res-X2} models the bounds on reserves exchange between two areas across link $e$. Following the current practice, we consider a zonal network representation for the reserve capacity markets and thus the transmission capacity ${T}_e$ of link $e$ is defined as the aggregated flow limit of all tie-lines $\ell \in \Lambda_{e}$ across link~$e$ calculated as
	${T}_e = \sum_{\ell \in \Lambda_{e}} \ {T}_{\ell},$ for all $e\in\mathcal E$.

Setting the transmission capacity allocation $\chi_e$ to any value different than zero, establishes practically a reserve exchange mechanism between the areas located at the two ends of the link and consequently it enables the exchange of balancing services during real-time operation. On the contrary, setting $\chi_e=0$ implies that there would be no reserve exchange at the procurement stage, i.e., the cross-border transmission capacity is fully allocated to day-ahead energy exchanges. In that case, we also prevent the exchange of balancing services and the imbalance netting between the adjacent areas, as we will formally describe in the balancing market model formulation below.

\subsubsection{Day-ahead market}

Given the optimal reserve procurement $\hat \Phi_{\text R} = \{ \hat r^{+}_{i},\hat r^{-}_{i}, \forall i ; \hat r_{e}^{+}, \hat r_{e}^{-}, \forall e  \}$, the day-ahead schedule is the solution to the following optimization problem:  
\begin{subequations}
	\label{mod:A-DA}
	\begin{flalign}
		&\underset{\Phi_{\text D}}{\text{minimize}} \quad  \sum_{i\in \mathcal I}{C}_{i}{p}_{i} \label{eq:DA-Obj}
	\end{flalign}\vspace{-1cm}
	\begin{flalign}	\vspace{-1cm}
		\text{subject\ to}\ & \sum_{j \in \mathcal{M}_n^{\mathcal J}} w_j + \sum_{i\in \mathcal{M}_n^{\mathcal I}}{p}_{i} - \sum_{\ell \in \mathcal L^{\text{AC}}} A_{\ell n} f_{\ell}   = D_n, \quad \forall n\in \mathcal N,  \label{eq:A-CV-da-bal-constr} \\
		&  \hat{r}_{i}^{-} \le p_i \leq \ {P}_{i}- \hat{r}_{i}^{+},   \quad \forall i\in\mathcal I,  \label{eq:A-CV-da-conv-max}\\
		%& \sum_{a} R_{ia}^{-} \le p_i,  \quad \forall i  \label{eq:A-CV-da-con-min}\\
		& 0 \le w_j \leq \overline{{W}}_{j},  \quad \forall j\in\mathcal J,  \label{eq:A-CV-da-wind-max}\\
		& f_{\ell} = B_{\ell} \sum_{n\in\mathcal N} A_{\ell n} \delta_n, \quad \forall \ell \in \mathcal L^{\text{AC}}, \label{eq:A-CV-da-flowAC-constr}\\
		& -(1- \chi_{\ell}) \ {T}_{\ell} \le f_{\ell} \le (1- \chi_{\ell}) \ {T}_{\ell}, \quad \forall \ell \in \mathcal L^{\text{AC}} , \label{eq:A-CV-da-flowAC-max}\\			
		& \delta_1 = 0, \quad \delta_n \; \text{free}, \quad \forall n \in \mathcal{N}, \label{eq:A-CV-da-RefBus-constr}
	\end{flalign}
\end{subequations}
where $\Phi_{\text D} = \{p_i, \forall i; w_{j}, \forall j; \delta_n, \forall n; f_{\ell}, \forall \ell \}$ is the set of optimization variables. We define $\chi_\ell= \chi_e$ for all tie-lines $\ell \in\Lambda_e$ and $ \chi_\ell=0$ for all intra-area lines. For the remainder, we strictly follow this notation. The objective function \eqref{eq:DA-Obj} to be minimized is the day-ahead cost of energy production. Constraints \eqref{eq:A-CV-da-bal-constr} enforce the day-ahead power balance for each node. The upper and lower production limits of dispatchable power plants are enforced by~\eqref{eq:A-CV-da-conv-max}, taking into account the reserve schedule from the previous trading floor. Constraints \eqref{eq:A-CV-da-wind-max} limit the stochastic production to a point forecast, typically the expected value of the stochastic process. Power flows are computed in \eqref{eq:A-CV-da-flowAC-constr} and then restricted by the capacity limits in \eqref{eq:A-CV-da-flowAC-max} considering that only $(1- \chi_{\ell})$ percent of the capacity is available for day-ahead energy trade. The voltage angle at node~$1$ is fixed to zero in \eqref{eq:A-CV-da-RefBus-constr} setting this as the reference node, whereas the remaining voltage angles are declared as free~variables. 
\looseness=-1
%Let the optimal solution be denoted by $\hat \Phi_{\text D} = \{\hat p_i, \forall i; \hat w_{j}, \forall j; \hat \delta_n, \forall n; \hat f_{\ell}, \forall \ell \}$.

\subsubsection{Balancing market}

Being close to real-time operation uncertainty realization $s'$ and actual wind power production $W_{js'}, \; \forall j \in \mathcal{J}$ are known. Any energy deviations from the optimal day-ahead schedule $\hat \Phi_{\text D} = \{\hat p_i, \forall i;\allowbreak \hat w_{j}, \forall j; \hat \delta_n, \forall n;\allowbreak \hat f_{\ell},\allowbreak \forall \ell \}$ must be contained using proper re-dispatch actions that respect the reserve procurement schedule $\hat \Phi_{\text R}$.
To determine the re-dispatch actions that minimize the balancing cost, the balancing market is cleared based on the following optimization problem: 
\begin{subequations}
	\label{mod:A-CV-rt}
	\begin{flalign}
		&\underset{ \Phi_{\text B}^{s'} }{\text{minimize}} \quad  \sum_{i \in \mathcal I}{C}_{i}\left({p}_{is'}^{+}-{p}_{is'}^{-}\right)+\sum_{n \in \mathcal N} {C}^{\text{sh}}l^{\text{sh}}_{ns'} \label{eq:A-CV-rt-Obj}
	\end{flalign}\vspace{-.1cm}
	{\medmuskip=.65mu\thinmuskip=.65mu\thickmuskip=.65mu\begin{flalign}\vspace{-.1cm}
		\text{subject to}\ & \sum_{i \in \mathcal{M}_n^{\mathcal I}} \left(p_{is'}^{+}-p_{is'}^{-} \right) + l^{\text{sh}}_{n s'} + \sum_{j \in \mathcal{M}_n^{\mathcal J}}\left(W_{j s'}-\hat{w}_j-w^{\text{spill}}_{j s'}\right)
		+ \sum_{\ell \in \mathcal{L}^{\text{AC}}} A_{\ell n} \left( \hat{f}_{\ell} - {f}_{\ell s'} \right) = 0, \ \forall n\in\mathcal N,   \label{eq:A-CV-rt-bal-constr} \\	
		& 0 \le p^+_{is'} \leq  \hat{r}_{i}^{+}, \quad \forall i\in\mathcal I, \quad 0 \le p^-_{is'} \leq  \hat{r}_{i}^{-} , \quad \forall i\in\mathcal I, \label{eq:A-CV-rt-conv-res-dn-max}\\
		& 0 \le l^{\text{sh}}_{n s'} \leq D_n, \quad \forall n\in\mathcal N,   \quad 0 \le w^{\text{spill}}_{j s'} \leq W_{j s'}, \quad \forall j\in\mathcal J,  \label{eq:A-CV-rt-wind-spill-max}\\
		& {f}_{\ell s'} = B_{\ell} \sum_{n\in\mathcal N} A_{\ell n} {\delta}_{n s'}, \quad \forall \ell \in \mathcal L^{\text{AC}},  \label{eq:A-CV-rt-flowAC-constr}\\
		& - T_{\ell} \le {f}_{\ell s'} \le T_{\ell}, \quad \forall \ell \in \mathcal L^{\text{AC}},  \label{eq:A-CV-rt-flowAC-max}\\
		& f_{\ell s'} = \hat f_\ell,\quad\forall \ell \in \cup_{e\in\mathcal E^-(\chi)} \Lambda_e, \label{eq:A-CV-rt-tieline}\\
		& {\delta}_{1 s'} = 0, \quad \delta_{n s'} \; \text{free}, \quad \forall n \in \mathcal{N}, \label{eq:A-CV-ref-node}
	\end{flalign}}
\end{subequations}
where $\Phi_{\text B}^{s'}\allowbreak =\allowbreak \{ {p}^{+}_{i s'},\allowbreak{p}^{-}_{is'}, \allowbreak\forall i; \allowbreak w^{\text{spill}}_{js'},\allowbreak \forall j; {l}^{\text{sh}}_{ns'},\allowbreak {\delta}_{ns'}, \forall n; {f}_{\ell s'},\forall \ell \}$ is the set of optimization variables.

 The objective function~\eqref{eq:A-CV-rt-Obj} is the cost of re-dispatch actions, i.e., reserve activation and load shedding. Up and down re-dispatches have the same prices since we previously assumed that no price premiums are required in the balancing stage to compensate for the opportunity cost of withdrawing capacity from the day-ahead stage, see Section~\ref{sec:seqassm}. Equality constraints \eqref{eq:A-CV-rt-bal-constr} ensure that all the nodes remain in balance after the re-dispatch of generation and any necessary wind power curtailment or load shedding. Constraints \eqref{eq:A-CV-rt-conv-res-dn-max} ensure that upward and downward reserve deployment respects the corresponding procured quantities. The upper bounds on load shedding and power spillage are set equal to the nodal demand and the realized wind power production by constraints \eqref{eq:A-CV-rt-wind-spill-max}. Real-time  power flows are modeled by \eqref{eq:A-CV-rt-flowAC-constr} and then restricted by the transmission capacity limits in \eqref{eq:A-CV-rt-flowAC-max}. Constraints \eqref{eq:A-CV-rt-tieline}, where $\mathcal E^-(\chi)=\{e\in\mathcal E \rvert  \chi_e=0\}$ denotes the set of inter-area links with no existing cross-border agreement across them, ensure that if $\chi_e=0$,  the real-time flows on the tie lines are fixed to their day-ahead values.
 This follows from the previous assumption that prevents any reserve sharing or imbalance netting during real-time operation across link $e$ if $\chi_e=0$.
Similarly to the day-ahead formulation, node~$1$ is set as the reference node in \eqref{eq:A-CV-ref-node} and the remaining voltage angles are declared as free variables. 
\looseness=-1
% Let $\hat \Phi_{\text B}^{s}\allowbreak =\allowbreak \{ \hat {p}^{+}_{is},\allowbreak\hat {p}^{-}_{is},\allowbreak \forall i;\allowbreak \hat w^{\text{spill}}_{js},\allowbreak \forall j;\allowbreak \hat {l}^{\text{sh}}_{ns}, \hat {\delta}_{ns}, \forall n; \hat {f}_{\ell s},\forall \ell \}$ denote the optimal solution. 

\section{Transmission capacity allocations for cross-border balancing}\label{sec:transcapalloc}

\subsection{Coordination schemes and transmission allocation arrangements for cross-border trading}

The transition to an integrated balancing market requires several organizational changes to the prevailing operational model, in which reserves are procured and deployed on an intra-area basis.
A prerequisite for the establishment of a well-functioning balancing framework is the standardization of the rules and
products as well as the definition of transparent mechanisms that will facilitate the cooperation among the TSOs \citep{hobbs2005more}. Below, we outline the main coordination schemes and transmission allocation arrangements as defined in the current European regulation~\citep{EC1}.

Inter-area reserve procurement can be organized as a reserve exchange scheme {and}/{or} as a reserve sharing agreement. Implementing the former scheme, regional TSOs can procure balancing capacity resources located in adjacent areas in order to meet their own area reserve requirements. 
Since the reserve requirements of each area remain unchanged, this coordination setup requires limited organizational changes, as it basically reallocates the reserve quantities towards areas with lower procurement costs. 
To improve also the dimensioning efficiency of the procurement process, a reserve sharing agreement allows a TSO to use available reserve capacity from adjacent TSOs. An implied prerequisite for this arrangement would be that the definition of regional reserve requirements is performed jointly by all TSOs that participate in the sharing agreement. 

In terms of coordination during reserves activation, the main organizational setups are the so-called imbalance netting and exchange of balancing energy. The first setup pertains to the inter-area exchange of imbalances with opposite sign, thus preventing the counteracting activation of balancing resources and reducing the total balancing energy volumes. 
In turn, the exchange of balancing energy enables the system-wide least-cost activation of reserves through a common merit-order list to meet the net imbalance of the joint TSO area. This improves the supply efficiency of balancing energy, at the expense of more extensive coordination requirements.

\begin{figure}[h]
	\centering
	\begin{tikzpicture}[scale=.84, every node/.style={scale=0.6}]
			\draw (4,1) circle (.35cm) node {\LARGE $a_1$};
			\draw (4,-1) circle (.35cm) node {\LARGE $a_2$};
			\draw[-,line width=.25mm] (4,.65) -- (4,-.65) node  at (4.4, -0.1) {\LARGE $\chi$};
			\node at (-.85,0) {\includegraphics[width=.6\textwidth, trim={3.4cm 9.7cm 3.9cm 9.8cm},clip]{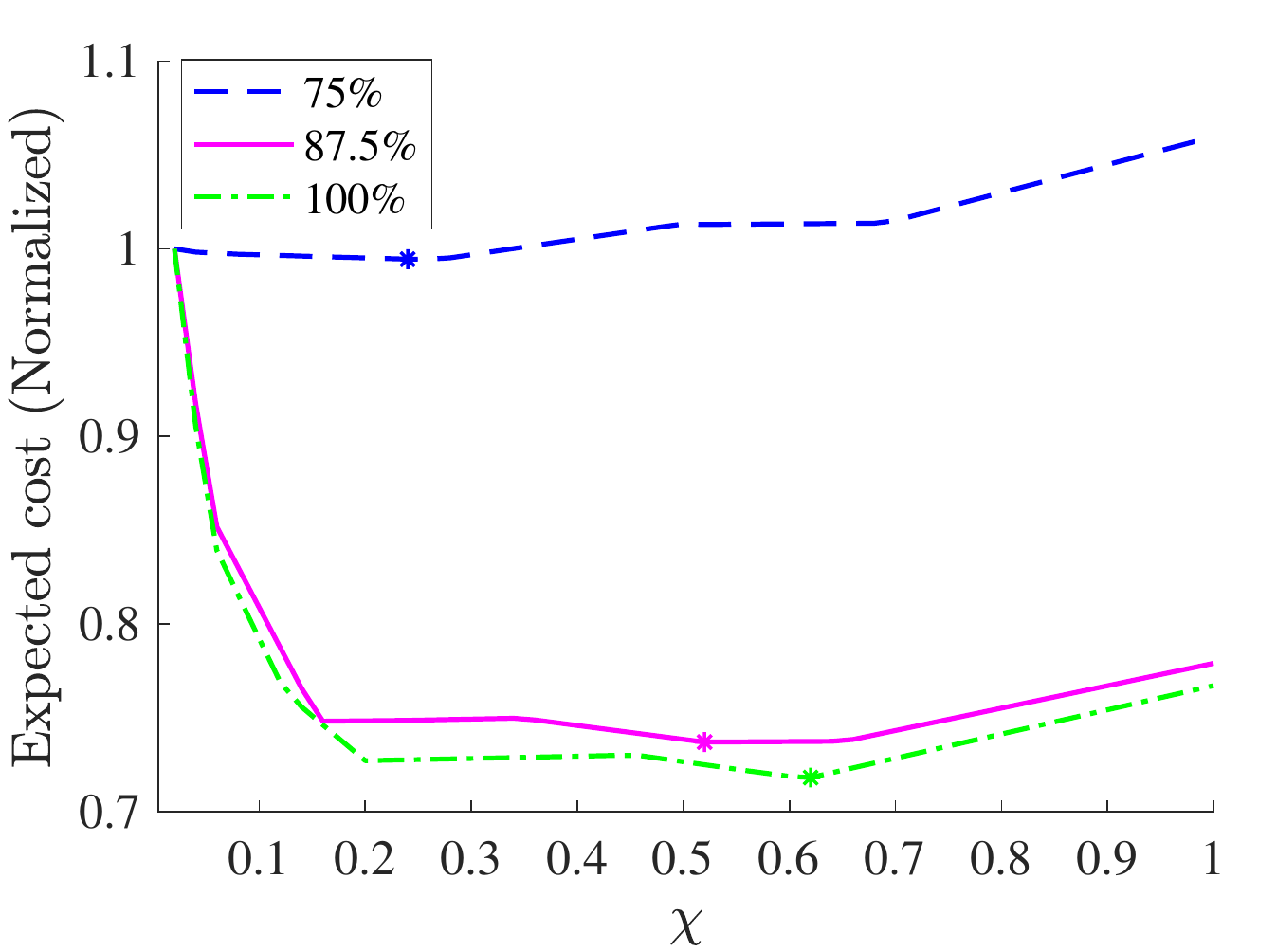}};
			\fill[fill=white] (-2.93,2.25) rectangle (-2.07,1.4);
			\node at (-2.45,1.58) {\footnotesize $312\; \text{MW}$};
			\node at (-2.45,1.87) {\footnotesize $273\; \text{MW}$};
			\node at (-2.45,2.15) {\footnotesize $234\; \text{MW}$};
			%\draw[->] (-2,-1.5) -- (2,-1.5) node[right] {\LARGE $\chi$};
			%\draw[->] (-2,-1.5) -- (-2,2) node at (-1.5, 2.2) {\large $\text{Expected cost of \eqref{mod:A-CV-rt}, \eqref{mod:A-DA}, and \eqref{mod:A-Res}}$};
			%\draw[scale=0.5,domain=-3:3,smooth,variable=\x,blue] plot ({\x},{\x*\x});
			%\draw[scale=0.5,domain=-3:3,smooth,variable=\y,red]  plot ({\y*\y},{\y});
	\end{tikzpicture}
	\caption{Expected operation cost (in \euro) as a function of transmission capacity allocated to inter-area reserves trading under different levels of wind power penetration (in MW).}\label{fig:plotta}
\end{figure}
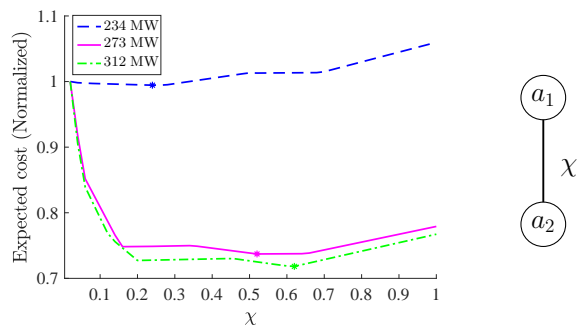

The establishment of any cross-border reserve procurement scheme requires the reservation of a certain share of the inter-area transmission capacity for the reserves and their activation. Before we describe the attributes of any specific transmission allocation mechanism, let us provide an illustrative example. This example highlights the seams issues pertaining to the ex-ante definition of transmission allocation between neighboring areas.
Figure~\ref{fig:plotta} shows the expected system cost, i.e., the sum of reserve procurement, day-ahead energy and expected balancing costs, as a function of the share of transmission capacity $\chi$ that is allocated to inter-area reserves trading. The data for this two-area power system is provided in~\ref{app:two_area_ex}. 
Even from this simple example, it becomes apparent that there generally exists an optimal allocation to be made, which however may dynamically vary depending on generation, load and system uncertainties.
Therefore, the efficiency of an integrated balancing market, in terms of expected system cost, is highly susceptible to the parameter~$\chi$. This in turn asks for a systematic method to optimally define~$\chi$, accounting for the market dynamics and the uncertainty involved in the operation of the power system. %We refer the interested reader to~\citep[\S 2.4]{delikaraoglou2018optimal} for a comprehensive discussion on the other open issues and potential caveats related to the existing transmission allocation models.

In this work, we focus on the prevailing market-based mechanism for the allocation of cross-border transmission capacity between energy and reserves. According to this methodology,
a share of inter-area transmission capacity is set aside for reserves based on the comparison of the market value of cross-zonal capacity for the exchange of balancing capacity or sharing of reserves and the market value of cross-zonal capacity for the exchange of energy. This methodology provides a reasonable trade-off between allocation efficiency and practical applicability within the current framework, albeit it still incurs the inherent drawbacks of the sequential market structure regarding the deterministic view of uncertainty and the separation of energy and reserve services. We refer the interested reader to \citep{EC1} for discussions on alternative transmission allocation mechanisms.
%Next, we bring in the preemptive tranmission allocation model.

\subsection{Preemptive transmission allocation model}
In this section, we first describe the preemptive transmission allocation model that was initially proposed in~\citep{delikaraoglou2018optimal} and is the key building block of the benefit allocation mechanisms proposed in this work. This work, however, defines the preemptive model in a more general framework, which allows us to consider coalitional deviations and in turn define the necessary incentives that support the solution proposed by the preemptive model. 

\begin{figure}[h]
	\centering
	\begin{tikzpicture}[scale=0.8, every node/.style={scale=0.56}]
	\draw[dashed, draw=black] (-5,-2.5) rectangle (-1.35,0.95);
	\draw[draw=black,rounded corners] (-4.8,-2.3) rectangle (-1.55,0.75) node at (-3.15,-0.4) {\Large Preemptive model} node at (-3.15,-1) {\Large Problem~\eqref{mod:B-Pre}};
	\draw[draw=black,rounded corners] (.7,-0.7) rectangle (3.95,0.75) node at (2.35,0.4) {\Large Res. cap. market} node at (2.35,-0.2) {\Large Problem~\eqref{mod:A-Res}};
	\draw[draw=black,rounded corners] (6.1,-1.3) rectangle (9.35,0.2)
	node at (7.75,-0.2) {\Large Day-ahead market} node at (7.75,-.8) {\Large Problem~\eqref{mod:A-DA}};
	\draw[draw=black,rounded corners] (11.7,-2.3) rectangle (14.95,0.75) node at (13.35,-0.1) {\Large Balancing market} node at (13.35,-0.65) {\Large Problem~\eqref{mod:A-CV-rt}} node at (13.35,-1.15) {\Large for scenario $s'$};
	\draw[->,line width=.1mm] (-1.55,0) -- (.7,0) node at (-.4,-.5) {\LARGE $\{\hat\chi_e, \forall e\}$};
	\draw[->,line width=.1mm] (3.95,-0) -- (6.1,-0) node at (5.1,-.5) {\LARGE $\left(\begin{subarray}{c} \hat{r}^{+}_{i},\hat{r}^{-}_{i}, \\  \hat{r}_{e}^{+}, \hat{r}_{e}^{-} \end{subarray} \right)$};
	\draw[->,line width=.1mm] (9.35,-0.75) -- (11.7,-.75) node at (10.7,-1.25) {\LARGE $\left(\begin{subarray}{c} \hat{p}_{i}, \hat{w}_{j}, \\  \hat{\delta}_n, \hat{f}_{\ell} \end{subarray} \right)$};
	\draw[->,line width=.1mm] (-1.55,-1.2) -- (6.1,-1.2) node at (2.75,-1.6) {\LARGE $\{1-\hat\chi_e, \forall e\}$};
	\draw[->,line width=.1mm] (3.95,.5) -- (11.7,.5) node at (10.7,0) {\LARGE $\left(\begin{subarray}{c} \hat{r}^{+}_{i},\hat{r}^{-}_{i}, \\  \hat{r}_{e}^{+}, \hat{r}_{e}^{-} \end{subarray} \right)$};
	\draw[->,line width=.1mm] (-1.55,-2) -- (11.7,-2) node at (10.4,-2.35) {\Large $\CC\subseteq\AC$};
	\end{tikzpicture}\vspace{-.2cm}
	\caption{Schematic representation of preemptive transmission allocation model}\label{fig:struct}
\end{figure}
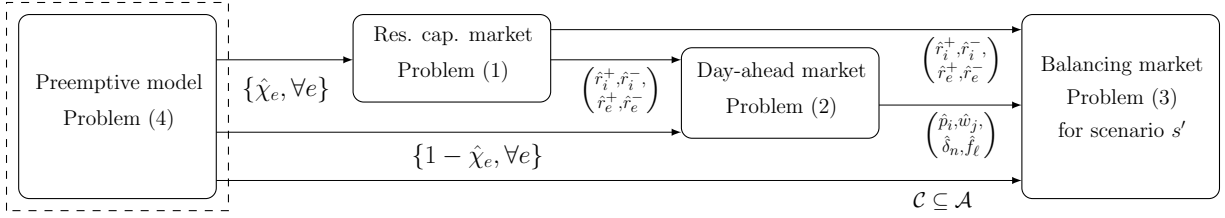

The preemptive transmission allocation model can be perceived as a decision-support tool for the system operators, which aims at defining the optimal shares of transmission capacity for inter-area trading of energy and reserves. Being fully aligned with the existing sequential market structure, the preemptive model is essentially a market-based allocation process that is performed prior to the reserve capacity and day-ahead energy markets to find the optimal transmission allocations $\{\hat\chi_e, \forall e\}$ that minimize the expected system cost, see Figure~\ref{fig:struct}.
% \footnote{This version of the preemptive model implements a reserve exchange scheme, since the reserve requirements $\{RR_a^+,RR_a^-,\forall a\}$ remain unchanged compared to the sequential model. We can easily extend this tool and the related discussions also to a reserve sharing scheme by incorporating the central optimization of the reserve requirements.} 
It is worth mentioning that the coordinated reserve exchange would require a transfer of some responsibilities to European bodies, even though some TSOs might be hesitant to assign some of their autonomy to a central authority.
\looseness=-1

In this work, the primary focus is on the establishment of coalitional agreements for the reserve procurement stage among a set of areas $\CC\subseteq\AC$. In the balancing market, we assume that all areas that participate in the coalition, exchange balancing energy in a perfectly coordinated setup.
 As a result, real-time tie-line flows among these areas are treated as free variables, allowing for deviations from the day-ahead schedule. %Notice that this is the case even when the optimal solution of the preemptive model allocates zero capacity to reserves exchange in the reserve capacity market. 
 The coalition-dependent version of the preemptive model is given by:
\begingroup
\allowdisplaybreaks
\begin{subequations} \label{mod:B-Pre}
{\medmuskip=.58mu\thinmuskip=.58mu\thickmuskip=.58mu\begin{align}
 {J(\CC)}&=\,\underset{ \Phi_{\text{PR}} }{\text{minimize}}\quad \sum_{i \in \mathcal I}\left(C^{+}_{i} r^{+}_{i} + C^{-}_i r^{-}_{i}\right) + \sum_{i\in \mathcal I}{C}_{i}{p}_{i} + \sum_{s\in\mathcal S} \pi_s \Big[ \sum_{i \in \mathcal I}{C}_{i}\left({p}_{is}^{+}-{p}_{is}^{-}\right)+\sum_{n \in \mathcal N} {C}^{\text{sh}}l^{\text{sh}}_{ns} \Big]  \label{eq:B-pre-obj}
\end{align}}\vspace{-1cm}
\begin{align}\vspace{-1cm}
\text{subject\ to}\ & \;\;\;  0\leq\chi_{e}'\leq 1,\quad \forall e\in\mathcal E_{\CC},\label{eq:B-pre-ta}\\
& \;\;\;  \chi_{e}'=\chi_{e},\quad \forall e\in\mathcal E\setminus\mathcal E_{\CC},\label{eq:B-pre-ta-fixed}\\
& \;\;\;  \text{Constraints} \;\; \eqref{eq:A-CV-rt-bal-constr} -\eqref{eq:A-CV-rt-flowAC-max}\text{ and }\eqref{eq:A-CV-ref-node},\quad \forall s\in \mathcal S,\label{eq:B-A-CV}\\
& \;\;\;  {{f}_{\ell s} = {f}_{\ell},\quad \forall \ell\in\cup_{e\in\mathcal E^-(\chi,\CC)} \Lambda_{e},\quad\forall s\in \mathcal S},\label{eq:B-A-CV-flow} \\
& \;\; \left(\begin{subarray}{c} r^{+}_{i},r^{-}_{i}, \\  r_{e}^{+}, r_{e}^{-} \end{subarray} \right)  \in \text{arg}
	\left\{\!\begin{aligned}
	&  \underset{\Phi_{\text{R}}}{\text{minimize}} \quad \eqref{eq:A-Res-obj} \quad \text{subject to}\ \text{constraints \eqref{eq:A-Res-req-up} - \eqref{eq:A-Res-X2}}
	\end{aligned}\right\}, \label{eq:LLR} \\
	& \;\;\; \left(\begin{subarray}{c} p_{i}, w_{j}, \\  \delta_n, f_{\ell} \end{subarray} \right)  \in \text{arg} \;
	\left\{\!\begin{aligned}
	&  \underset{\Phi_{\text{D}}}{\text{minimize}} \quad \eqref{eq:DA-Obj} \quad \text{subject to}\
	\text{constraints \eqref{eq:A-CV-da-bal-constr} - \eqref{eq:A-CV-da-RefBus-constr}}
	\end{aligned}\right\}, \label{eq:LLD} 
\end{align}
\end{subequations}
\endgroup
where $\Phi_{\text{PR}}\allowbreak=\allowbreak\{\chi_e',\allowbreak \forall e \cup \Phi_{\text R} \cup \Phi_{\text D} \cup  \Phi_{\text B}^{s}, \allowbreak\forall s\}$ is the set of primal optimization variables. For the sake of {brevity}, the Lagrange multipliers of the lower-level optimization problems are omitted here, but the complete set of Karush-Kuhn-Tucker (KKT) conditions of problems \eqref{mod:A-Res} and \eqref{mod:A-DA} are listed in~\ref{app:a}. Unless stated otherwise, the preemptive model refers to problem~\eqref{mod:B-Pre} associated with $J(\AC)$ and the corresponding optimal transmission allocation is denoted by $\{\hat\chi_e, \forall e\}$.

Model \eqref{mod:B-Pre} is a stochastic bilevel optimization problem, since wind power uncertainty is described by a finite set of scenarios~$\mathcal{S}$, with $W_{js}$ being the realization of stochastic generation of wind farm $j$ in scenario $s$ and $\pi_s$ the corresponding probability of occurrence. The objective function~\eqref{eq:B-pre-obj} to be minimized is the expected system cost according to the sequential market structure described in Section \ref{sec:seq}.
Constraint~\eqref{eq:B-pre-ta} bounds the share of transmission capacity ${\chi}'_{\textcolor{blue}{e}}\in[0,\,1]^{\mathcal E}$ allocated to reserve exchange for links $e\in\mathcal E_{\CC}$, where $\mathcal E_{\CC}\allowbreak=\allowbreak\{e\,|\, \mathcal{H}(e,a)=0,\,\forall a\in \AC\setminus\CC\}$ is the set of links among only the areas in the coalition $\CC$. The transmission capacity of the remaining links $e\in\mathcal E\setminus\mathcal E_{\CC}$, i.e., links that are connected to areas that are not members of coalition $\mathcal{C}$, is fixed according to existing cross-border agreements ${\chi_e}\in[0,\,1]^{\mathcal E}$ in~\eqref{eq:B-pre-ta-fixed}. 
Constraints~\eqref{eq:B-A-CV} ensure feasibility of re-dispatch actions for each scenario, whereas constraint \eqref{eq:B-A-CV-flow} restricts the real-time tie-line flows to the respective day-ahead values, for links $e \in \mathcal E^-(\chi,\CC)=\{e\in\mathcal E \rvert  \chi_e=0,\ \text{and}\ \exists a\in\AC\setminus\CC\ \text{such}\ \text{that}\ \mathcal{H}(e,a)\neq0\}$, which is the set of links that are, at least on one end, connected to an area that is not in the coalition~$\CC$ and does not have an existing cross-border agreements for reserves exchange, i.e., $\chi_e=0$. In other words, constraint~\eqref{eq:B-A-CV-flow} prevents any imbalance netting or exchange of balancing energy between areas that are not members of coalition $\CC$. Here, notice that constraints~\eqref{eq:B-A-CV-flow} are removed from problem~\eqref{mod:B-Pre} associated with $J(\AC)$, since $\mathcal E(\chi,\AC)=\emptyset$ for any~$\chi$. 

The lower-level problems \eqref{eq:LLR} and \eqref{eq:LLD} are identical to models \eqref{mod:A-Res} and \eqref{mod:A-DA} implementing the shares of transmission capacities in~${\chi}'$. Having this bilevel model \eqref{mod:B-Pre} ensures by construction that the reserve capacity, day-ahead and balancing markets are cleared in consecutive and independent auctions, in accordance with the sequential market design. This structure allows the definition of $\{\hat\chi_e, \forall e\}$ anticipating the impact of these parameters in all subsequent trading floors. From a computational perspective, to obtain a solvable instance of the model \eqref{mod:B-Pre}, we can equivalently replace the lower-level problems \eqref{eq:LLR} and \eqref{eq:LLD} by the corresponding KKT conditions, given that \eqref{eq:LLR} and \eqref{eq:LLD} are linear programs. 
The resulting problem is a single-level mathematical program with equilibrium constraints (MPEC) that involves the non-convex complementary slackness constraints, which can be transformed into a mixed-integer linear program (MILP) using disjunctive constraints. We refer to~\ref{app:a} for the KKT conditions of problems \eqref{mod:A-Res} and \eqref{mod:A-DA}. 

Regarding the structure of the bilevel model \eqref{mod:B-Pre}, it should be noted that, in contrast to the reserve capacity and the day-ahead markets, the balancing market clearing is modeled in the upper-level problem by the last term of the objective~function~\eqref{eq:B-pre-obj} and constraints \eqref{eq:B-A-CV}-\eqref{eq:B-A-CV-flow}. 
The proposed structure leverages the fact that
 the variables of the balancing market in~\eqref{mod:A-CV-rt} do not enforce any restriction on the upper-level variables $\{ \chi_e',\,\forall e\}$ and also they do not impact the lower-level problems \eqref{eq:LLR}~and~\eqref{eq:LLD}. Following this observation, the proposed formulation reduces the computational complexity of the final MILP, since it avoids the integer reformulation of the balancing market complementarity conditions for each scenario~$s\in\mathcal{S}$. For similar applications of bilevel programming in the~electricity markets, the interested reader is referred to~\citep{pineda2016capacity,morales2014electricity,dvorkin2018setting,jensen2017cost}.

Having defined the main properties of the preemptive model, the following comments are in order.
Define the sum of the costs~\eqref{eq:A-Res-obj},~\eqref{eq:DA-Obj} and~\eqref{eq:A-CV-rt-Obj}, as $J^{s'}(\emptyset)$ and $J^{s'}(\AC)$, when $\chi_e, \forall e$ are fixed to the existing cross-border arrangements and to the optimal $\{\hat\chi_e, \forall e\}$ from the preemptive model, respectively. From an economic intuition, it follows that the preemptive model reduces the total expected cost since the establishment of broader coalitions enlarges the pool of available reserves and balancing resources, and enables the more efficient allocation of available generation capacity between these services.
This can be mathematically stated as $J(\AC)=\mathbb E_s[J^s(\AC)] \leq J(\emptyset) =\mathbb E_s[J^s(\emptyset)]$, where $\mathbb{E}_s[\cdot] $ is the expectation calculated over the scenario set $\mathcal{S}$.
Using a similar reasoning, it follows that $J$ is a nonincreasing function with respect to the number of areas participating in the coalition, i.e., $J(\CC)\leq J(\hat \CC)$ for all $\hat \CC\subseteq\CC.$ (Note that the preemptive model does not guarantee $J^s(\CC)\leq J^s(\hat \CC),$ $\forall\hat \CC\subseteq\CC$ for each scenario independently.) 
It becomes apparent that the implementation of the preemptive model results in a different {system} cost than the one under the current sequential market. Next, we discuss cost allocation in this new sequential market.

\subsection{Cost allocation for the preemptive model}

The preemptive model \eqref{mod:B-Pre} implements a centralized transmission allocation mechanism, under the implicit assumption that all areas are willing to accept the $\{\hat\chi_e, \forall e\}$ solution that by construction minimizes the system-wide expected cost. However, this model does not suggest an area-specific cost allocation that guarantees sufficient incentives for all areas to remain in the grand coalition~$\AC$.

The task of setting the transmission shares and allocating the resulting costs could be accomplished through the establishment of a new trading floor, cleared before the reserve capacity market, in which regional operators would place their bids/offers for the reservation of inter-area transmission capacities, akin to the decision variable of the preemptive model, $\{\chi_e', \forall e\}$. This new market would constitute an ideal benchmark of the market-based allocation process described in \citep{EC1}, implementing a \textit{complete} market for inter-area transmission allocations in which capacities would be traded based on bids/offers that reflect the valuations from regional operators.
In this work, we follow an alternative path to promote the formation of stable coalitions for the exchange of reserves. %and ensure the efficient allocation of the resulting benefits. 
Our approach builds an ex-post benefit allocation mechanism on top of the preemptive model, aiming to realize the necessary conditions that accomplish the coordination requirements of this model. The proposed mechanism is readily compatible with the existing market structure and does not require the establishment of any new marketplace.
In the remainder of this section, we outline the main concepts related to benefit allocations for the preemptive model and we discuss the desirable properties that we want to achieve. 
\looseness=-1

Let $J^s_a(\emptyset)$ denote the cost allocated to area~$a$ in scenario~$s$ in the existing sequential market. As previously discussed, the current implementation of the sequential market provides a cost allocation method that satisfies budget balance under every scenario, i.e., $J^s(\emptyset)=\sum_{a\in\mathcal A}J^s_a(\emptyset)$. 
The implementation of the preemptive model requires a new method to allocate costs to the areas that participate in this arrangement. This task can equivalently be viewed as allocating benefits based on the change in the total cost as a discount or a mark-up on the original cost allocation of each area defined by $J^s_a(\emptyset)$. 
While choosing these benefits, our main goal is to ensure that all areas in~$\AC$ are willing to use the preemptive model as a decision-support tool, since otherwise some areas may opt for having their own reserve {exchange} agreement. Moreover, we should aim to form coalitions as large as possible in order to achieve the highest reduction in the expected system cost. 

To achieve these objectives, we will treat the preemptive model in a coalitional game framework, which allows us to approach the benefit allocation problem in two ways. First, we can allocate the expected cost reduction, $J(\emptyset)-J(\AC)\geq 0$, to all areas as benefits. Allocating benefits this way achieves budget balance in expectation, which implies that there is no deficit or surplus if the preemptive model is used repeatedly and the uncertainty modeling is accurate enough.\footnote{This conclusion requires that the scenario set fully describes the uncertainties. In practice, this is inevitably an approximation to the real world. There are various results showing asymptotic guarantees in a convex optimization framework as long as the scenario set is rich enough, see \citep{birge2011introduction}} However, this method does not guarantee that the resulting allocation satisfies budget balance in every scenario, thus requiring a large financial reserve to buffer the fluctuations in the budget in case of surplus or deficit for some uncertainty realizations.
The second approach is to allocate the scenario-specific cost variation, $J^s(\emptyset)-J^s(\AC)$. This method would guarantee budget balance for every scenario.  {However}, the areas participating in the coalition would {collect} benefits that vary under different scenarios, {possibly raising also risk considerations}.

\section{Coalitional Game Theory Framework}\label{sec:3}

% In this section, we bring in the preliminaries for benefit allocation mechanisms and {their} desirable properties. We then review {existing} mechanisms from the literature and discuss whether they {attain} these properties. {To facilitate the exposition}, this section treats a general coalitional game. 
%\looseness=-1%{whereas} in Section~\ref{sec:4}, we focus on two coalitional games arising from the preemptive model {that lay the foundations for the proposed} methods.

A coalitional game is {defined} by a set of players and {the so-called} coalitional value function, that maps from the subsets of players to the values, i.e., the total benefits created by these players~\citep{osborne1994course}. In the preemptive model, the set of players are given by the set of areas $\AC$, whereas the coalitional value function $v:2^{\AC}\rightarrow \mathbb{R}$ can be defined 
either as the expected cost reduction achieved, i.e., $\bar v(\CC) = J(\emptyset) - J(\CC),$ for all $\CC\subseteq \AC$ or based on the resulting change in the cost of the realized scenario $s\in\mathcal{S}$, i.e., $v^s(\CC) = J^s(\emptyset) - J^s(\CC),$ for all $\CC\subseteq \AC$.
% in two ways. The first approach is to define $v$ as the expected cost reduction achieved, i.e., $\bar v(\CC) = J(\emptyset) - J(\CC),$ for all $\CC\subseteq \AC$. {A second definition of the value function is based on} the resulting change in the cost of the realized scenario $s\in\mathcal{S}$, i.e., $v^s(\CC) = J^s(\emptyset) - J^s(\CC),$ for all $\CC\subseteq \AC$.
Clearly, {it holds that} $\bar v(\CC)=\mathbb{E}_s[v^s(\CC)].$ Later, we will see that these functions yield different structures for the coalitional game. In the remainder of this section, we study a generic function~$v$ for the preemptive model satisfying $v(\CC)=0$, for all $|\CC|\leq 1$. This assumption holds since coordination is not possible in the preemptive model without the participation of at least two adjacent areas.

Given the coalitional value function $v$, \textit{a benefit allocation mechanism} defines the benefit received by each area $a\in\AC$ with $\beta_a(v)\in \R$. The cost allocated to area $a$ under the preemptive transmission allocation model would then be given by $J^s_a(\AC)=J^s_a(\emptyset)-\beta_a(v)$. Depending on its sign, the benefit can be considered as a discount or a mark-up on the original cost allocation.

{When designing benefit allocation mechanisms, there are three }fundamental properties we want to guarantee, namely, efficiency, individual rationality, and stability. %Notice that these properties pertain to the economic side, and hence they can also be referred to as the economic properties. Moreover, we desire an additional computational property, which will be explained later on.
A benefit allocation $\beta(v) = \{\beta_a(v)\}_{a\in\mathcal A}\in\mathbb{R}^{\mathcal A}$ is \textit{efficient} if {the whole} value created by the grand coalition, {i.e., $\CC = \AC$}, is allocated to the {member}-areas, i.e., $\sum_{a\in\AC}\beta_a(v)=v(\AC)$.\footnote{In coalitional games, efficiency is also often referred to as budget balance. For clarity, we use efficiency for the benefit allocation, and the term budget balance is reserved for the cost allocation.} A benefit allocation ensures \textit{individual rationality} if all areas {obtain} nonnegative benefits, i.e., $\beta_a(v)\geq 0,$ for all $a\in\AC$. If this property does not hold, the coordination arrangement would yield increased costs for some areas. As a result, these areas may decide not to participate in the preemptive model. Finally, a benefit allocation attains \textit{stability} (in other words, group rationality) if it eliminates the incentives of the areas to form {sub-}coalitions, i.e., $\nexists \mathcal C\subset\AC$ such that $v(\CC)>\sum_{a\in\CC}\beta_a(v)$. This last property is crucial for the preemptive model, since otherwise some areas may opt for having their own reserve {exchange} agreement by excluding the remaining areas. This coincides with our aforementioned goal of ensuring that all areas participate in the preemptive model.

In coalitional game theory, these properties are known to be attained if the benefit allocation lies in the \textit{core} {defined as} $\beta(v) \in K_{\text{Core}}(v)$, where
$K_\text{Core}(v)=\{\beta\in \mathbb{R}^\AC\,|\, \sum_{a\in\AC}\beta_a=v(\AC),\ \sum_{a\in\CC}\beta_a\geq v(\CC),\ \forall \CC\subset\AC\}.$ %closed polytope if nonempty
In this definition, the equality constraint ensures efficiency, while inequality constraints guarantee stability, i.e., there is no subset of areas $\CC\subset\AC$ that can yield higher total benefits for its members compared to the benefit allocation under the grand coalition.
The inequality constraints also include $\beta_a(v)\geq v(a)=0$ for all $a\in\AC$.\footnote{For the sake of simplicity, singleton sets are denoted by $a$ instead of $\{a\}.$} This restriction ensures individual rationality. 

The core is a closed polytope involving $2^{|\AC|}$ linear constraints. %This polytope can be guaranteed to be nonempty only in restricted settings.
This polytope is nonempty if and only if the coalitional game is balanced~\citep{shapley1967balanced}. Such settings include the cases in which the coalitional value function exhibits supermodularity\footnote{Supermodularity is attained if for any set the participation of an area results in a larger value increment when compared to the subsets of the set under consideration, i.e., $v(\CC\cup\{a\})-v(\CC)\geq v(\CC'\cup\{a\})-v(\CC')$, $\forall a\notin\CC,\CC'\subset\CC\subseteq\AC.$}~\citep{shapley1971cores} and the cases in which the coalitional value function can be modeled by a concave exchange economy~\citep{shapley1969market}, a linear production game~\citep{owen1975core} or a risk-sharing game~\citep{csoka2009stable}. In their most general form, the coalitional value functions in these works are given by an optimization problem minimizing a convex objective subject to linear constraints. In the problem at hand, coalitional value functions are associated with solutions to the general non-convex optimization problem \eqref{mod:B-Pre}. As a result, previous works on the nonemptiness of the core are not applicable to our setup.
\looseness=-1

In case the core is empty, we need to devise a method to approximate a core allocation. To this end, we bring in the {notion of} strong $\epsilon$-core, defined in \citep{shapley1966quasi} as 
$K_\text{Core}(v,\epsilon)=\{\beta\in \mathbb{R}^\AC\,|\, \sum_{a\in\AC}\beta_a=v(\AC),\ \sum_{a\in\CC}\beta_a\geq v(\CC)-\epsilon,\ \forall \CC\subset\AC\}.$
This definition can be interpreted as follows. If organizing a coalitional deviation entails an additional cost of $\epsilon\in\mathbb{R}$, coalition values would be given by $v(\CC)-\epsilon$ for all $\CC\neq\AC$. Then, the resulting core would correspond to the strong $\epsilon$-core. For $\epsilon=0$, we retrieve the original core definition, i.e., $K_\text{Core}(v,0)=K_\text{Core}(v)$. 

Let $\epsilon^*(v)$ be the critical value of $\epsilon$ such that the strong $\epsilon$-core is nonempty, {which is mathematically defined as} $\epsilon^*(v)=\min\{\epsilon\,|\,K_\text{Core}(v,\epsilon)\not=\emptyset\}$. The value $\epsilon^*(v)$ is guaranteed to be finite for any function~$v$ and the set $K_\text{Core}(v,\epsilon^*(v))$ is called the \textit{least-core}~\citep{maschler1979geometric}. 
Let the excess of a coalition be defined by $\theta(v,\beta,\CC)=v(\CC)-\sum_{a\in\CC}\beta_a,$ for any nonempty $\CC\subset\AC$. In other words, the set $K_\text{Core}(v,\epsilon^*(v))$ is the set of all efficient benefit allocations minimizing the maximum excess.
If the core is empty, the maximum excess is the maximum violation of a stability constraint. This implies that the least-core achieves an approximate stability property. As a remark, the least-core relaxes also the inequality constraints corresponding to singleton sets $\beta_a\geq v(a)-\epsilon^*(v)=-\epsilon^*(v)$ for all $a\in\AC$ since $\epsilon^*(v)> 0$, and hence it yields approximate individual rationality. Finally, if the core is not empty, we have $\epsilon^*(v)\leq 0$ and the least-core is a subset of the core. %, which can be considered as centrally located

With the discussion above, we conclude that whenever the core is empty, we can use the least-core to achieve the second best outcome available, i.e., a benefit allocation which is efficient, approximately individually rational and approximately stable. Observe that there are generally many points to choose from the least-core (or the core if it is nonempty) achieving the same fundamental properties. In this case, it could be desirable to require additional intuitively acceptable properties to pick a unique benefit allocation. Later, we revisit this idea in our proposed methods.
 
Apart from the aforementioned fundamental properties that pertain to the economic side of the problem, computational tractability is also a practical concern, considering that we may need the complete list of coalition values $v(\CC)$ for all $\CC\subseteq\AC$ to fully describe the core and the least-core.
For the coalitional games arising from the preemptive model, each coalition value requires another solution to MILP in~\eqref{mod:B-Pre}, which is NP-hard in general. 
% Other than problem~\eqref{mod:B-Pre} itself being hard to solve, this computation would grow exponentially in the number of areas and quickly become intractable for large instances.
Hence, our goal is to find a core or a least-core benefit allocation that can be computed with limited queries to the coalitional value function. Next, we briefly review two benefit allocation mechanisms that are widely used in the literature. %and we discuss their properties with respect to the core and least-core definitions as well as their computational efficiency.

\subsection{Shapley value}
The benefit assigned by the Shapley value is given by
$\beta_a^{\text{Shapley}}(v) = \sum_{\CC\subseteq\AC} \frac{(|\CC|-1)!(|\AC|-|\CC|)!}{|\AC|!}(v(\CC)-v(\CC\setminus a)).$
This benefit is the average of the marginal contribution of the area $a$ under all coalitions, considering also all possible orderings of areas. % \citep{shapley1953value}
The Shapley value results in an efficient benefit allocation. Individual rationality is also satisfied if the coalitional value function is nondecreasing, since the marginal contributions would be nonnegative. On the other hand, the Shapley value is guaranteed to lie in the core only when the coalitional value function is supermodular~\citep{shapley1971cores}. This is a restrictive condition that is not applicable to our problem. In addition, when the core is empty, the Shapley value does not necessarily lie in the least-core, making it incompatible with the fundamental properties we desire~\citep{maschler1979geometric}. 
In terms of the computational performance, the calculation of the Shapley value requires the exhaustive enumeration of coalition values $v(\CC)$ for all $\CC\subseteq\AC$. Finally, it should be noted that the Shapley value is the unique efficient benefit allocation that satisfies dummy player, symmetry, and additivity properties simultaneously.
Dummy player property requires $\beta_a=0$ for all $a$ for which $v(\CC)-v(\CC\setminus a)=0$ for all~$\CC\subseteq\AC$. In other words, an area incapable of contributing to any coalition~$\CC$ ends up with zero benefits. Next, we show the relation between the previously discussed properties and the dummy player property.
%The core provides us with a more restrictive version of the dummy player property, i.e., $\beta_a=0$ for all $a$ such that $v(\AC)-v(\AC\setminus\{a\})=0$. For the least-core, we have the following result.
%It is straightforward to check that for such areas 

\begin{proposition}\label{prop:dummy} For the core and the least-core, we have, (i) if $a'$ satisfies $v(\AC)-v(\AC\setminus a')=0$, then $K_\text{Core}(v)\subset\{\beta\,|\,\beta_{a'}=0\}$, (ii) if $a'$ satisfies $v(\CC)-v(\CC\setminus a')=0$ for all~$\CC\subseteq\AC$, then $K_\text{Core}(v,\epsilon^*(v))\subset\{\beta\,|\,\beta_{a'}=0\}$.
\end{proposition}

This result shows that the core attains a more restrictive version of the dummy player property, i.e., an area incapable of contributing to the set $\AC$ ends up with zero benefits. Finally, the least-core attains the dummy player property in the same way that it is defined for the Shapley value.
The proof and the discussions on symmetry and additivity are relegated to~\ref{app:a2}.

\subsection{Nucleolus allocation}

 Among all efficient benefit allocations, the nucleolus allocation is the unique benefit allocation that minimizes the excesses of all coalitions in a lexicographic manner~\citep{schmeidler1969nucleolus}. Nucleolus allocation lies in the least-core and hence attains the desirable economic properties. In terms of practical implementation, the lexicographic minimization is computationally demanding in the general case. Nucleolus allocation can be computed by solving a sequence of ${\mathcal{O}(|\AC|)}$ linear programs with constraint sets that are parametrized versions of the core $K_{\text{Core}}(v)$, see~\citep{kopelowitz1967computation,fromen1997reducing}.
%  The formulation in \citep{kopelowitz1967computation} was shown to require solving $\mathcal{O}(|\AC|)$ linear programs to obtain linearly independent equalities that uniquely characterize the nucleolus~\citep{fromen1997reducing}. 
 However, each linear program requires the complete list of coalition values. In case the coalition values are given implicitly by the objective value of a single linear optimization problem with constraints depending on the participants of the coalition, the work by~\citet{hallefjord1995computing} proposes using constraint generation algorithms. 
In this approach, each linear program is solved by a separate constraint generation algorithm that iteratively generates coalitional values on demand. Nevertheless, we may still need to generate all possible coalition values~\citep{hallefjord1995computing,kimms2012approximate}. When the number of areas is large, this approach involving the execution of the constraint generation algorithm $\mathcal{O}(|\AC|)$ times becomes computationally prohibitive for our application.
 {For the sake of completeness,}~\ref{app:a2-2} provides the mathematical definition for the nucleolus allocation and its comparison with the Shapley value.

\section{Benefit Allocation Mechanisms for Preemptive Transmission Allocation}\label{sec:4}
%In this section, we propose two benefit allocation mechanisms that can be applied at the day-ahead or at the real-time market trading floors. 
In this section, the first benefit allocation mechanism, which is an ex-ante process with respect to the uncertainty realization, employs as coalitional value function the expected cost reduction from the preemptive transmission allocation model. The second mechanism is an ex-post process that can be applied only when the scenario $s\in\SC$ is unveiled, since it uses as coalitional value function the scenario-specific cost variation.
For each coalitional game, we first discuss the implications of the structure of the game on the tools we reviewed from the literature in Section~\ref{sec:3}. 
We then propose our approach, the least-core selecting mechanism. Our mechanism comes with a computational improvement over the nucleolus and it can approximate additional criteria, for instance, fairness.

\subsection{Benefit allocations for expected cost reduction}\label{sec:4a}

{For the ex-ante allocation mechanism, }the coalitional value function $\bar v(\CC) = J(\emptyset) - J(\CC)\geq 0$ for all $\CC\subseteq \AC$ is nondecreasing, since $J$ is nonincreasing.
Given the function~$\bar v$, an efficient benefit allocation, $\sum_{a\in\AC} \beta_a(\bar v)=\bar v(\AC)$, would result in a cost allocation that is budget balanced in expectation, since $J(\AC)= J(\emptyset) - \bar v(\AC) = \mathbb{E}_s\big[\sum_{a\in\AC}J^s_a(\emptyset)\big] - \sum_{a\in\AC} \beta_a(\bar v)= \mathbb{E}_s\big[\sum_{a\in\AC}J^s_a(\AC)\big].$

While designing a benefit allocation mechanism, our goal is to achieve the three fundamental properties, i.e., efficiency, individual  rationality and stability, associated with the core $K_{\text{Core}}(\bar v)$. However, as already mentioned, the previous results on the nonemptiness of the core are not applicable to our problem. Thus, we have the following condition that is applicable to some specialized instances of the coalitional value function~$\bar v$.

\begin{proposition}\label{prop:corempt}
	$K_\text{Core}(\bar v)$ is nonempty if there exists an area $a'\in\AC$ such that $\bar v(\AC\setminus a')=0$.
\end{proposition}

The proof is relegated to~\ref{app:a3}. 
Note that this condition can only be found in specialized instances of the preemptive  model. For instance, in the case of a star graph $(\mathcal A,\mathcal E)$, the central area would satisfy this condition, since it is indispensable for enabling any reserve exchange.
%It is not possible in general to guarantee that the graph has a star formation.
However, in a general graph, the core could potentially be empty and we focus on this case in the illustrative example provided in Section~\ref{sec:emptycore}.

In case of an empty core, our goal is to achieve {a least-core solution}, {which can be perceived as} the second best outcome in our context. Other than approximating the stability property, the least-core also approximates the individual rationality property by relaxing the inequality constraints for singleton sets, i.e., $\beta_a\geq \bar v(a)-\epsilon^*(\bar v)=-\epsilon^*(\bar v)$ for all $a\in\AC$. The following proposition shows that the least-core is compatible with the individual rationality property for the function~$\bar v$.

\begin{proposition}\label{prop:leastcore_volp}
	$K_\text{Core}(\bar v,\epsilon^*(\bar v))$ lies in $\mathbb{R}^\AC_+$.
\end{proposition}

The proof is relegated to~\ref{app:a4}. It relies on the observation that any least-core allocation violating the individual rationality would imply the existence of an $\epsilon <\epsilon^*$, such that $K_\text{Core}(\bar v,\epsilon)$ is nonempty, contradicting the definition of the least-core. Thus, we can use the least-core to achieve efficiency, individual rationality and approximate stability, whenever the core is empty.

For this coalitional game, the Shapley value satisfies efficiency and individual rationality, but stability (or approximate stability) and computational tractability are not attained. We provide an example for stability violation in Section~\ref{sec:shapno}. The nucleolus allocation, on the other hand, lies in the least-core and it satisfies efficiency, individual rationality and approximate stability.\footnote{Note that the function~$\bar v$ is given implicitly by the optimization problem~\eqref{mod:B-Pre}. For this problem, constraints \eqref{eq:B-pre-ta} and \eqref{eq:B-pre-ta-fixed} change depending on the participants of the coalition. This would allow us to generalize the computation of the nucleolus via the sequential constraint generation algorithms provided in~\citep{hallefjord1995computing}.} 
\looseness=-1

Based on the previous discussions, we propose a least-core selecting mechanism that is mathematically formulated as the following optimization problem:
\begin{equation}%+ \epsilon\, ||\beta-\beta^{\text{c}}||_2^2
\label{eq:ba_opt}
\underset{  \epsilon,\, \beta}{\text{minimize}}\ \epsilon\quad \text{subject to}\ \epsilon\geq 0,\ \beta\in K_{\text{Core}}(\bar v,\epsilon).
\end{equation}
Let $\hat{\epsilon}$ denote the optimal value of $\epsilon$ for this problem. If the core is empty, we have $\hat\epsilon=\epsilon^*(\bar v)>0$ and problem~\eqref{eq:ba_opt} finds a least-core allocation. On the other hand, if the core is nonempty, we have $\hat\epsilon=0$ and problem~\eqref{eq:ba_opt} finds instead a core benefit allocation, which attains properties of efficiency, individual rationality and  stability.
The nucleolus allocation always forms an optimal solution pair with $\hat{\epsilon}$ to problem~\eqref{eq:ba_opt}, since it lies in the least-core. In fact, there are in general many optimal solutions to this problem. To this end, we will propose an additional criterion for tie-breaking. 

Let $\beta^{\text{c}}$ be a desirable and a fair benefit allocation that is easy to compute but not necessarily in the core or in the least-core. An example could be the marginal contribution of each area $\beta^{\text{m}}:\ \beta_a^{\text{m}} = \bar v(\AC)-\bar v(\AC\setminus a)$ for all $a\in\AC$ which requires $|\AC|+2$ calls to problem~\eqref{mod:B-Pre}. %This choice implies that the areas contributing the most in expectation should receive larger benefits. 
Receiving the marginal contribution can be regarded as a fair outcome.\footnote{This allocation coincides with the Vickrey-Clarke-Groves mechanism, which is equivalent to the externality of each area~\citep{krishna2009auction}. In an auction, this allocation ensures that truthfully reporting the preferences is a dominant strategy Nash equilibrium. Note that we assumed no strategic behaviour from the regional operators.} This allocation satisfies individual rationality, dummy player and symmetry properties. However, it is generally not efficient, see~\citep{krishna2009auction}, and not stable, see~\citep{karaca2018core}.
Another example could be $\beta^{\text{eq}}=(\bar v(\AC)/|\AC|)\bm {1}^\top$ which assigns equal importance to each area. This choice requires two calls to problem~\eqref{mod:B-Pre} and it satisfies efficiency and individual rationality. However, this allocation also violates stability. 

Starting from such a desirable benefit allocation, we can solve the following problem 
\begin{equation}%+ \epsilon\, ||\beta-\beta^{\text{c}}||_2^2
\label{eq:ba_opt_tie}
\underset{\beta}{\text{minimize}}\ ||\beta-\beta^{\text{c}}||_2^2\quad \text{subject to}\ \beta\in K_{\text{Core}}(\bar v, \hat \epsilon),
\end{equation}
to obtain a unique benefit allocation for problem~\eqref{eq:ba_opt}.
The uniqueness follows from having a strictly convex objective.
Let $\hat{\beta}(\bar v,\beta^{\text{c}})$ denote the optimal value of $\beta$ in problem~\eqref{eq:ba_opt_tie}. We define the benefit allocation $\hat{\beta}(\bar v,\beta^{\text{c}})$ as the\textit{ least-core selecting mechanism}. {This allocation} achieves economic properties of the least-core, and also the core if the core is nonempty, while approximating an additional criterion defined by $\beta^{\text{c}}$. For instance, if the marginal contribution $\beta^{\text{m}}$ is chosen, problem~\eqref{eq:ba_opt_tie} would pick the allocation $\hat{\beta}(\bar v,\beta^{\text{m}})$ approximating the fairness of the marginal contribution.

Characterizing the constraint sets of problems~\eqref{eq:ba_opt}~and~\eqref{eq:ba_opt_tie} still requires exponentially many solutions to \eqref{mod:B-Pre}. We show that~\eqref{eq:ba_opt}~and~\eqref{eq:ba_opt_tie} can be solved by a single constraint generation algorithm. 

\textit{Constraint generation algorithm:} 

Here, we describe the steps of the constraint generation algorithm at iteration $k\geq 1$. Let $\mathcal{F}^k\subset 2^{\AC}$ denote the family of coalitions for which we have already generated the coalition values. The algorithm first obtains a candidate solution by solving a relaxed version of problem~\eqref{eq:ba_opt} as follows:
\begin{equation}
\label{mod:candidate_k}
\underset{  \epsilon,\, \beta}{\text{minimize}}\ \epsilon\quad
\text{subject\ to}\ \epsilon\geq 0,\quad \beta\geq0,\quad\sum_{a\in\AC}\beta_a=\bar v(\AC),\quad \sum_{a\in\CC}\beta_a\geq \bar v(\CC)-\epsilon,\quad \forall \CC\in\mathcal{F}^k. 
\end{equation}
Let the optimal solution be denoted by $\epsilon^{k}$. Clearly, it satisfies $\hat\epsilon\geq\epsilon^{k}.$
Then, we solve the following problem as a tie-breaker:\begin{equation}
	\label{mod:candidate_k_unique}
\underset{\beta}{\text{minimize}}\ ||\beta-\beta^{\text{c}}||_2^2\quad
\text{subject\ to}\ \beta\geq0,\quad\sum_{a\in\AC}\beta_a=\bar v(\AC),
\quad\sum_{a\in\CC}\beta_a\geq \bar v(\CC)-\epsilon^k,\quad \forall \CC\in\mathcal{F}^k.
\end{equation}
Denote the optimal benefit allocation for problem~\eqref{mod:candidate_k_unique} by~$\beta^{k}$. This allocation would form an optimal solution pair to problem~\eqref{mod:candidate_k} with~$\epsilon^{k}$. 
In principle, $\mathcal{F}^1$ can potentially be chosen as an empty set, by setting $\epsilon^{1}$  equal to zero and removing the last set of constraints in~\eqref{mod:candidate_k_unique}.
% Note that $\mathcal{F}^1$ can potentially be chosen as an empty set. In this case, we set $\epsilon^{k}$ to be $0$, and remove the constraints \eqref{eq:modified_uniq} from problem~\eqref{mod:candidate_k_unique}.

Given a candidate allocation $\beta^k$, we can then generate the coalition with the maximum stability violation by solving the following problem, which treats $\beta^k$ as a fixed parameter:
\begin{flalign}\label{eq:ba_findviolation}
	&\underset{ \CC\subseteq\AC}{\text{maximize}}\ \bar v(\CC)-\sum_{a\in\CC}\beta_a^k .
\end{flalign}
Denote the optimal solution by $\CC^k$, and the optimal value by~$\eta^k$.
Using the fact that the coalitional value function~$\bar v$ is given implicitly by the MILP version of the preemptive model~\eqref{mod:B-Pre}, we can show that problem \eqref{eq:ba_findviolation} has an equivalent MILP reformulation, which can be solved by off-the-shelf optimization solvers. Using this approach, we eliminate the need for evaluating $\bar v(\CC)$ for all $\CC\subseteq \AC$ to solve problem~\eqref{eq:ba_findviolation}.
This MILP reformulation is relegated to~\ref{app:a6}. 
Finally, this problem generates the coalition value $\bar v(\CC^k)=\eta^k+\sum_{a\in\CC^k}\beta_a^k,$ and we add $\CC^k$ to $\mathcal{F}^k$ for the next iteration.

In order to define a stopping criterion ensuring that the iterative solution of problems~\eqref{mod:candidate_k},~\eqref{mod:candidate_k_unique} and~\eqref{eq:ba_findviolation} converges to the optimal solution of problem~\eqref{eq:ba_opt_tie}, we need the following two observations. First, if $\epsilon^k>0$ in problem~\eqref{mod:candidate_k}, then there exists a set $\CC\in\mathcal{F}^k$ such that $\bar v(\CC)-\sum_{a\in\CC}\beta_a^k=\epsilon^k$, which implies that $\eta^k\geq\epsilon^k$. On the other hand, if $\epsilon^k=0$, by setting $\CC=\AC$ we can show that $\eta^k\geq\bar v(\AC)-\sum_{a\in\AC}\beta_a^k=0=\epsilon^k$. Based on these remarks, we have $\eta^k\geq\epsilon^k$ for any iteration $k$.

The iterative solution of problems~\eqref{mod:candidate_k},~\eqref{mod:candidate_k_unique}~and~\eqref{eq:ba_findviolation} terminates when $\eta^k=\epsilon^k$. In this case,
 problem \eqref{eq:ba_findviolation} provides a certificate that the pair $(\epsilon^k,\beta^k)$ is a feasible solution to problem~\eqref{eq:ba_opt}. Note that this pair is also optimal to a relaxed version of problem~\eqref{eq:ba_opt}, given by problem~\eqref{mod:candidate_k}. This concludes that $(\epsilon^k,\beta^k)$ is optimal for problem~\eqref{eq:ba_opt}.
Observing that $\epsilon^k=\hat \epsilon$ and using a similar reasoning, we conclude that $\beta^k$ is the optimal solution to problem~\eqref{eq:ba_opt_tie}. Hence, the algorithm converges. 
\looseness=-1

In the intermediate solution points of the iterative process, we have $\eta^k>\epsilon^k$ and consequently $\CC^k\notin \mathcal{F}^k$ according to the optimality of~$\epsilon^k$ for problem~\eqref{mod:candidate_k}. We then extend the family of generated coalitions by $\mathcal{F}^{k+1}=\mathcal{F}^k\cup\CC^k$ until convergence is achieved.
Since there are finitely many coalitions to be generated, the algorithm converges after a finite number of iterations. As a remark, the algorithm would not generate any set from $\mathcal{F}^1$, the full set, the empty set, and the singleton sets (since we enforce $\beta\geq0$ in problem~\eqref{mod:candidate_k}). In practice, even when there are many areas, the algorithm requires the generation of only several coalition values. We show this in the numerical results. 
Finally, we summarized this iterative algorithm in~\ref{app:a6}.

\subsection{Benefit allocations per scenario}

As previously explained, allocating benefits for the expected cost reduction does not guarantee that the resulting cost allocation satisfies budget balance in every scenario. Having a surplus or a deficit might be undesirable, since this may necessitate a large financial reserve to buffer the fluctuations in the budget. To address this issue, here we focus on the allocation of the {the scenario-specific cost variation}, $J^{s}(\emptyset)-J^{s}(\AC)$.
The coalitional value function in this case is given by $v^{s}(\CC) = J^{s}(\emptyset) - J^{s}(\CC),$ for all $\CC\subset \AC$. Observe that the set function $v^{s}$ is not necessarily nondecreasing, while it can also map to negative reals, since the preemptive transmission allocation model does not guarantee that $J^{s}(\CC)\leq J^{s}(\emptyset)$ holds. Given the function $v^{s}$, an efficient benefit allocation mechanism, $\sum_{a\in\AC} \beta_a(v^{s})=v^{s}(\AC)$, would result in a cost allocation that is budget balanced in scenario~${s}$, since $J^{s}(\AC)= J^{s}(\emptyset) - v^{s}(\AC) = \sum_{a\in\AC}J^{s}_a(\emptyset) - \sum_{a\in\AC} \beta_a(v^{s})= \sum_{a\in\AC}J^{s}_a(\AC).$

Aiming at establishing a per-scenario benefit allocation, our goal now is to achieve the properties of the scenario-specific core $K_\text{Core}(v^{s})$. However, neither the previous results nor Proposition~\ref{prop:corempt} apply to this core to prove that it is nonempty 
as it can be affirmed by the following result.
\begin{proposition}\label{prop:empty}
	$K_\text{Core}(v^{s})$ is empty if there exists $\CC\subset\AC$ such that $v^{s}(\AC)< v^{s}(\CC)$.
\end{proposition}

The proof is relegated to~\ref{app:a7}. In practice, the condition above would prevent the formation of the grand coalition $\AC$,
as shown in the illustrative example of Section \ref{sec:shapno}. 
The coalition value $v^{s}(\AC)$ being negative is a special case of Proposition~\ref{prop:empty}, since we would then have $v^{s}(\AC)< v^{s}(a)=0$ for all $a\in \AC$. %Moreover, if $v^{s}(\AC)<0$ for some $s\in\mathcal S$, then it is straightforward to conclude that the individual rationality property is not compatible with achieving budget balance in this scenario. 
We see that it may not be realistic to achieve all three fundamental properties, and we should instead aim for the least-core $K_\text{Core}(v^s,\epsilon^*(v^s)).$ Note that in this case Proposition~\ref{prop:leastcore_volp} is not applicable and the least-core would instead achieve efficiency, approximate individual rationality, and approximate stability.

For the coalitional game arising from the function~$v^s$, the Shapley value satisfies efficiency, but individual rationality, stability, and computational tractability are not attained. On the other hand, the nucleolus allocation provides a least-core allocation. Note that, in contrast to the expected coalitional value function~$\bar v$, the function~$v^s$ is not implicitly given by an optimization problem. Instead, it is an ex-post calculation from the sequential electricity market after the uncertainty realization. As a result, the coalitional value function~$v^s$ is not amenable to a constraint generation approach. This would prevent the iterative computations of the nucleolus and also the least-core selecting mechanism from problems~\eqref{eq:ba_opt}~and~\eqref{eq:ba_opt_tie} in case they are formulated for $v^s$.

Next, we look at an alternative approach that can be computed in a computationally tractable manner. This approach will extend our results from Section~\ref{sec:4a}, showing that any efficient benefit allocation for the expected cost reduction gives rise to an efficient scenario-specific benefit allocation that results in budget balance in every scenario. Moreover, using the proposed approach we can achieve the fundamental properties associated with the original benefit allocation in expectation.

Let $\beta(\bar v)\in\mathbb{R}^{\AC}_+$ be an efficient individually rational benefit allocation for the expected cost reduction, computed prior to the uncertainty realization. We then define the following scenario-specific benefit allocation, $\beta(v^{s},\beta(\bar v))\in\mathbb{R}^{\AC}$,
\begin{equation} \label{eq:defscendep}
\beta_{a}(v^{s},\beta(\bar v))=\dfrac{\beta_a(\bar v)}{\bar v(\AC)} v^{s}(\AC),\quad \forall a\in\AC.
\end{equation}
The benefit $\beta_{a}(v^{s},\beta(\bar v))$ for each area $a$ is computed based on an ex-post computation of $v^{s}(\AC)$ for the specific uncertainty realization $s$. Given $\beta(\bar v)$, this definition does not require any further solutions to the preemptive transmission allocation model or the sequential market. The term ${\beta_a(\bar v)}/{\bar v(\AC)}\in[0,1]$ can be considered as a percentage share of profits/losses depending on the sign of $v^s(\AC)$. (This also holds for any other weighting from the $|\AC|$-simplex.)
Notice that since $\sum_{a\in\AC}\beta_{a}(v^s,\beta(\bar v))=v^s(\AC)$, the efficiency property holds. Moreover, {having} $\mathbb E [\beta_{a}(v^s,\beta(\bar v))]=\beta_a(\bar v)$ implies that the scenario-specific benefit allocation $\beta(v^s,\beta(\bar v))$ satisfies in expectation the other fundamental properties of the original benefit allocation~$\beta(\bar v)$. As a remark, in case uncertainty $s\notin \mathcal{S}$ is unveiled, we can still compute the scenario-specific benefits, however, it would not be possible to obtain any guarantees on the properties other than efficiency. %Note that for convex optimization problems there are asymptotic guarantees for out-of-sample realizations, see~\citep{birge2011introduction,zakeri2018pricing}. 

Given the above reasoning, we propose a \textit{scenario-specific least-core selecting mechanism}, which builds upon the least-core selecting benefit allocation mechanism from problems~\eqref{eq:ba_opt}~and~\eqref{eq:ba_opt_tie} to define $\beta(v^{s},\hat{\beta}(\bar v,\beta^{\text{c}}))\in\mathbb{R}^{\AC}$ according to the procedure in \eqref{eq:defscendep}. 
We have previously showed that the allocation $\hat{\beta}(\bar v,\beta^{\text{c}})$ satisfies individual rationality and approximate stability, while enabling a tractable computation via a constraint generation algorithm. In a similar vein, the scenario-specific version $\beta(v^s,\hat{\beta}(\bar v,\beta^{\text{c}}))$ will satisfy individual rationality and approximate stability in expectation, while still enabling a tractable computation. We illustrate this approach in Section \ref{sec:CaseStudies}. As a remark, it is possible to use the Shapley value and the nucleolus allocation in a similar manner. The comparisons of these mechanisms in the previous section would remain unchanged.\footnote{In contrast to the least-core selecting mechanism, $\hat{\beta}(\bar v,\beta^{\text{c}})$, the scenario-specific benefit allocation $\beta(v^s,\hat{\beta}(\bar v,\beta^{\text{c}}))$ would lie in the least-core only in expectation and thus the benefits that each area receives would vary under different scenarios. This implies, in turn, that coalition member areas are now exposed to risk and in case some of these areas are risk-averse, they may not be willing to participate in the benefit allocation mechanism and the preemptive model. Defining ways to incorporate risk measures is part of our ongoing work.}

\section{Numerical Case Studies}
\label{sec:CaseStudies} 
In this section, we first use an illustrative three-area nine-bus system to provide and discuss benefit allocation mechanisms under different system configurations and stochastic renewable in-feed. We then apply the models and the benefit allocation mechanisms in a more realistic case study. All problems are solved with GUROBI 7.5~\citep{gurobi} called through MATLAB on a computer equipped with 32 GB RAM and a 4.0 GHz Intel i7 processor. 
\subsection{Illustrative three-area examples}\label{sec:illu_ex}
\begin{figure*}[h]
	\centering
	\begin{tikzpicture}[scale=0.64, every node/.style={scale=0.5}]
\draw[-,black!80!blue,line width=.1mm] (+0.7,1) -- (+0.7,1.4);
\draw[-,black!80!blue,line width=.1mm] (0.7,1.4) -- (7.3,1.4);
\draw[-,black!80!blue,line width=.1mm] (7.3,1) -- (7.3,1.42);
\draw[-,black!80!blue,line width=.1mm] (+8.7,1) -- (+8.7,1.4);
\draw[-,black!80!blue,line width=.1mm] (8.7,1.4) -- (15.3,1.4);
\draw[-,black!80!blue,line width=.1mm] (15.3,1) -- (15.3,1.4);
\draw[-,black!80!blue,line width=.1mm] (+0.3,1) -- (+0.3,1.45);
\draw[-,black!80!blue,line width=.1mm] (0.3,1.44) -- (1.55,2.4);
\draw[-,black!80!blue,line width=.1mm] (1.55,2.4) -- (2,2.4);
\draw[-,black!80!blue,line width=.1mm] (+8.3,1) -- (+8.3,1.45);
\draw[-,black!80!blue,line width=.1mm] (8.3,1.44) -- (9.55,2.4);
\draw[-,black!80!blue,line width=.1mm] (9.55,2.4) -- (10,2.4);
\draw[-,black!80!blue,line width=.1mm] (+16.3,1) -- (+16.3,1.45);
\draw[-,black!80!blue,line width=.1mm] (16.3,1.44) -- (17.55,2.4);
\draw[-,black!80!blue,line width=.1mm] (17.55,2.4) -- (18,2.4);
\draw[-,black!80!blue,line width=.1mm] (-0.3,1) -- (-0.3,1.45);
\draw[-,black!80!blue,line width=.1mm] (-0.3,1.44) -- (-1.55,2.4);
\draw[-,black!80!blue,line width=.1mm] (-1.55,2.4) -- (-2,2.4);
\draw[-,black!80!blue,line width=.1mm] (7.7,1) -- (7.7,1.45);
\draw[-,black!80!blue,line width=.1mm] (7.7,1.44) -- (6.45,2.4);
\draw[-,black!80!blue,line width=.1mm] (6.45,2.4) -- (6,2.4);
\draw[-,black!80!blue,line width=.1mm] (15.7,1) -- (15.7,1.45);
\draw[-,black!80!blue,line width=.1mm] (15.7,1.44) -- (14.45,2.4);
\draw[-,black!80!blue,line width=.1mm] (14.45,2.4) -- (14,2.4);
\draw[-,black!80!blue,line width=.1mm] (-2,3.6) -- (2,3.6);
\draw[-,black!80!blue,line width=.1mm] (6,3.6) -- (10,3.6);
\draw[-,black!80!blue,line width=.1mm] (14,3.6) -- (18,3.6);
\draw[-,black!80!blue,line width=.1mm] (2,3.8) -- (6,3.8);
\draw[-,black!80!blue,line width=.1mm] (10,3.8) -- (14,3.8);
\draw[-,line width=.5mm] (-1.2,1) -- (1.2,1) node[anchor=west]  {\LARGE $n_3$};
\draw[-,line width=.5mm] (6.8,1) -- (9.2,1) node[anchor=west]  {\LARGE $n_6$};
\draw[-,line width=.5mm] (14.8,1) -- (17.2,1) node[anchor=west]  {\LARGE $n_9$};
\draw[-,line width=.5mm] (-2,1.8) -- (-2,4.2) node[anchor=west]  {\LARGE $n_1$};
\draw[-,line width=.5mm] (6,1.8) -- (6,4.2) node[anchor=west]  {\LARGE $n_4$};
\draw[-,line width=.5mm] (14,1.8) -- (14,4.2) node[anchor=west]  {\LARGE $n_7$};
\draw[-,line width=.5mm] (2,2) -- (2,4.2) node[anchor=west]  {\LARGE $n_2$};
\draw[-,line width=.5mm] (10,2) -- (10,4.2) node[anchor=west]  {\LARGE $n_5$};
\draw[-,line width=.5mm] (18,2) -- (18,4.2) node[anchor=west]  {\LARGE $n_8$};
\draw[<-,line width=.1mm] (0.8,.5) -- (0.8,1) node at (0.8,0.1) {\LARGE $D_3$};
\draw[<-,line width=.1mm] (8.8,.5) -- (8.8,1) node at (8.8,0.1) {\LARGE $D_6$};
\draw[<-,line width=.1mm] (16.8,.5) -- (16.8,1) node at (16.8,.1) {\LARGE $D_9$};
\draw[-,line width=.1mm] (-0.8,.5) -- (-0.8,1);
\draw (-0.8,.15) circle (.35cm) node {\LARGE $i_3$};
\draw[-,line width=.1mm] (7.2,.5) -- (7.2,1);
\draw (7.2,.15) circle (.35cm) node {\LARGE $i_6$};
\draw[-,line width=.1mm] (15.2,.5) -- (15.2,1);
\draw (15.2,.15) circle (.35cm) node {\LARGE $i_9$};
\draw[-,line width=.1mm] (-2.5,3) -- (-2,3);
\draw (-2.85,3) circle (.35cm) node {\LARGE $i_1$};
\draw[-,line width=.1mm] (5.5,3) -- (6,3);
\draw (5.15,3) circle (.35cm) node {\LARGE $i_4$};
\draw[-,line width=.1mm] (13.5,3) -- (14,3);
\draw (13.15,3) circle (.35cm) node {\LARGE $i_7$};
\draw[-,line width=.1mm] (2.5,3) -- (2,3);
\draw (2.85,3) circle (.35cm) node {\LARGE $i_2$};
\draw[-,line width=.1mm] (10.5,3) -- (10,3);
\draw (10.85,3) circle (.35cm) node {\LARGE $i_5$};
\draw[-,line width=.1mm] (18.5,3) -- (18,3);
\draw (18.85,3) circle (.35cm) node {\LARGE $i_8$};
\draw[-,line width=.1mm] (0,.5) -- (0,1);
\draw (0,.15) circle (.35cm) node {\LARGE $j_3$};
\draw[-,line width=.1mm] (8,.5) -- (8,1);
\draw (8,.15) circle (.35cm) node {\LARGE $j_6$};
\draw[-,line width=.1mm] (16,.5) -- (16,1);
\draw (16,.15) circle (.35cm) node {\LARGE $j_9$};
	\draw[dashed,draw=black] (-3.35,-.5) rectangle (3.35,4.6) node at (0, 5) {\huge $\text{Area}\ 1:\, a_1$};
	\draw[dashed,draw=black] (4.65,-.5) rectangle (11.35,4.6) node at (8, 5) {\huge $\text{Area}\ 2:\, a_2$};
	\draw[dashed,draw=black] (12.65,-.5) rectangle (19.35,4.6) node at (16, 5) {\huge $\text{Area}\ 3:\, a_3$};
	\draw[dashed,draw=black] (3.7,-0.5) rectangle (4.3,4.6) node at (4, 5) {\LARGE $\text{Link}\ e_1$};
	\draw[dashed,draw=black] (11.7,-0.5) rectangle (12.3,4.6) node at (12, 5) {\LARGE $\text{Link}\ e_2$};
	\end{tikzpicture}
	\caption{Nine-node three-area interconnected power system}\label{fig:three_area}
\end{figure*}
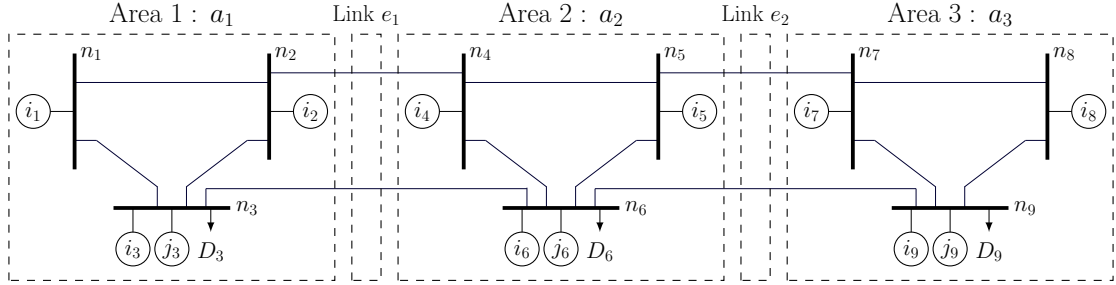

We describe a base model, which will be subject to several modifications in the system configuration and penetration of stochastic renewables to discuss the resulting changes in the benefit allocations described in Sections~\ref{sec:3} and~\ref{sec:4}.
We consider the nine-bus system depicted in Figure~\ref{fig:three_area} which comprises three areas. The intra-area transmission network consists of AC lines with capacity and reactance equal to $100$ MW and $0.13$ p.u., respectively. The four tie lines between areas~$1$ and~$2$, and between areas~$2$ and~$3$ are AC lines with capacity of $20$ MW, and reactance of $0.13$ p.u. each. 
\looseness=-1

The day-ahead price offers and the generation capacities of conventional units are provided in Table~\ref{tab:three_generators}. Units $i_1$, $i_4$, and $i_7$ are inflexible, i.e., these units cannot change their generation level during real-time operation, while all remaining units are flexible offering half of their capacity for upward and downward reserves provision at a cost equal to $10\%$ of their day-ahead energy offer $C$. The cost of load shedding $C^\text{sh}$ is equal to $1000$\euro$/$MWh for the inelastic electricity demands $D_3=220$ MW, $D_6=190$ MW, and $D_9=220$ MW. In addition, there are three wind power plants, $j_3,$ $j_6,$ and $j_9$, with installed capacities $50$, $80$, and $50$ MW, respectively. The stochastic wind power generation is modeled using two scenarios, $s_1$ and $s_2$, listed in Table~\ref{tab:three_wind} with probability of occurance $0.6$ and $0.4$, respectively. Hence, the corresponding expected wind power production $\overline{W}_j$ for $j_3$ is equal to $42$ MW, for $j_6$ is equal to $70.4$ MW, and for $j_9$ is equal to $42$ MW. Wind power price offers and subsequently the wind power spillage costs are considered to be zero.

\begin{table}[h]\footnotesize
\begin{minipage}{.6\linewidth}
\centering
\caption{Generator data}
\resizebox{.95\textwidth}{!}{\begin{tabular}{|c||c|c|c|c|c|c|c|c|c|}
\hline
 Unit &$i_1$ &$i_2$&$i_3$&$i_4$&$i_5$&$i_6$&$i_7$&$i_8$&$i_9$ \\ \hline
 $C$ (\euro$/$MWh) & $20$ &$30$ &$40$ &$30$ &$40$ &$50$ &$25$ &$35$ &$45$ \\ \hline
 $P^\text{cap}$ (MW) & $120$ &$50$ &$50$ &$120$ &$50$ &$50$ &$120$ &$50$ &$50$ \\ \hline
 Flexible & No & Yes & Yes & No & Yes & Yes & No & Yes & Yes \\ \hline
\end{tabular}}\label{tab:three_generators}
\end{minipage}
\begin{minipage}{.35\linewidth}
\centering
\caption{Wind scenarios as percentage of the nominal value of the power plant }
\resizebox{.91\textwidth}{!}{\begin{tabular}{|c||c|c|c|}
\hline
 Wind power plant & $j_3$ &$j_6$ & $j_9$ \\ \hline
 Scenario $s_1$&$1$ &$0.8$&$1$ \\ \hline
 Scenario $s_2$&$0.6$ &$1$&$0.6$ \\ \hline
\end{tabular}}\label{tab:three_wind}
\end{minipage}
\end{table}

Following the prevailing approach in which regional capacity markets are cleared separately, we set the percentage of transmission capacity allocated to reserves exchange equal to $\chi=0$. Reserve requirements for each area are calculated based on the probabilistic forecasts, such that the largest negative and positive deviations from the expected wind power production foreseen in the scenario set are covered by domestic resources. For instance, the upward and downward area reserve requirements for area $1$ are calculated~as, $RR^+_{a_1}  = \overline{W}_{j_3} - \min\{W_{j_3s_1},W_{j_3s_2}\}= 12\, \text{MW},$ and $RR^-_{a_1}  = \max\{W_{j_3s_1},W_{j_3s_2}\}-\overline{W}_{j_3} = 8\, \text{MW}.$
For the other areas, these values are given as $RR^+_{a_2}  =6.4\, \text{MW}$, $RR^-_{a_2}  =9.6\, \text{MW}$, $RR^+_{a_3}  =12\, \text{MW}$, $RR^-_{a_3}  =8\, \text{MW}$.

The market costs and transmission allocations resulting from the preemptive model are provided in Table~\ref{tab:three_costs}. The preemptive model reallocated transmission resources from the day-ahead energy trading to the reserve capacity trading, increasing the costs in the day-ahead market. This reallocation yields an expected system cost of $13{,}238.0$\euro, which translates to $25.9\%$ reduction compared to the cost of $17{,}871.2$\euro\ from the existing sequential market. Under the existing setup with $\chi=0$, the uncertainty realization $s_2$ leads to significant load shedding in the balancing stage. In this scenario, even though we have enough reserve capacity, we are not able to deploy it due to network congestion. This problem is avoided by enabling reserve exchange when the preemptive model is implemented. Quantities assigned to each generator at all trading floors are provided in~\ref{app:market_outcome}. 
\looseness=-1

We now provide a budget balanced cost allocation method for the existing sequential market. For this method, we assume that all three trading floors are cleared by marginal pricing mechanisms (zonal prices for the reserve capacities, nodal prices for the day-ahead and balancing energy services), albeit, similar methods can be applied also to other payment mechanisms. 
This method assigns producer and consumer surpluses, and congestion rents of the intra-area lines to their corresponding areas, and divides the congestion rents of the tie lines equally between the adjacent areas, see~\citep{kristiansen2018mechanism}. Budget balance holds since the market cost is given by the opposite of the sum of producer and consumer surpluses, and congestion rent for each trading floor.  These values are summarized in~\ref{app:market_outcome}. We refer to Table~\ref{tab:budget_out} for the resulting cost allocations. Area~1 is allocated a large cost in scenario $s_2$ because of the load shedding in node~3.
\looseness=-1

\begin{table}[h]\footnotesize
\begin{minipage}{.64\linewidth}
\centering
\caption{Comparison of market costs (in \euro)}
\resizebox{.9\textwidth}{!}{\begin{tabular}{|c||c|c|}
\hline
 Model & Existing Seq. Market & Preemptive Model \\ \hline
 $[\chi_{e_1},\,\chi_{e_2}]$ & $[0,\,0]$ &$[0,\,0.0592]$\\ \hline 
 Reserve capacity cost &$194.0$ &$191.6$\\ \hline 
 Day-ahead cost &$13{,}087.2$ &$13{,}120.2$\\ \hline 
 Balancing cost in $s_1$ &$1{,}150.0$ &$-410.7$\\ \hline 
 Balancing cost in $s_2$ &$9{,}750.0$ &$431.5$\\ \hline 
 Total cost in $s_1$ &$14{,}431.2$ &$12{,}901.2$\\ \hline 
 Total cost in $s_2$ &$23{,}031.2$ &$13{,}743.4$\\ \hline 
\end{tabular}}\label{tab:three_costs}
\end{minipage}
\begin{minipage}{.35\linewidth}
\centering
\caption{Cost allocation for each area in the existing sequential market (in \euro)}
\resizebox{.95\textwidth}{!}{\begin{tabular}{|c||c|c|c|}
\hline
 Areas & Area $1$  & Area $2$ & Area $3$\\ \hline
 $J^{s_1}_a(\emptyset)$& $4{,}348.4$  & $9{,}853.8$ & $229.0$\\ \hline
 $J^{s_2}_a(\emptyset)$& $16{,}348.4$  & $3{,}453.8$ & $3{,}229.0$\\ \hline
\end{tabular}}\label{tab:budget_out}
\end{minipage}\
\end{table}

%Next, we are ready to bring in benefit allocation mechanisms for this problem.

\subsubsection{Comparison of the different benefit allocations}\label{sec:shapno}

\begin{figure}[h]
\begin{minipage}[t]{0.49\linewidth}
	\centering
	\begin{tikzpicture}[scale=1, every node/.style={scale=0.5}]
			\node at (0,0) {\includegraphics[width=2.2\textwidth]{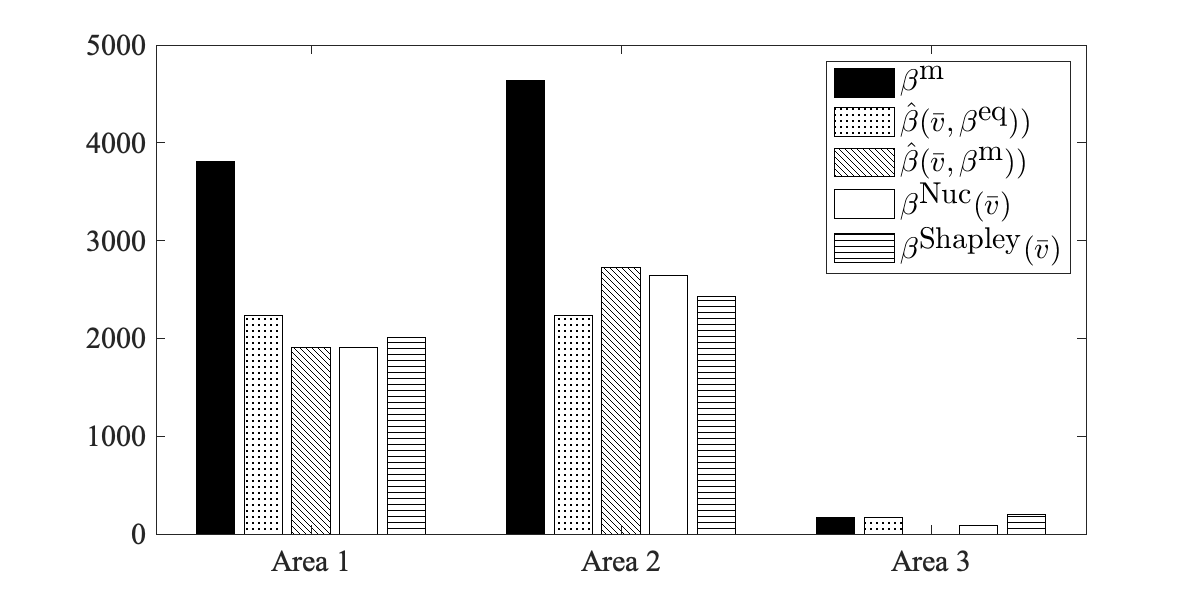}};
	\end{tikzpicture}
	\caption{Benefit allocations for the three-area system (in \euro)}\label{fig:three_benefit}
\end{minipage}
\begin{minipage}[t]{0.49\linewidth}
	\centering
	\begin{tikzpicture}[scale=1, every node/.style={scale=0.5}]
			\node at (0,0) {\includegraphics[width=2.2\textwidth]{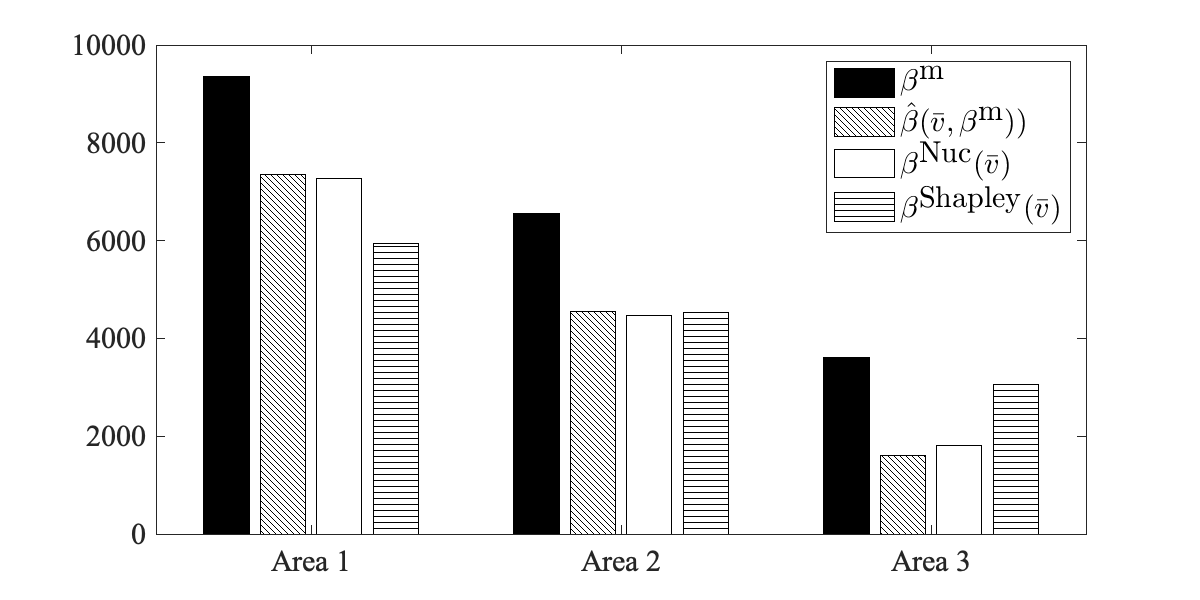}};
	\end{tikzpicture}
	\caption{Benefit allocations after connecting areas 1 and 3 (in \euro)}\label{fig:three_connected}
\end{minipage}
\end{figure}

Benefit allocation mechanisms for the expected cost reduction are provided in Figure~\ref{fig:three_benefit}. The core is nonempty since area 2 satisfies the veto condition in Proposition~\ref{prop:corempt}. Marginal contribution benefit allocation $\beta^\text{m}$ is not in the core since it is not efficient. We provide this allocation since it can be regarded as a fair outcome. Observe that the marginal contributions of areas 1 and 2 are larger than that of area 3. This is because area 1 has low cost generators and area 2 is indispensable for any coordination considering that in the current network configuration, in which areas 1 and 3 are not directly interconnected, area 2 has to act as an intermediary for any reserves exchange.

The Shapley value $\beta^\text{Shapley}$ is not in the core. Among the core constraints, combining $\sum_{a\in\AC}\beta_a=v(\AC)$ with $\sum_{a\in\AC\setminus\hat{a}}\beta_a\geq v(\AC\setminus\hat{a})$ implies that $\beta_{\hat a}\leq\beta^\text{m}_{\hat a}=v(\AC)-v(\AC\setminus\hat{a})$, or equivalently, no area can receive more than its marginal contribution in the core. This condition is violated for the Shapley value assigned to area 3. The coalitional value function is also not supermodular, since
$\bar v(\{1,2,3\})-\bar v(\{1,2\})\not\geq \bar v(\{2,3\})-\bar v(\{2\})\implies
    4{,}633.1-4{,}460.5 \not\geq 826.8-0.$ On the other hand, the nucleolus $\beta^\text{Nuc}$ is in the core, however, the lexicographic minimization results in allocating benefits to area~3. We later see that there is a core allocation that better approximates the marginal contribution in terms of minimizing the Euclidean distance by allocating no benefits to~area~3.

Finally, we employ our approach approximating two different criteria, i.e., marginal contribution and equal shares, with corresponding allocations being denoted as $\hat{\beta}(\bar{v},\beta^{\text{m}}))$ and $\hat{\beta}(\bar{v},\beta^{\text{eq}}))$, respectively. These two outcomes are different from each other, and they approximate their respective fairness consideration in an effective manner. This criterion should be decided either by the regulator or it should be based on the consensus of participating areas. In the following, we will approximate the marginal contribution, since similar discussions can be made for any other criteria.
\looseness=-1

Next, we study the budget balance per scenario for the cost allocation in the preemptive model. For all efficient benefit allocations of the expected cost reduction, i.e., all methods except the marginal contribution allocation, the budget  $\sum_{a\in\AC}J^s_a(\AC)-J^{s}(\AC)$ remains unchanged. In scenario $s_1$, there is a deficit of $3{,}103.1$\euro, whereas in scenario $s_2$ there is a surplus of $4{,}654.7$\euro, thus budget balance is obtained in expectation. 
In the coalitional game arising from the scenario-specific cost variation, despite that $K_\text{Core}(\bar v)$ is nonempty, the core $K_\text{Core}(v^{s_1})$ is empty, since the condition in Proposition~\ref{prop:empty} is satisfied by $v^{s_1}(\{1,2\})=J^{s_1}(\emptyset) - J^{s_1}(\{1,2\}) > J^{s_1}(\emptyset) - J^{s_1}(\{1,2,3\})=v^{s_1}(\{1,2,3\})\!\!\! \implies \!\!\!
    14{,}431.2-12{,}884.6 > 14{,}431.2-12{,}901.2.$
For scenario $s_2$, this condition is not found and $K_\text{Core}(v^{s_2})$ is nonempty, since the coalitional game is supermodular.

To address the budget balance, we now employ {the proposed  scenario-specific least-core selecting mechanism.} The scenario-specific allocations generated by $\hat{\beta}(\bar v,\beta^{\text{m}})$ for the expected cost reduction are given by $\beta(v^{s_1},\hat{\beta}(\bar v,\beta^{\text{m}}))=[628.5,\, 901.5,\, 0]^\top$ and $\beta(v^{s_2},\hat{\beta}(\bar v,\beta^{\text{m}}))=[3{,}815.2,\, 5{,}472.6,\, 0]^\top$. These allocations result in a budget balanced cost allocation under both scenarios, since they sum up to the scenario-specific cost variations in Table~\ref{tab:three_costs}. In the remainder, we will focus our efforts on the coalitional game arising from the expected cost reduction, since we can always map the benefit allocations to the scenario-specific case using our proposed approach. 
\looseness=-1

\subsubsection{Impact of the uncertainties on benefit allocations}\label{sec:threearea_windprof}
Here, we aim to assess the impact of the spatial correlation of the wind power forecast errors on the outcome of the different benefit allocation mechanisms that we consider in this work. In order to eliminate the impact of the network topology (cf. \ref{app:three_area_netw_top}), we connect areas 1 and 3 via two AC lines. The first connects nodes 1 and 8, and the second connects nodes 3 and 9, each with transmission capacity of $20$ MW, and reactance of $0.13$ p.u. The area graph is not a star anymore, and Proposition~\ref{prop:corempt} is not applicable. However, we verified that the core is still nonempty. 

The resulting benefit allocations for the expected cost reduction are provided in Figure~\ref{fig:three_connected}, which shows that areas~1 and~2 receive most of the benefit under every allocation mechanism. This outcome can be explained considering that these areas have complementary wind power production scenarios, i.e., the corresponding wind power scenarios exhibit negative correlation. Moreover, area~1 has low cost generation.
Finally, notice that the total benefits are greater than the ones from the example in Section~\ref{sec:shapno}. This follows since compared to Section~\ref{sec:shapno} expected system cost is increased by $49$\% ($26{,}687.9$\euro) in the existing sequential market due to additional network dependencies, whereas this cost is decreased by $0.01$\% ($13{,}161.0$\euro) in the preemptive model.
%Finally, the Shapley value is in the core since the coalition value function is 

\subsubsection{Benefit allocations in the case of an empty core}\label{sec:emptycore}
A natural question that arises in the context of this work is how the different benefit allocation mechanisms perform when we have an empty core, which can occur when the condition in Proposition~\ref{prop:corempt} is not found.
To this end, we modify the example in Section~\ref{sec:threearea_windprof} by changing the wind scenarios. The stochastic wind power generation is modeled using two scenarios, $s_1$ and $s_2$ with probability of occurance $0.8$ and $0.2$. We have $1$ and $0.8$ for $j_3$, $0.4$ and $1$ for $j_6$, $0.4$ and $1$ for $j_9$ as the percentages of the nominal values of the plants, respectively. Hence, the corresponding expected wind power productions for $j_3$ is equal to $48$ MW, for $j_6$ is equal to $41.6$ MW, and for $j_9$ is equal to $26$ MW. The reserve requirements are recomputed accordingly. Since the uncertainty is significantly increased, we allow the units $i_1$, $i_4$, and $i_7$ to be flexible in order to ensure feasibility.

\begin{figure}[h]
	\centering
	\begin{tikzpicture}[scale=1, every node/.style={scale=0.48}]
			\node at (0,0) {\includegraphics[width=1.2\textwidth]{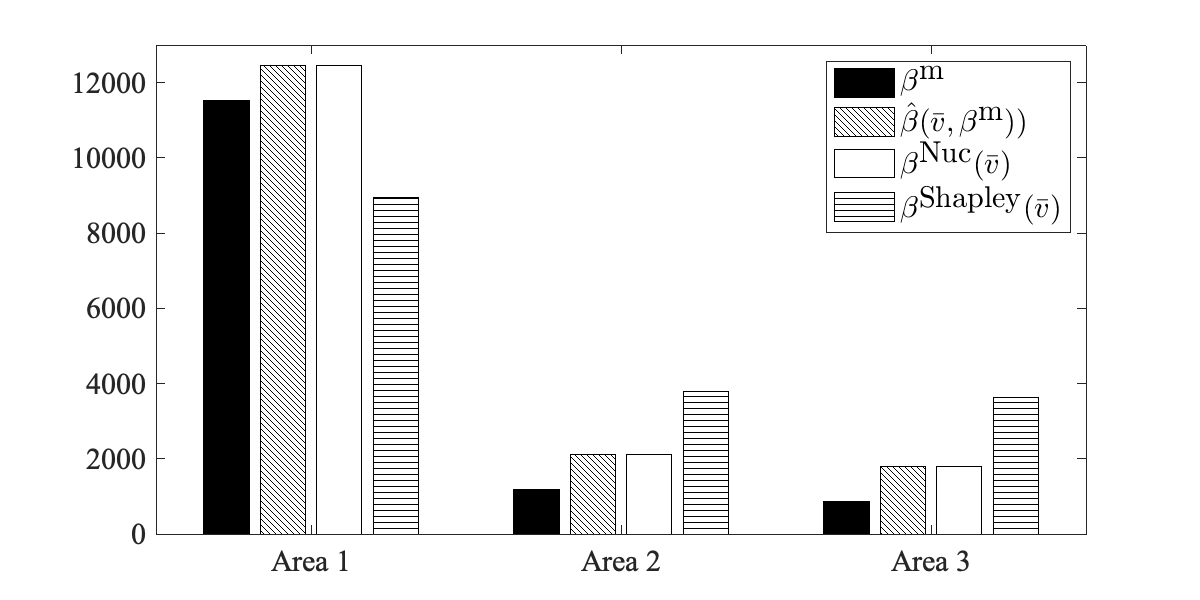}};
	\end{tikzpicture}
	\caption{Benefit allocations in case the core is empty (in \euro)}\label{fig:three_empty}
\end{figure}

The resulting benefit allocations for the expected cost reduction are provided in Figure~\ref{fig:three_empty}.
We observe that the nucleolus and the least-core selecting benefit allocation coincide. For both allocations, the maximum violation of a stability constraint is given by $\epsilon^*=924.9$\euro, where $\epsilon^*(\bar v)=\min\{\epsilon\,|\,K_\text{Core}(\bar v,\epsilon)\not=\emptyset\}$.
On the other hand, the maximum stability violation for the Shapley value is $2{,}752.0$\euro. In other words, if the Shapley value is utilized, there are $3$ times the profits to be made by not participating in the preemptive model compared to the case implementing a least-core allocation. We see that all benefit allocation mechanisms allocated the most benefits to area 1, since it has low cost generation and also its wind profile complements the wind profiles of areas 2 and 3.
\looseness=-1

\subsection{Case study based on the IEEE RTS}\label{sec:IEEE_RTS}
We now consider a six-area power system that is based on the modernized version of the IEEE Reliability Test System (RTS) presented in \citep{pandzicunit}. 
The definitions of the areas correspond to the ones proposed by~\citet{dvorkin2018setting} and \citet{jensen2017cost}, and they are provided in~\ref{app:IEEERTS96_area}.  The nodal positions, types, generation capacities and offers from conventional and wind power generators, and transmission line parameters are provided in~\citep{karacadata}. Due to their limited flexibility, nuclear, coal and integrated gasification combined cycle (IGCC) units do not provide any reserves. On the other hand, open and combined cycle gas turbines (OCGT and CCGT) offer $50\%$ of their capacity for upward and downward reserves at a cost equal to $20\%$ of their day-ahead energy offer. Wind power production is modeled using a set of 10 equiprobable scenarios that capture the spatial correlation of forecast errors over the different wind farm locations provided in~\citep{bukhsh2017data}. The demand is inelastic with the cost of load shedding equal to $1{,}000$\euro$/$MWh. In the existing sequential market, the percentage of transmission capacity allocated to reserves exchange is set to $\chi=0$, while the reserve requirements are calculated according to the methodology discussed in Section~\ref{sec:illu_ex}. 

Table~\ref{tab:rts_costs} compares the system costs and transmission allocations resulting from the existing market with $\chi$ fixed to zero and the preemptive model where $\chi$ is a decision variable. The preemptive model yields an expected cost of $87{,}594.8$\euro, which translates to $2.3\%$ reduction compared to the cost of $89{,}696.4$\euro\ from the existing market.
This can be explained by $2{,}097.6$\euro\ reduction in the expected balancing cost obtained by eliminating the need for load shedding. Using the approach in Section~\ref{sec:illu_ex}, we provide the expected values for a budget balanced cost allocation for the existing sequential market in Table~\ref{tab:rts_budget_out}.

\begin{table}[h]\footnotesize
\centering\vspace{-.1cm}
\caption{Comparison of market costs (in \euro)}
\resizebox{.66\textwidth}{!}{\begin{tabular}{|c||c|c|}
\hline
 Model & Existing Seq. Market & Preemptive Model \\ \hline
 $\chi=[\chi_{e_1},\ldots,\chi_{e_7}]$ & $[0,\,0,\,0,\,0,\,0,\,0,\,0]$ &$[0,\,0.359,\,0,\,0.038,\,0,\,0,\,0]$\\ \hline 
 Reserve capacity cost & $2{,}392.5$ & $2{,}389.4$\\ \hline 
 Day-ahead cost & $90{,}734.2$& $90{,}733.3$\\ \hline 
 {Expected} balancing cost  &$-3{,}430.3$& $-5{,}527.9$\\ \hline 
 {Expected} cost & $89{,}696.4$ & $87{,}594.8$\\ \hline 
\end{tabular}}\label{tab:rts_costs}
\end{table}
\begin{table}[h]\footnotesize
\centering
\caption{Expected cost allocation for each area in the existing sequential market (in \euro)}
\resizebox{.66\textwidth}{!}{\begin{tabular}{|c||c|c|c|c|c|c|}
\hline
 Areas & Area $1$  & Area $2$ & Area $3$ &Area $4$  & Area $5$ & Area $6$\\ \hline
 $\mathbb{E}_s[J^{s}_a(\emptyset)]$ &$8{,}966.1$ & $25{,}252.6$    & $10{,}085.2$ & $10{,}208.2$ &  $25{,}151.1$ &  $10{,}033.2$\\ \hline
\end{tabular}}\label{tab:rts_budget_out}
\end{table}

\begin{figure}[h]
\begin{minipage}[t]{0.48\linewidth}
	\centering
	\begin{tikzpicture}[scale=1, every node/.style={scale=0.48}]
			\node at (0,0) {\includegraphics[width=2.2\textwidth]{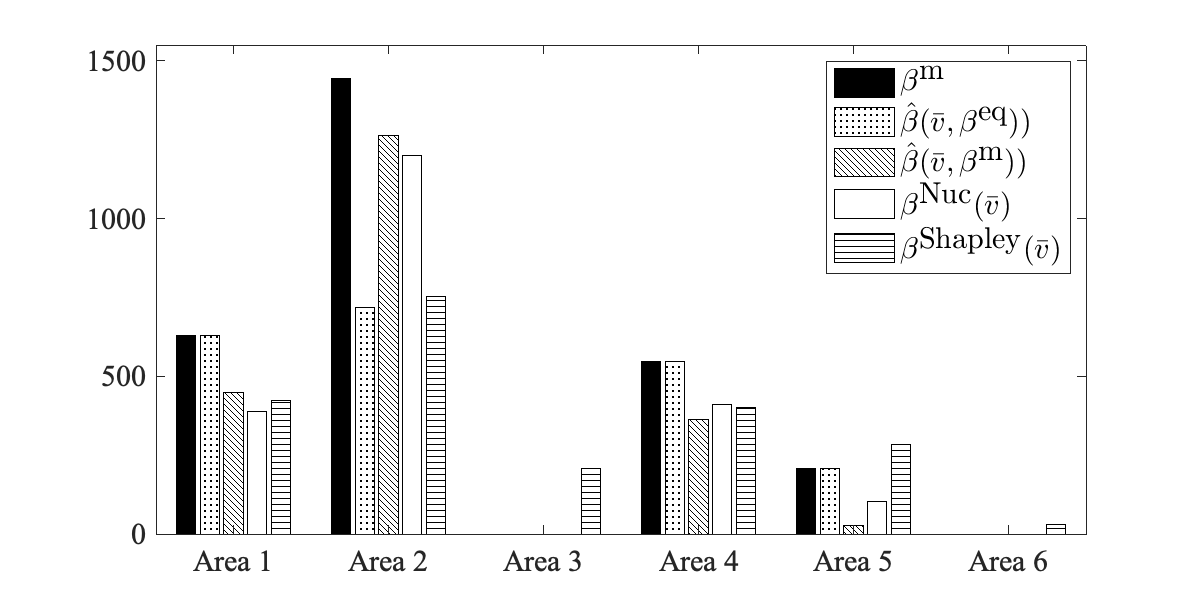}};
	\end{tikzpicture}
	\caption{Benefit allocations for the IEEE RTS case study (in \euro)}\label{fig:rts_benefit}
\end{minipage}
\begin{minipage}[t]{0.48\linewidth}
    \centering
	\begin{tikzpicture}[scale=1, every node/.style={scale=0.5}]
			\node at (0,0) {\includegraphics[width=1.65\textwidth]{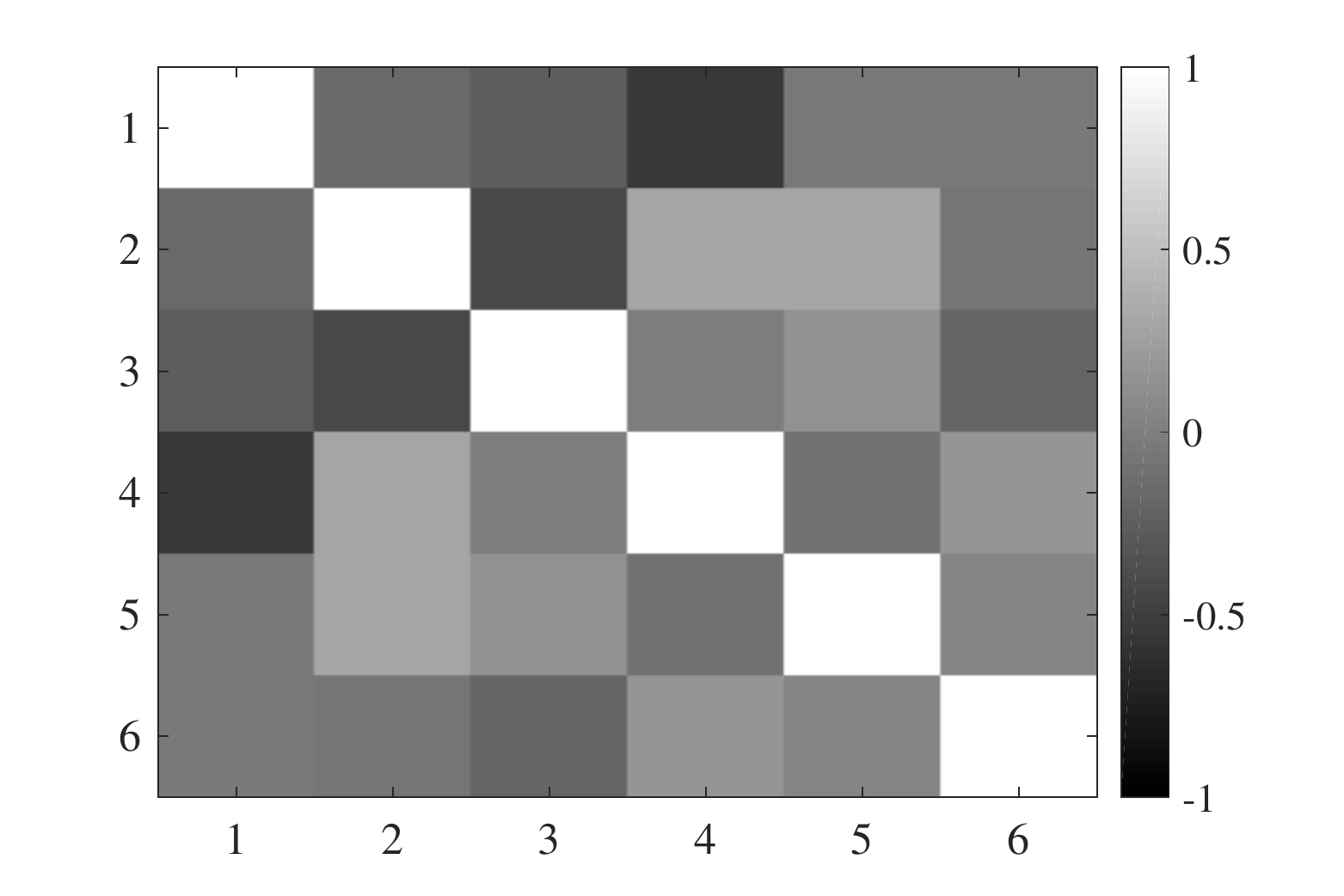}};
	\end{tikzpicture}
	\caption{The Pearson correlation coefficients for the wind profiles of all areas}\label{fig:rts_corr}
\end{minipage}
\end{figure}

The results of the different benefit allocation mechanisms for the expected cost reduction are provided in Figure~\ref{fig:rts_benefit}. 
We verified that the core is nonempty even though the condition in Proposition~\ref{prop:corempt} is not satisfied. Marginal contribution benefit allocation is not in the core since it is not efficient, whereas the Shapley value is not in the core since areas 3, 5, and~6 receive more than their marginal contributions. The nucleolus and the least-core selecting benefit allocation mechanisms result in core allocations. Notice that our approach provides a different benefit allocation depending on the criteria considered. The nucleolus allocation is not consistent with the marginal contribution allocation since it allocates more benefits to area~4 compared to area~1. All mechanisms allocated the most benefits to area 2, since it has a central role by being well-connected in the area graph. On the other hand, areas 1 and 4 are also allocated a significant amount, since they are the two largest areas with wind profiles complementing each other as it is shown in Figure~\ref{fig:rts_corr}.

% To see this complementarity, we provide the Pearson correlation coefficients as a heat map for the wind profiles of all areas in Figure~\ref{fig:rts_corr}.

We now provide the computational comparison for the different benefit allocation mechanisms. The coalitions $J(\emptyset)$ and $J(\AC)$ are precomputed to obtain $v(\AC)$, in $19.4$ and $35.4$ seconds, respectively. 
The calculation of the marginal contribution allocation, which involves solving the preemptive model~\eqref{mod:B-Pre} for coalitions $\{\AC\setminus\{a\}\}_{a\in\AC}$, requires $119.7$ seconds. 
The Shapley value requires solving the preemptive model~\eqref{mod:B-Pre} for all coalitions except the singleton sets, the empty set, and the full set, i.e., $2^6-6-1-1=58$ coalitions, and the resulting computational time is $1{,}264.6$ seconds. 
The least-core selecting mechanism with the marginal contribution criteria requires only a single iteration from the constraint generation algorithm, which takes $84.8$ seconds. This constraint generation algorithm converges fast, since the initial family of coalitions $\mathcal{F}^1$ is not an empty set; instead, it is given by the coalitions that were used to compute the marginal contribution allocation. 
On the other hand, the least-core selecting mechanism with the equal shares criteria requires four iterations from the constraint generation algorithm, which takes $150.4$ seconds. Notice that, in this case, the initial family of coalitions $\mathcal{F}^1$ is an empty set. 
Finally, using the method proposed in~\citep{kopelowitz1967computation,fromen1997reducing}, the nucleolus is computed by solving $15$ linear programs sequentially to find $21$ coalitional equality constraints that fully describe the nucleolus allocation. This computation takes $1{,}266.2$ seconds, since it needs the complete list of coalition values.
As an alternative, the iterative method in~\citep{hallefjord1995computing} would require running $15$ separate constraint generation algorithms, {increasing significantly the computational time compared to} the least-core selecting mechanism, since each algorithm run requires at least one iteration of constraint generation. In~\ref{app:modif_rts}, we provide modifications to the IEEE RTS case study verifying our observations from the illustrative three-area nine-node system.

\section{Conclusion}\label{sec:conc}

In this work, we formulated a coalition-dependent version of the preemptive transmission allocation model that defines the optimal inter-area transmission capacity allocation {between energy and reserves} for a given {set of areas participating in the coalition}. We then accompanied this general model with benefit allocation mechanisms such that all {coalition member areas} have sufficient incentives to accept the solution proposed by the preemptive model. We formulated the coalitional game both as an ex-ante and as an ex-post process with respect to the uncertainty realization {and we} showed that the former results in budget balance in expectation, whereas the latter results in budget balance in every {uncertainty realization}. Applying the prevailing benefit allocation mechanisms to a larger case study, we showed that they are unable to find a benefit allocation with minimal stability violation within a reasonable computational timeframe. To address this issue for both coalitional games under consideration, we proposed the least-core selecting benefit allocation mechanism and we formulated an iterative constraint generation algorithm for the efficient computation of this allocation. 
Considering that this work aims to contribute to the ongoing discussion towards the design of the transmission allocation model, our benefit allocation mechanism can be adapted to different plausible fairness criteria that may be imposed by the regulatory authorities, moving towards the full integration of the European balancing markets.

Our future work will explore the development of decentralized schemes that enable the coordination of areas in terms of transmission allocation, while preserving privacy of the areas with minimal exchange of intra-area information. 
% As an extension, we aim to study the impact of strategic behaviour of the areas and the generators on the transmission allocations and the resulting benefit allocation mechanisms. 
Finally, analyzing the impact of out-of-sample uncertainty realizations on the benefit allocation mechanisms would be interesting.

%\section*{References}

%\bibliographystyle{elsarticle-harv}
\bibliography{library}	

\begin{thebibliography}{57}
\expandafter\ifx\csname natexlab\endcsname\relax\def\natexlab#1{#1}\fi
\providecommand{\url}[1]{\texttt{#1}}
\providecommand{\href}[2]{#2}
\providecommand{\path}[1]{#1}
\providecommand{\DOIprefix}{doi:}
\providecommand{\ArXivprefix}{arXiv:}
\providecommand{\URLprefix}{URL: }
\providecommand{\Pubmedprefix}{pmid:}
\providecommand{\doi}[1]{\href{http://dx.doi.org/#1}{\path{#1}}}
\providecommand{\Pubmed}[1]{\href{pmid:#1}{\path{#1}}}
\providecommand{\bibinfo}[2]{#2}
\ifx\xfnm\relax \def\xfnm[#1]{\unskip,\space#1}\fi
%Type = Misc
\bibitem[{{Agency for the Cooperation of Energy Regulators
  (ACER)}(2012)}]{ACER}
\bibinfo{author}{{Agency for the Cooperation of Energy Regulators (ACER)}}
  (\bibinfo{year}{2012}).
\newblock \bibinfo{title}{Framework guidelines on electricity balancing}.
\newblock \URLprefix
  \url{http://www.acer.europa.eu/official_documents/acts_of_the_agency/framework_guidelines/framework%20guidelines/framework%20guidelines%20on%20electricity%20balancing.pdf}.
%Type = Article
\bibitem[{Ahlstrom et~al.(2015)Ahlstrom, Ela, Riesz, O'Sullivan, Hobbs,
  O'Malley, Milligan, Sotkiewicz \& Caldwell}]{ahlstrom2015evolution}
\bibinfo{author}{Ahlstrom, M.}, \bibinfo{author}{Ela, E.},
  \bibinfo{author}{Riesz, J.}, \bibinfo{author}{O'Sullivan, J.},
  \bibinfo{author}{Hobbs, B.~F.}, \bibinfo{author}{O'Malley, M.},
  \bibinfo{author}{Milligan, M.}, \bibinfo{author}{Sotkiewicz, P.}, \&
  \bibinfo{author}{Caldwell, J.} (\bibinfo{year}{2015}).
\newblock \bibinfo{title}{The evolution of the market: Designing a market for
  high levels of variable generation}.
\newblock {\it \bibinfo{journal}{IEEE Pow. and En. Mag.}\/},  {\it
  \bibinfo{volume}{13}\/}, \bibinfo{pages}{60--66}.
%Type = Article
\bibitem[{Aumann \& Maschler(1985)}]{aumann1985game}
\bibinfo{author}{Aumann, R.~J.}, \& \bibinfo{author}{Maschler, M.}
  (\bibinfo{year}{1985}).
\newblock \bibinfo{title}{Game theoretic analysis of a bankruptcy problem from
  the {Talmud}}.
\newblock {\it \bibinfo{journal}{J. of {E}con. {T}h.}\/},  {\it
  \bibinfo{volume}{36}\/}, \bibinfo{pages}{195--213}.
%Type = Inproceedings
\bibitem[{Avramiotis-Falireas et~al.(2018)Avramiotis-Falireas, Margelou \&
  Zima}]{avramiotis2018investigations}
\bibinfo{author}{Avramiotis-Falireas, I.}, \bibinfo{author}{Margelou, S.}, \&
  \bibinfo{author}{Zima, M.} (\bibinfo{year}{2018}).
\newblock \bibinfo{title}{Investigations on a fair {TSO-TSO} settlement for the
  imbalance netting process in european power system}.
\newblock In {\it \bibinfo{booktitle}{15th Int. Conf. on the EEM}\/} (pp.
  \bibinfo{pages}{1--6}).
\newblock \bibinfo{organization}{IEEE}.
%Type = Article
\bibitem[{Baeyens et~al.(2013)Baeyens, Bitar, Khargonekar \&
  Poolla}]{baeyens2013coalitional}
\bibinfo{author}{Baeyens, E.}, \bibinfo{author}{Bitar, E.~Y.},
  \bibinfo{author}{Khargonekar, P.~P.}, \& \bibinfo{author}{Poolla, K.}
  (\bibinfo{year}{2013}).
\newblock \bibinfo{title}{Coalitional aggregation of wind power}.
\newblock {\it \bibinfo{journal}{IEEE Trans. on Pow. Syst.}\/},  {\it
  \bibinfo{volume}{28}\/}, \bibinfo{pages}{3774--3784}.
%Type = Book
\bibitem[{Birge \& Louveaux(2011)}]{birge2011introduction}
\bibinfo{author}{Birge, J.~R.}, \& \bibinfo{author}{Louveaux, F.}
  (\bibinfo{year}{2011}).
\newblock {\it \bibinfo{title}{Introduction to stochastic programming}\/}.
\newblock \bibinfo{publisher}{Springer Sci. \& B. Med.}
%Type = Techreport
\bibitem[{Bondy et~al.(2014)Bondy, Tarnowski, Heussen \&
  Hansen}]{bondy2014operational}
\bibinfo{author}{Bondy, D. E.~M.}, \bibinfo{author}{Tarnowski, G.},
  \bibinfo{author}{Heussen, K.}, \& \bibinfo{author}{Hansen, L.~H.}
  (\bibinfo{year}{2014}).
\newblock {\it \bibinfo{title}{Operational scenario: Manual regulating
  power}\/}.
\newblock \bibinfo{type}{Technical Report} iPower Consortium.
%Type = Article
\bibitem[{Bouffard et~al.(2005)Bouffard, Galiana \&
  Conejo}]{bouffard2005market}
\bibinfo{author}{Bouffard, F.}, \bibinfo{author}{Galiana, F.~D.}, \&
  \bibinfo{author}{Conejo, A.~J.} (\bibinfo{year}{2005}).
\newblock \bibinfo{title}{Market-clearing with stochastic security-part {I}:
  {F}ormulation}.
\newblock {\it \bibinfo{journal}{IEEE Trans. on Pow. Syst.}\/},  {\it
  \bibinfo{volume}{20}\/}, \bibinfo{pages}{1818--1826}.
%Type = Misc
\bibitem[{Bukhsh(2017)}]{bukhsh2017data}
\bibinfo{author}{Bukhsh, W.} (\bibinfo{year}{2017}).
\newblock \bibinfo{title}{Data for stochastic multiperiod optimal power flow
  problem}.
\newblock
  \bibinfo{howpublished}{\url{https://sites.google.com/site/datasmopf/home}}.
%Type = Article
\bibitem[{Churkin et~al.(2019)Churkin, Pozo, Bialek, Korgin \&
  Sauma}]{churkin2019can}
\bibinfo{author}{Churkin, A.}, \bibinfo{author}{Pozo, D.},
  \bibinfo{author}{Bialek, J.}, \bibinfo{author}{Korgin, N.}, \&
  \bibinfo{author}{Sauma, E.} (\bibinfo{year}{2019}).
\newblock \bibinfo{title}{Can cross-border transmission expansion lead to fair
  and stable cooperation? {Northeast Asia} case analysis}.
\newblock {\it \bibinfo{journal}{Energy Econ.}\/},  (p.
  \bibinfo{pages}{104498}).
%Type = Article
\bibitem[{Cs{\'o}ka et~al.(2009)Cs{\'o}ka, Herings \&
  K{\'o}czy}]{csoka2009stable}
\bibinfo{author}{Cs{\'o}ka, P.}, \bibinfo{author}{Herings, P. J.-J.}, \&
  \bibinfo{author}{K{\'o}czy, L.~{\'A}.} (\bibinfo{year}{2009}).
\newblock \bibinfo{title}{Stable allocations of risk}.
\newblock {\it \bibinfo{journal}{Games and Econ. Beh.}\/},  {\it
  \bibinfo{volume}{67}\/}, \bibinfo{pages}{266--276}.
%Type = Article
\bibitem[{Day \& Raghavan(2007)}]{day2007fair}
\bibinfo{author}{Day, R.~W.}, \& \bibinfo{author}{Raghavan, S.}
  (\bibinfo{year}{2007}).
\newblock \bibinfo{title}{Fair payments for efficient allocations in public
  sector combinatorial auctions}.
\newblock {\it \bibinfo{journal}{Management Sci.}\/},  {\it
  \bibinfo{volume}{53}\/}, \bibinfo{pages}{1389--1406}.
%Type = Article
\bibitem[{Delikaraoglou \& Pinson(2018)}]{delikaraoglou2018optimal}
\bibinfo{author}{Delikaraoglou, S.}, \& \bibinfo{author}{Pinson, P.}
  (\bibinfo{year}{2018}).
\newblock \bibinfo{title}{Optimal allocation of {HVDC} interconnections for
  exchange of energy and reserve capacity services}.
\newblock {\it \bibinfo{journal}{Energy Syst.}\/},  (pp.
  \bibinfo{pages}{1--41}).
%Type = Article
\bibitem[{Deng \& Papadimitriou(1994)}]{deng1994complexity}
\bibinfo{author}{Deng, X.}, \& \bibinfo{author}{Papadimitriou, C.~H.}
  (\bibinfo{year}{1994}).
\newblock \bibinfo{title}{On the complexity of cooperative solution concepts}.
\newblock {\it \bibinfo{journal}{Math. of Oper. Res.}\/},  {\it
  \bibinfo{volume}{19}\/}, \bibinfo{pages}{257--266}.
%Type = Article
\bibitem[{Denholm \& Hand(2011)}]{denholm2011grid}
\bibinfo{author}{Denholm, P.}, \& \bibinfo{author}{Hand, M.}
  (\bibinfo{year}{2011}).
\newblock \bibinfo{title}{Grid flexibility and storage required to achieve very
  high penetration of variable renewable electricity}.
\newblock {\it \bibinfo{journal}{Energy Pol.}\/},  {\it
  \bibinfo{volume}{39}\/}, \bibinfo{pages}{1817--1830}.
%Type = Article
\bibitem[{Drechsel \& Kimms(2010)}]{drechsel2010computing}
\bibinfo{author}{Drechsel, J.}, \& \bibinfo{author}{Kimms, A.}
  (\bibinfo{year}{2010}).
\newblock \bibinfo{title}{Computing core allocations in cooperative games with
  an application to cooperative procurement}.
\newblock {\it \bibinfo{journal}{Int. J. of Production Econ.}\/},  {\it
  \bibinfo{volume}{128}\/}, \bibinfo{pages}{310--321}.
%Type = Article
\bibitem[{Dvorkin et~al.(2018)Dvorkin, Delikaraoglou \&
  Morales}]{dvorkin2018setting}
\bibinfo{author}{Dvorkin, V.}, \bibinfo{author}{Delikaraoglou, S.}, \&
  \bibinfo{author}{Morales, J.~M.} (\bibinfo{year}{2018}).
\newblock \bibinfo{title}{Setting reserve requirements to approximate the
  efficiency of the stochastic dispatch}.
\newblock {\it \bibinfo{journal}{IEEE Trans. on Pow. Syst.}\/},  {\it
  \bibinfo{volume}{34}\/}, \bibinfo{pages}{1524--1536}.
%Type = Misc
\bibitem[{{EC}(2016)}]{EC2}
\bibinfo{author}{{EC}} (\bibinfo{year}{2016}).
\newblock \bibinfo{title}{Integration of electricity balancing markets and
  regional procurement of balancing reserves.}
\newblock \URLprefix
  \url{https://ec.europa.eu/energy/sites/ener/files/documents/dg_ener_balancing_-_161021_-_final_report_-_version_27.pdf}.
%Type = Misc
\bibitem[{{EC}(2017)}]{EC1}
\bibinfo{author}{{EC}} (\bibinfo{year}{2017}).
\newblock \bibinfo{title}{{C}ommission {R}egulation ({EU}) 2017/2195
  {E}stablishing a guideline on electricity balancing.}
\newblock \URLprefix
  \url{https://eur-lex.europa.eu/legal-content/EN/TXT/PDF/?uri=CELEX:32017R2195&from=EN}.
%Type = Misc
\bibitem[{Energinet.dk(2010)}]{energinet}
\bibinfo{author}{Energinet.dk} (\bibinfo{year}{2010}).
\newblock \bibinfo{title}{Reservation of capacity of ancillary services -
  general principles and the sk4 case}.
\newblock
  \bibinfo{howpublished}{\url{http://www.elforsk.se/Documents/Market\%20Design/seminars/Interconnectors/4
  Energinetdk.pdf}}.
%Type = Misc
\bibitem[{{ENTSO-E}(2018{\natexlab{a}})}]{ENTSOE1}
\bibinfo{author}{{ENTSO-E}} (\bibinfo{year}{2018}{\natexlab{a}}).
\newblock \bibinfo{title}{An overview of the european balancing market and
  electricity balancing guideline}.
\newblock \URLprefix
  \url{https://docstore.entsoe.eu/Documents/Network\%20codes\%20documents/NC\%20EB/entso-e_balancing_in\%20_europe_report_Nov2018_web.pdf}.
%Type = Misc
\bibitem[{{ENTSO-E}(2018{\natexlab{b}})}]{ENTSOE2}
\bibinfo{author}{{ENTSO-E}} (\bibinfo{year}{2018}{\natexlab{b}}).
\newblock \bibinfo{title}{Commission regulation ({EU}) 2017/2195 - {EB}
  {D}eliverables}.
\newblock \URLprefix \url{https://www.entsoe.eu/network_codes/eb/}.
%Type = Article
\bibitem[{Fromen(1997)}]{fromen1997reducing}
\bibinfo{author}{Fromen, B.} (\bibinfo{year}{1997}).
\newblock \bibinfo{title}{Reducing the number of linear programs needed for
  solving the nucleolus problem of n-person game theory}.
\newblock {\it \bibinfo{journal}{Eur. J. of Oper. Res.}\/},  {\it
  \bibinfo{volume}{98}\/}, \bibinfo{pages}{626--636}.
%Type = Phdthesis
\bibitem[{Gebrekiros(2015)}]{gebrekiros2015analysis}
\bibinfo{author}{Gebrekiros, Y.~T.} (\bibinfo{year}{2015}).
\newblock {\it \bibinfo{title}{Analysis of integrated balancing markets in
  northern Europe under different market design options}\/}.
\newblock Ph.D. thesis NTNU.
%Type = Misc
\bibitem[{{Gurobi Optimization Inc.}(2016)}]{gurobi}
\bibinfo{author}{{Gurobi Optimization Inc.}} (\bibinfo{year}{2016}).
\newblock \bibinfo{title}{Gurobi optimizer reference manual}.
%Type = Article
\bibitem[{Hallefjord et~al.(1995)Hallefjord, Helming \&
  J{\o}rnsten}]{hallefjord1995computing}
\bibinfo{author}{Hallefjord, {\AA}.}, \bibinfo{author}{Helming, R.}, \&
  \bibinfo{author}{J{\o}rnsten, K.} (\bibinfo{year}{1995}).
\newblock \bibinfo{title}{Computing the nucleolus when the characteristic
  function is given implicitly: A constraint generation approach}.
\newblock {\it \bibinfo{journal}{Int. J. of Game Th.}\/},  {\it
  \bibinfo{volume}{24}\/}, \bibinfo{pages}{357--372}.
%Type = Article
\bibitem[{Hobbs et~al.(2005)Hobbs, Rijkers \& Boots}]{hobbs2005more}
\bibinfo{author}{Hobbs, B.~F.}, \bibinfo{author}{Rijkers, F.~A.}, \&
  \bibinfo{author}{Boots, M.~G.} (\bibinfo{year}{2005}).
\newblock \bibinfo{title}{The more cooperation, the more competition? a
  {C}ournot analysis of the benefits of electric market coupling}.
\newblock {\it \bibinfo{journal}{The Energy J.}\/},  (pp.
  \bibinfo{pages}{69--97}).
%Type = Article
\bibitem[{Jensen et~al.(2017)Jensen, Kazempour \& Pinson}]{jensen2017cost}
\bibinfo{author}{Jensen, T.~V.}, \bibinfo{author}{Kazempour, J.}, \&
  \bibinfo{author}{Pinson, P.} (\bibinfo{year}{2017}).
\newblock \bibinfo{title}{Cost-optimal {ATC}s in zonal electricity markets}.
\newblock {\it \bibinfo{journal}{IEEE Trans. on Pow. Syst.}\/},  {\it
  \bibinfo{volume}{33}\/}, \bibinfo{pages}{3624--3633}.
%Type = Misc
\bibitem[{Karaca et~al.(2019{\natexlab{a}})Karaca, Delikaraoglou, Hug \&
  Kamgarpour}]{karacadata}
\bibinfo{author}{Karaca, O.}, \bibinfo{author}{Delikaraoglou, S.},
  \bibinfo{author}{Hug, G.}, \& \bibinfo{author}{Kamgarpour, M.}
  (\bibinfo{year}{2019}{\natexlab{a}}).
\newblock \bibinfo{title}{{System data---Enabling inter-area reserves exchanges
  through stable benefit allocation mechanisms}}.
\newblock \bibinfo{howpublished}{\url{https://doi.org/10.5281/zenodo.3579082}}.
%Type = Article
\bibitem[{Karaca \& Kamgarpour(2018)}]{karaca2018core}
\bibinfo{author}{Karaca, O.}, \& \bibinfo{author}{Kamgarpour, M.}
  (\bibinfo{year}{2018}).
\newblock \bibinfo{title}{Core-selecting mechanisms in electricity markets}.
\newblock {\it \bibinfo{journal}{IEEE Trans. on Smart Grid,
  arXiv:1811.09646}\/},  {\it \bibinfo{volume}{to appear}\/}.
%Type = Article
\bibitem[{Karaca et~al.(2019{\natexlab{b}})Karaca, Sessa, Walton \&
  Kamgarpour}]{karaca2019designing}
\bibinfo{author}{Karaca, O.}, \bibinfo{author}{Sessa, P.~G.},
  \bibinfo{author}{Walton, N.}, \& \bibinfo{author}{Kamgarpour, M.}
  (\bibinfo{year}{2019}{\natexlab{b}}).
\newblock \bibinfo{title}{Designing coalition-proof reverse auctions over
  continuous goods}.
\newblock {\it \bibinfo{journal}{IEEE Trans. on Automatic Control}\/},  {\it
  \bibinfo{volume}{64}\/}, \bibinfo{pages}{4803--4810}.
%Type = Article
\bibitem[{Kimms \& {\c{C}}etiner(2012)}]{kimms2012approximate}
\bibinfo{author}{Kimms, A.}, \& \bibinfo{author}{{\c{C}}etiner, D.}
  (\bibinfo{year}{2012}).
\newblock \bibinfo{title}{Approximate nucleolus-based revenue sharing in
  airline alliances}.
\newblock {\it \bibinfo{journal}{Eur. J. of Oper. Res.}\/},  {\it
  \bibinfo{volume}{220}\/}, \bibinfo{pages}{510--521}.
%Type = Article
\bibitem[{Kohlberg(1972)}]{kohlberg1972nucleolus}
\bibinfo{author}{Kohlberg, E.} (\bibinfo{year}{1972}).
\newblock \bibinfo{title}{The nucleolus as a solution of a minimization
  problem}.
\newblock {\it \bibinfo{journal}{SIAM J. on App. Math.}\/},  {\it
  \bibinfo{volume}{23}\/}, \bibinfo{pages}{34--39}.
%Type = Techreport
\bibitem[{Kopelowitz(1967)}]{kopelowitz1967computation}
\bibinfo{author}{Kopelowitz, A.} (\bibinfo{year}{1967}).
\newblock {\it \bibinfo{title}{Computation of the kernels of simple games and
  the nucleolus of n-person games}\/}.
\newblock \bibinfo{type}{Technical Report} The Hebrew University of Jerusalem,
  Department of Math.
%Type = Book
\bibitem[{Krishna(2009)}]{krishna2009auction}
\bibinfo{author}{Krishna, V.} (\bibinfo{year}{2009}).
\newblock {\it \bibinfo{title}{Auction theory}\/}.
\newblock \bibinfo{publisher}{Academic press}.
%Type = Article
\bibitem[{Kristiansen et~al.(2018)Kristiansen, Mu{\~n}oz, Oren \&
  Korp{\aa}s}]{kristiansen2018mechanism}
\bibinfo{author}{Kristiansen, M.}, \bibinfo{author}{Mu{\~n}oz, F.~D.},
  \bibinfo{author}{Oren, S.}, \& \bibinfo{author}{Korp{\aa}s, M.}
  (\bibinfo{year}{2018}).
\newblock \bibinfo{title}{A mechanism for allocating benefits and costs from
  transmission interconnections under cooperation: A case study of the north
  sea offshore grid}.
\newblock {\it \bibinfo{journal}{The Energy J.}\/},  {\it
  \bibinfo{volume}{39}\/}.
%Type = Article
\bibitem[{Maschler et~al.(1979)Maschler, Peleg \&
  Shapley}]{maschler1979geometric}
\bibinfo{author}{Maschler, M.}, \bibinfo{author}{Peleg, B.}, \&
  \bibinfo{author}{Shapley, L.~S.} (\bibinfo{year}{1979}).
\newblock \bibinfo{title}{Geometric properties of the kernel, nucleolus, and
  related solution concepts}.
\newblock {\it \bibinfo{journal}{Math. of Oper. Res.}\/},  {\it
  \bibinfo{volume}{4}\/}, \bibinfo{pages}{303--338}.
%Type = Article
\bibitem[{Maschler et~al.(1992)Maschler, Potters \& Tijs}]{maschler1992general}
\bibinfo{author}{Maschler, M.}, \bibinfo{author}{Potters, J.~A.}, \&
  \bibinfo{author}{Tijs, S.~H.} (\bibinfo{year}{1992}).
\newblock \bibinfo{title}{The general nucleolus and the reduced game property}.
\newblock {\it \bibinfo{journal}{Int. J. of Game Th.}\/},  {\it
  \bibinfo{volume}{21}\/}, \bibinfo{pages}{85--106}.
%Type = Misc
\bibitem[{Mitridati et~al.(2019)Mitridati, Kazempour \&
  Pinson}]{mitridatidesign}
\bibinfo{author}{Mitridati, L.}, \bibinfo{author}{Kazempour, J.}, \&
  \bibinfo{author}{Pinson, P.} (\bibinfo{year}{2019}).
\newblock \bibinfo{title}{Design and game-theoretical analysis of
  community-based market mechanisms in heat and electricity systems}.
\newblock
  \bibinfo{howpublished}{\url{http://pierrepinson.com/docs/Mitridatietal2019.pdf}}.
\newblock \bibinfo{note}{WP.}
%Type = Article
\bibitem[{Morales et~al.(2012)Morales, Conejo, Liu \&
  Zhong}]{morales2012pricing}
\bibinfo{author}{Morales, J.~M.}, \bibinfo{author}{Conejo, A.~J.},
  \bibinfo{author}{Liu, K.}, \& \bibinfo{author}{Zhong, J.}
  (\bibinfo{year}{2012}).
\newblock \bibinfo{title}{Pricing electricity in pools with wind producers}.
\newblock {\it \bibinfo{journal}{IEEE Trans. on Pow. Syst.}\/},  {\it
  \bibinfo{volume}{27}\/}, \bibinfo{pages}{1366--1376}.
%Type = Article
\bibitem[{Morales et~al.(2014)Morales, Zugno, Pineda \&
  Pinson}]{morales2014electricity}
\bibinfo{author}{Morales, J.~M.}, \bibinfo{author}{Zugno, M.},
  \bibinfo{author}{Pineda, S.}, \& \bibinfo{author}{Pinson, P.}
  (\bibinfo{year}{2014}).
\newblock \bibinfo{title}{Electricity market clearing with improved scheduling
  of stochastic production}.
\newblock {\it \bibinfo{journal}{Eur. J. of Oper. Res.}\/},  {\it
  \bibinfo{volume}{235}\/}, \bibinfo{pages}{765--774}.
%Type = Article
\bibitem[{Neuhoff(2005)}]{neuhoff2005large}
\bibinfo{author}{Neuhoff, K.} (\bibinfo{year}{2005}).
\newblock \bibinfo{title}{Large-scale deployment of renewables for electricity
  generation}.
\newblock {\it \bibinfo{journal}{Ox. Rev. of Econ. P.}\/},  {\it
  \bibinfo{volume}{21}\/}, \bibinfo{pages}{88--110}.
%Type = Book
\bibitem[{Osborne \& Rubinstein(1994)}]{osborne1994course}
\bibinfo{author}{Osborne, M.~J.}, \& \bibinfo{author}{Rubinstein, A.}
  (\bibinfo{year}{1994}).
\newblock {\it \bibinfo{title}{A course in game theory}\/}.
\newblock \bibinfo{publisher}{MIT press}.
%Type = Article
\bibitem[{Owen(1974)}]{owen1974note}
\bibinfo{author}{Owen, G.} (\bibinfo{year}{1974}).
\newblock \bibinfo{title}{A note on the nucleolus}.
\newblock {\it \bibinfo{journal}{Int. J. of Game Th.}\/},  {\it
  \bibinfo{volume}{3}\/}, \bibinfo{pages}{101--103}.
%Type = Article
\bibitem[{Owen(1975)}]{owen1975core}
\bibinfo{author}{Owen, G.} (\bibinfo{year}{1975}).
\newblock \bibinfo{title}{On the core of linear production games}.
\newblock {\it \bibinfo{journal}{Math. Programming}\/},  {\it
  \bibinfo{volume}{9}\/}, \bibinfo{pages}{358--370}.
%Type = Misc
\bibitem[{Pandzic et~al.(2014)Pandzic, Dvorkin, Qiu, Wang \&
  Kirschen}]{pandzicunit}
\bibinfo{author}{Pandzic, H.}, \bibinfo{author}{Dvorkin, Y.},
  \bibinfo{author}{Qiu, T.}, \bibinfo{author}{Wang, Y.}, \&
  \bibinfo{author}{Kirschen, D.} (\bibinfo{year}{2014}).
\newblock \bibinfo{title}{{Library of the Renewable Energy Analysis Lab (REAL),
  University of Washington, Seattle, USA}}.
\newblock
  \bibinfo{howpublished}{\url{https://labs.ece.uw.edu/real/library.html}}.
%Type = Article
\bibitem[{Pineda \& Morales(2016)}]{pineda2016capacity}
\bibinfo{author}{Pineda, S.}, \& \bibinfo{author}{Morales, J.~M.}
  (\bibinfo{year}{2016}).
\newblock \bibinfo{title}{Capacity expansion of stochastic power generation
  under two-stage electricity markets}.
\newblock {\it \bibinfo{journal}{Computers \& Oper. Res.}\/},  {\it
  \bibinfo{volume}{70}\/}, \bibinfo{pages}{101--114}.
%Type = Article
\bibitem[{Pritchard et~al.(2010)Pritchard, Zakeri \&
  Philpott}]{pritchard2010single}
\bibinfo{author}{Pritchard, G.}, \bibinfo{author}{Zakeri, G.}, \&
  \bibinfo{author}{Philpott, A.} (\bibinfo{year}{2010}).
\newblock \bibinfo{title}{A single-settlement, energy-only electric power
  market for unpredictable and intermittent participants}.
\newblock {\it \bibinfo{journal}{Oper. Res.}\/},  {\it \bibinfo{volume}{58}\/},
  \bibinfo{pages}{1210--1219}.
%Type = Article
\bibitem[{Ruiz \& Contreras(2007)}]{ruiz2007effective}
\bibinfo{author}{Ruiz, P.~A.}, \& \bibinfo{author}{Contreras, J.}
  (\bibinfo{year}{2007}).
\newblock \bibinfo{title}{An effective transmission network expansion cost
  allocation based on game th.}
\newblock {\it \bibinfo{journal}{IEEE Trans. on Pow. Syst.}\/},  {\it
  \bibinfo{volume}{22}\/}, \bibinfo{pages}{136--144}.
%Type = Article
\bibitem[{Schmeidler(1969)}]{schmeidler1969nucleolus}
\bibinfo{author}{Schmeidler, D.} (\bibinfo{year}{1969}).
\newblock \bibinfo{title}{The nucleolus of a characteristic function game}.
\newblock {\it \bibinfo{journal}{SIAM J. on {A}pp. {M}ath.}\/},  {\it
  \bibinfo{volume}{17}\/}, \bibinfo{pages}{1163--1170}.
%Type = Article
\bibitem[{Shapley(1967)}]{shapley1967balanced}
\bibinfo{author}{Shapley, L.~S.} (\bibinfo{year}{1967}).
\newblock \bibinfo{title}{On balanced sets and cores}.
\newblock {\it \bibinfo{journal}{Naval Res. Logistics Quarterly}\/},  {\it
  \bibinfo{volume}{14}\/}, \bibinfo{pages}{453--460}.
%Type = Article
\bibitem[{Shapley(1971)}]{shapley1971cores}
\bibinfo{author}{Shapley, L.~S.} (\bibinfo{year}{1971}).
\newblock \bibinfo{title}{Cores of convex games}.
\newblock {\it \bibinfo{journal}{Int. J. of Game Th.}\/},  {\it
  \bibinfo{volume}{1}\/}, \bibinfo{pages}{11--26}.
%Type = Article
\bibitem[{Shapley \& Shubik(1966)}]{shapley1966quasi}
\bibinfo{author}{Shapley, L.~S.}, \& \bibinfo{author}{Shubik, M.}
  (\bibinfo{year}{1966}).
\newblock \bibinfo{title}{Quasi-cores in a monetary economy with nonconvex
  preferences}.
\newblock {\it \bibinfo{journal}{Econ.}\/},  (pp. \bibinfo{pages}{805--827}).
%Type = Article
\bibitem[{Shapley \& Shubik(1969)}]{shapley1969market}
\bibinfo{author}{Shapley, L.~S.}, \& \bibinfo{author}{Shubik, M.}
  (\bibinfo{year}{1969}).
\newblock \bibinfo{title}{On market games}.
\newblock {\it \bibinfo{journal}{J. of Econ. Theory}\/},  {\it
  \bibinfo{volume}{1}\/}, \bibinfo{pages}{9--25}.
%Type = Article
\bibitem[{Sobolev(1975)}]{sobolev1975characterization}
\bibinfo{author}{Sobolev, A.~I.} (\bibinfo{year}{1975}).
\newblock \bibinfo{title}{The characterization of optimality principles in
  cooperative games by functional equations}.
\newblock {\it \bibinfo{journal}{Math. {M}ethods in the {S}ocial {S}ci.}\/},
  {\it \bibinfo{volume}{6}\/}, \bibinfo{pages}{151}.
%Type = Article
\bibitem[{Young(1985)}]{young1985monotonic}
\bibinfo{author}{Young, H.~P.} (\bibinfo{year}{1985}).
\newblock \bibinfo{title}{Monotonic solutions of cooperative games}.
\newblock {\it \bibinfo{journal}{Int. J. of Game Th.}\/},  {\it
  \bibinfo{volume}{14}\/}, \bibinfo{pages}{65--72}.
%Type = Article
\bibitem[{Young et~al.(1982)Young, Okada \& Hashimoto}]{young1982cost}
\bibinfo{author}{Young, H.~P.}, \bibinfo{author}{Okada, N.}, \&
  \bibinfo{author}{Hashimoto, T.} (\bibinfo{year}{1982}).
\newblock \bibinfo{title}{Cost allocation in water resources development}.
\newblock {\it \bibinfo{journal}{Water {R}esources {R}es.}\/},  {\it
  \bibinfo{volume}{18}\/}, \bibinfo{pages}{463--475}.

\end{thebibliography}

\newpage
\appendix
\section{Notation}\label{app:not}
Main notation is stated below. An illustration is provided in Figure~\ref{fig:notation}. Additional symbols are defined throughout the paper when~needed.

\begin{itemize}
\setstretch{1}
\setlength{\itemsep}{0pt}
\item[] \textbf{Sets and indices}
	\item [$\AC$] Set of areas indexed by $a$
% 	\item [$a_r(e)$] Receiving-end area of link $e$.
% 	\item [$a_s(e)$] Sending-end area of link $e$.
	\item [$\mathcal E$] Set of inter-area links indexed by $e$
	\item [$\mathcal L^{\text{AC}}$] Set of AC transmission lines indexed by~$\ell$
	\item [$a_r(e/\ell)$] Receiving-end area of link $e$/line $\ell$ 
	\item [$a_s(e/\ell)$] Sending-end area of link $e$/line $\ell$ 
	\item [$ \Lambda_e$] Set of tie-lines across link~$e$
	%\item [$ \Lambda$] Set of lines that are across some link~$e\in \mathcal E$
	\item [$\mathcal I$] Set of dispatchable power plants indexed by $i$
	\item [$\mathcal J$] Set of stochastic power plants indexed by $j$
% 	\item [$L^{\text{DC}}$] Set of HVDC transmission lines.
	\item [$\mathcal N$] Set of network nodes (buses) indexed by $n$.
	%\item [$ \mathcal{M}_{a}^{\AC}$] Set of areas that belong in $\AC$ and are interconnected with area $a$
	\item [$ \mathcal{M}_{n}^{\mathcal I}$] Set of dispatchable power plants  $i$ located at node $n$
	\item [$ \mathcal{M}_{n}^{\mathcal J}$] Set of stochastic power plants $j$ located at node $n$
	\item [$ \mathcal{M}_{a}^{\mathcal I}$] Set of dispatchable power plants  $i$ located in area $a$
	\item [$ \mathcal{S}$] Set of stochastic power production scenarios indexed by~$s$
% \item[$N^{\text{G}}$] Set of natural gas network nodes $n_{g}$

% PARAMETERS
\item[] \textbf{Parameters}
\item [$\overline{{W}}_{j}$] Expected power production of stochastic power plant~$j$ [MW]
\item [$W_{js}$] Power production by stochastic power plant $j$ in scenario~$s$ [MW]
\item [$\pi_{s}$] Probability of occurrence of scenario $s$
\item [$A_{\ell n}$] Line-to-bus incidence matrix
\item [$B_{\ell}$] Absolute value of the susceptance of AC line $\ell$
\item [$D_{n}$] Demand at node $n$ [MW]
\item [$C_{i}$] Energy offer price of power plant $i$ [\euro/MWh]
\item [$C^{+/-}_{i}$] Up/down reserve capacity offer price of power plant $i$ [\euro/MW]
\item [$C^{\text{sh}}$] Value of involuntarily shed load [\euro/MWh]
\item [$P_{i}$] Capacity of dispatchable power plant $i$ [MW]
\item [$R^{+/-}_{i}$] Up/down reserve capacity offer quantity of power plant $i$ [\euro/MW]
\item [$RR_a^{+/-}$] Up/down reserve capacity requirements of area $a$ [MW]
\item [$T_{e/\ell}$] Transmission capacity of link $e$/line $\ell$ [MW]

% VARIABLES
\item[] \textbf{Variables}

\item [$\delta_n$] Voltage angle at node $n$ at day-ahead stage [rad]
\item [${\delta}_{ns}$] Voltage angle at node $n$ in scenario $s$ [rad]
\item [$\chi_{e/\ell}$] Transmission allocation of link $e$/line $\ell$, i.e., percentage of inter-area interconnection capacity of link $e$/line $\ell$ allocated to reserves exchange
\item [$f_{\ell}$] Power flow in AC line $\ell$ at day-ahead stage [MW]
\item [${f}_{\ell s}$] Power flow in AC line $\ell$ in scenario $s$ [MW]
\item [$l^{\text{sh}}_{ns}$] Load shedding at node $n$ in scenario $s$ [MW]
\item [${p}_{i}$] Day-ahead schedule of dispatchable power plant $i$ [MW]
\item [${p}_{is}^{+/-}$]Up/down regulation provided by dispatchable power plant $i$ in scenario $s$ [MW]
\item [$r^{+/-}_{e}$] Up/down reserve capacity (`exported') from area $a_s(e)$ to area $a_r(e)$ [MW] (equivalently imported to area $a_r(e)$ from area $a_s(e)$)
%\item [$r^{+/-}_{a' a}$] Up/down reserve capacity (`exported') from area $a'$ to area $a$ [MW] (equivalently imported to area $a$ from area $a'$)
\item [${w}_{j}$] Day-ahead schedule of stochastic power plant $j$ [MW]
\item [$w^{\text{spill}}_{js}$] Power spilled by stochastic power plant $j$ in scenario $s$ [MW]
% \item [$z_{\ell}$] Power flow in HVDC line $\ell$ at day-ahead stage [MW].
% \item [$\tilde{z}_{\ell s}$] Power flow in HVDC line $\ell$ in scenario $s$ [MW].
\setstretch{1.2}
\end{itemize}
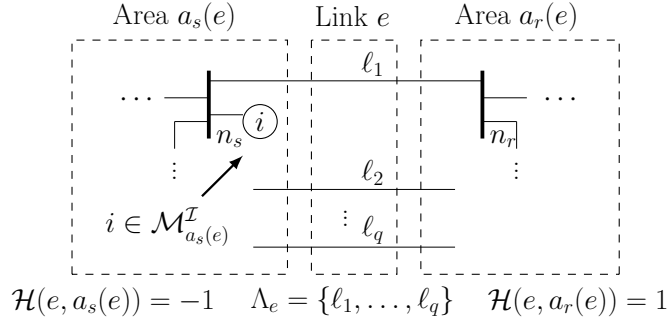
\begin{figure}[h]
	\centering
	\begin{tikzpicture}[scale=0.9, every node/.style={scale=0.63}]
	\draw[-,line width=.5mm] (-2,.5) -- (-2,-.5) node[anchor=west]  {\LARGE $n_s$};
	\draw[-,line width=.5mm] (2,.5) -- (2,-.5) node[anchor=west]  {\LARGE $n_r$};	
	\draw[-,line width=.1mm] (-2,0.35) -- (2,0.35) node at (.4, .6) {\LARGE $\ell_1$};
	\draw[-,line width=.1mm] (-1.35,-1.25) -- (1.6,-1.25) node at (.4, -1) {\LARGE $\ell_2$}; %connect
	\draw[-,line width=.1mm] (-1.35,-2.1) -- (1.6,-2.1) node at (0, -1.6) {\LARGE $\vdots$} node at (.4, -1.85) {\LARGE $\ell_q$}; %connect
	\draw[-,line width=.1mm] (-2,-0.15) -- (-1.5,-0.15); %go to i
	\draw[-,line width=.1mm] (-2,0.1) -- (-2.65,0.1) node[anchor=west] at (2.8, 0.075)  {\LARGE $\cdots$}; % expand to the right
	\draw[-,line width=.1mm] (2.65,0.1) -- (2,0.1) node[anchor=west] at (-3.4, 0.075) {\LARGE $\cdots$} ; % expand to the left
	\draw[-,line width=.1mm] (-2,-0.25) -- (-2.5,-0.25); % expand down left
	\draw[-,line width=.1mm] (2.5,-0.25) -- (2,-0.25); % expand down right
	\draw[-,line width=.1mm] (-2.5,-0.65) -- (-2.5,-0.25) node at (-2.5, -.9) {\LARGE $\vdots$} ;
	\draw[-,line width=.1mm] (2.5,-0.25) -- (2.5,-0.65) node at (2.5, -.9) {\LARGE $\vdots$} ;
	\draw (-1.25,-0.25) circle (.25cm) node {\LARGE $i$} node at (-2.6, -1.8) {\LARGE $i\in\mathcal M^{\mathcal I}_{a_s(e)}$} ;
	\draw[->,line width=.3mm] (-2.1,-1.35) -- (-1.5,-.75);
	\draw[dashed, draw=black] (-4,-2.5) rectangle (-0.85,0.95) node at (-2.5, 1.3) {\LARGE $\text{Area}\ a_s(e)$} node at (-3.4, -2.9) {\LARGE $\mathcal{H}(e,a_s(e))=-1$};
	\draw[dashed,draw=black] (4,-2.5) rectangle (1.1,0.95) node at (2.5, 1.3) {\LARGE $\text{Area}\ a_r(e)$} node at (3.4, -2.9) {\LARGE $\mathcal{H}(e,a_r(e))=1$};
	\draw[dashed,draw=black] (-0.5,0.95) rectangle (.75,-2.5) node at (.1, 1.3) {\LARGE $\text{Link}\ e$} node at (.1, -2.9) {\LARGE $\Lambda_e=\{\ell_1,\ldots,\ell_q\}$};
	\end{tikzpicture}\vspace{-.2cm}
	\caption{An illustration of the main notation used for graph $(\mathcal A,\mathcal E)$ }\label{fig:notation}
\end{figure}

\section{Data for two-area example}\label{app:two_area_ex}
For this example, we removed area 3 from the example in Section~\ref{sec:illu_ex} while all the other parameters are kept unchanged.
The different levels of wind power penetration are modeled by changing the installed capacities of the wind power plants $j_3$ and $j_6$. $312$ MW corresponds to $120$ and $192$ MW, respectively. $273$ MW corresponds to $105$ and $168$ MW, whereas $234$ MW corresponds to $90$ and $144$ MW. Normalization is done by dividing the expected operation cost by the cost under the same wind power penetration but with $\chi=\epsilon$ where $\epsilon>0$ is a small positive number.
\section{KKT conditions}\label{app:a}
We provide the complete set of KKT conditions for the reserve and day-ahead markets in~\eqref{eq:LLR} and~\eqref{eq:LLD}, that appear in the lower level of the stochastic bilevel optimization in \eqref{mod:B-Pre}. The dual multipliers of inequality constraints are listed to the right of the complementarity relationships denoted by $\perp$. For the equality constraints, the dual multipliers are listed after a colon.

The KKT conditions for the reserve market in~\eqref{eq:LLR} are:
\begin{align*}
		& 0 \le  {R}^{+}_{i}- r_{i}^{+}\perp \mu_i^{\text{R}^{+}}\geq 0, \quad \forall i,\\
		& 0 \le  {R}^{-}_{i}-r_{i}^{-}\perp \mu_i^{\text{R}^{-}}\geq 0, \quad \forall i, \\
		 &0 \le \sum_{i \in \mathcal{M}_{a}^{\mathcal I}} r_{i}^{+} + \sum_{e\in\mathcal E} \mathcal{H}(e,a) r_{e}^{+} - RR_a^{+}\perp \mu_a^{\text{RR}^{+}}\geq0,\quad \forall a, \\
	    &0 \le \sum_{i \in \mathcal{M}_{a}^{\mathcal I}} r_{i}^{-} + \sum_{e\in\mathcal E} \mathcal{H}(e,a) r_{e}^{-} - RR_a^{-}\perp \mu_a^{\text{RR}^{-}}\geq0,\quad \forall a, \\
		& 0 \le r_{e}^{+}+\chi_{e}' T_{e}\perp \zeta_e^{\text{L}^{+}}\geq0,\quad \forall e, \\
		& 0 \le  \chi_{e}' T_{e} - r_{e}^{+}\perp \zeta_e^{\text{U}^{+}}\geq0,\quad \forall e, \\
		& 0 \le r_{e}^{-}+\chi_{e}' T_{e}\perp \zeta_e^{\text{L}^{-}}\geq0,\quad \forall e, \\
		& 0 \le \chi_{e}' T_{e}-r_{e}^{-}\perp \zeta_e^{\text{U}^{-}}\geq0,\  \forall e,\\
		&0\le C_i^++\mu_i^{\text{R}^{+}}-\sum_{a: i \in \mathcal{M}_{a}^{\mathcal I}}\mu_a^{\text{RR}^{+}}\perp r_{i}^{+}\geq 0,\quad \forall i, \\
		&0\le C_i^-+\mu_i^{\text{R}^{-}}-\sum_{a: i \in \mathcal{M}_{a}^{\mathcal I}}\mu_a^{\text{RR}^{-}}\perp r_{i}^{-}\geq 0,\quad \forall i, \\
		&\sum_{a} \mu_a^{\text{RR}^{+}}\mathcal{H}(e,a) + \zeta_e^{\text{L}^{+}} - \zeta_e^{\text{U}^{+}} = 0,\quad \forall e,\\
		&\sum_{a} \mu_a^{\text{RR}^{-}}\mathcal{H}(e,a) + \zeta_e^{\text{L}^{-}} - \zeta_e^{\text{U}^{-}} = 0,\quad \forall e.
\end{align*}

The KKT conditions for the day-ahead market in~\eqref{eq:LLD} are:
\begin{align*}
    		& \sum_{j \in \mathcal{M}_n^{\mathcal J}} w_j + \sum_{i\in \mathcal{M}_n^{\mathcal I}}{p}_{i} - \sum_{\ell \in \mathcal L^{\text{AC}}} A_{\ell n} f_{\ell}   = D_n: \lambda_n\ \text{free},\quad \forall n,\\
		&  0 \le p_i-{r}_{i}^{-}\perp\mu_i^{\text{PL}}\geq 0,   \quad \forall i, \\
		&  0 \leq  {P}_{i}- {r}_{i}^{+}-p_i\perp\mu_i^{\text{PU}}\geq 0,   \  \forall i, \\
		& 0 \leq \overline{{W}}_{j}-w_j\perp\mu_j^{\text{WU}}\geq 0,  \quad \forall j,\\
		& f_{\ell} = B_{\ell} \sum_{n\in\mathcal N} A_{\ell n} \delta_n: \lambda_\ell^\text{F}\ \text{free},\quad \forall \ell, \\
		& 0 \le f_{\ell}+(1- \chi_{\ell}') \ {T}_{\ell}\perp\zeta_\ell^{\text{L}}\geq0,\quad \forall \ell ,\\
		&0 \le (1- \chi_{\ell}') \ {T}_{\ell}-f_{\ell}\perp\zeta_\ell^{\text{U}}\geq0,\quad \forall \ell,\\	
		& \delta_1 = 0: \lambda^\text{REF}\ \text{free},\\
		& C_i + \sum_{n:i\in \mathcal{M}_n^{\mathcal I}}\lambda_n-\mu_i^{\text{PL}}+\mu_i^{\text{PU}}=0,\quad \forall i,\\
		&0\leq  \sum_{n:j \in \mathcal{M}_n^{\mathcal J}}\lambda_n +\mu_j^{\text{WU}}\perp w_j\geq 0, \quad \forall j,\\
		&\sum_{n}A_{\ell n}\lambda_n-\lambda_\ell^\text{F}-\zeta_\ell^{\text{L}}+\zeta_\ell^{\text{U}}=0, \quad \forall \ell,\\
		&\sum_{\ell}\lambda_\ell^\text{F}B_{\ell}  A_{\ell n}=0, \quad \forall n\neq 1,\\
		&\sum_{\ell}\lambda_\ell^\text{F}B_{\ell} A_{\ell n}-\lambda^\text{REF}=0, \ n=1.
\end{align*}
Note that the conditions above involve the solutions of the reserve market ${r}_{i}^{+}$ and ${r}_{i}^{-}$, and the transmission allocations $\chi'$ from the optimization variables of preemptive model~\eqref{mod:B-Pre}.
\section{Unique properties of the Shapley value}\label{app:a2}
In this section, we analyze each of the unique properties of the Shapley value and how they relate to the fundamental properties associated with the core and the least-core.
    \subsection{Dummy player} Dummy player property requires $\beta_a=0$ for all $a$ such that $v(\CC)-v(\CC\setminus a)=0$ for all~$\CC\subseteq\AC$. In other words, an area incapable of contributing to any coalition~$\CC$ ends up with zero benefits. We now reiterate Proposition~\ref{prop:dummy} and then prove the two claims.
%The core provides us with a more restrictive version of the dummy player property, i.e., $\beta_a=0$ for all $a$ such that $v(\AC)-v(\AC\setminus\{a\})=0$. For the least-core, we have the following result.
%It is straightforward to check that for such areas 
\begin{proposition}\label{prop:dummy2} For the core and the least-core, we have,
	\begin{itemize}
		\item[(i)] 	If $a'$ satisfies $v(\AC)-v(\AC\setminus a')=0$, then $K_\text{Core}(v)\subset\{\beta\,|\,\beta_{a'}=0\}$,
		\item[(ii)]	If $a'$ satisfies $v(\CC)-v(\CC\setminus a')=0$ for all~$\CC\subseteq\AC$, then $K_\text{Core}(v,\epsilon^*(v))\subset\{\beta\,|\,\beta_{a'}=0\}$.
	\end{itemize}
\end{proposition}
\begin{pf}
\setstretch{1.2}
	(i) Assume core is nonempty, since otherwise the proof is trivial. Combining the equality constraint with the inequality constraint corresponding to $\AC\setminus a'$, we obtain $\beta_{a'}\leq v(\AC)-v(\AC\setminus a')=0$ for any $\beta\in K_\text{Core}(v)$. Combining this with $K_\text{Core}(v)\subset\mathbb{R}^\AC_+$ gives us $\beta_{a'}=0$ for any $\beta\in K_\text{Core}(v)$.
	
	(ii) Assume core is empty, $\epsilon^*(v)>0$, since otherwise part~(i) concludes that $K_\text{Core}(v,\epsilon^*(v))\subseteq K_\text{Core}(v)\subset\{\beta\,|\,\beta_{a'}=0\}.$ Next, we prove by contradiction that first $\beta_{a'}>0$ is not possible and then $\beta_{a'}<0$ is not possible.
	
	Let $\hat \beta \in K_\text{Core}(v,\epsilon^*(v))$ be a benefit allocation with $\hat \beta_{a'}>0$. We now show that there exists $\epsilon <\epsilon^*(v)$ such that $K_\text{Core}(v,\epsilon)\neq \emptyset$. This would contradict the definition of the least-core. 
	
	For any $\CC\ni a'$, we have $\sum_{a\in\CC}\hat\beta_a>\sum_{a\in\CC\setminus a'}\hat\beta_a\geq v(\CC\setminus a')-\epsilon^*(v)=v(\CC)-\epsilon^*(v)$. Notice that we can always find a small positive number $\delta$ such that  $\sum_{a\in\CC}\hat\beta_a - (|\AC|-|\CC|)\delta> v(\CC)-\epsilon^*(v)+\delta$ holds for any $\CC\ni a'$. Next, we show that $K_\text{Core}(v,\epsilon^*(v)-\delta)$ is nonempty for this particular choice. 
	
	Define $\bar{\beta}$ such that $\bar\beta_a=\hat \beta_a + \delta$ for all $a\neq a'$and $\bar{\beta}_{a'}=\hat{\beta}_{a'}-(|\AC|-1)\delta$. This new allocation $\bar{\beta}$ clearly satisfies the equality constraint in $K_\text{Core}(v,\epsilon^*(v)-\delta)$. For inequality constraints $\CC\ni a'$, we have $\sum_{a\in\CC}\bar\beta_a=\sum_{a\in\CC}\hat\beta_a - (|\AC|-|\CC|)\delta> v(\CC)-\epsilon^*(v)+\delta$, where the strict inequality follows from the definition of $\delta$. For inequality constraints $\CC\not\ni a'$, we have $\sum_{a\in\CC}\bar\beta_a\geq\sum_{a\in\CC}\hat\beta_a+\delta\geq v(\CC)-\epsilon^*(v)+\delta$. Hence, $\bar \beta \in K_\text{Core}(v,\epsilon^*(v)-\delta)$, in other words, $K_\text{Core}(v,\epsilon^*(v)-\delta)\neq\emptyset.$ This contradicts $K_\text{Core}(v,\epsilon^*(v))$ being the least-core. Hence, $\hat \beta_{a'}\not>0$. 
	
	Next, let $\hat \beta \in K_\text{Core}(v,\epsilon^*(v))$ be a benefit allocation with $\hat \beta_{a'}<0$. We again show that there exists $\epsilon <\epsilon^*(v)$ such that $K_\text{Core}(v,\epsilon)\neq \emptyset$. 
	
	Since $\hat \beta \in K_\text{Core}(v,\epsilon^*(v))$ and $v(a')=0$, we have $0>\beta_{a'}\geq v(a') -\epsilon^*(v)=-\epsilon^*(v)$. Notice that, for any $\CC\not\ni a'$, we have $\sum_{a\in\CC}\hat\beta_a\geq v(\CC\cup a')-\epsilon^*(v)-\beta_{a'}$ by adding and subtracting $\beta_{a'}$, and by using the fact that $\epsilon^*(v)>0$ for the special case corresponding to $\CC\cup a'=\AC$. Since we have $v(\CC\cup a')=v(\CC)$ and $\beta_{a'}<0$, we obtain $\sum_{a\in\CC}\hat\beta_a> v(\CC)-\epsilon^*(v)$. 
	Notice that we can always find a small positive number $\delta$ such that  $\sum_{a\in\CC}\hat\beta_a - |\CC|\delta> v(\CC)-\epsilon^*(v)+\delta$ holds for any $\CC\not\ni a'$. Next, we show that $K_\text{Core}(v,\epsilon^*(v)-\delta)$ is nonempty for this particular choice. 
	
	Define $\bar{\beta}$ such that $\bar\beta_a=\hat \beta_a - \delta$ for all $a\neq a'$and $\bar{\beta}_{a'}=\hat{\beta}_{a'}+(|\AC|-1)\delta$. This new allocation $\bar{\beta}$ clearly satisfies the equality constraint in $K_\text{Core}(v,\epsilon^*(v)-\delta)$. For inequality constraints $\CC\not\ni a'$, we have $\sum_{a\in\CC}\bar\beta_a=\sum_{a\in\CC}\hat\beta_a - |\CC|\delta> v(\CC)-\epsilon^*(v)+\delta$, where the strict inequality follows from the definition of $\delta$. For inequality constraints $\CC\ni a'$, we have $\sum_{a\in\CC}\bar\beta_a\geq\sum_{a\in\CC}\hat\beta_a+\delta\geq v(\CC)-\epsilon^*(v)+\delta$. Hence, $\bar \beta \in K_\text{Core}(v,\epsilon^*(v)-\delta)$, in other words, $K_\text{Core}(v,\epsilon^*(v)-\delta)\neq\emptyset.$ This contradicts $K_\text{Core}(v,\epsilon^*(v))$ being the least-core.  Hence, $\hat \beta_{a'}\not<0$. This concludes that $\hat \beta_{a'}=0$.\hfill \QEDA
\end{pf}

The proposition above provides a missing link in the comparisons of the Shapley value, the core, and the least-core in a generic coalitional game. It shows that the core attains a more restrictive version of the dummy player property, i.e., $\beta_a=0$ for all $a$ such that $v(\AC)-v(\AC\setminus a)=0$. In other words, an area incapable of contributing to the set of all areas $\AC$ ends up with zero benefits. Finally, the least-core attains the dummy player property in the same way that it is defined for the Shapley value.

\subsection{Symmetry} Symmetry property is achieved if the benefit allocations of two areas are the same whenever their marginal contributions to any coalition~$\CC$ are the same. It can be verified that this property does not hold for every benefit allocation from the core and the least-core. However, it is always possible to find a benefit allocation satisfying the symmetry property in any nonempty strong $\epsilon$-core (and hence both in the core and the least-core) since the linear inequality constraints imposed by the convex polytope $K_\text{Core}(v,\epsilon)$ on two such areas are identical. Notice that the symmetry property is computationally hard to check since it would require evaluating the function $v$ for all coalitions. We can instead aim for a more restrictive version of the symmetry property by considering only the marginal contributions to the set of all areas $\AC$. This stronger condition would be computationally tractable. %This idea will be used later in our proposals.

\subsection{Additivity} Additivity property is given by $\beta(\hat v + v)=\beta(\hat v)+\beta(v)$ for all $\hat v,v:2^\AC\rightarrow \mathbb{R}$. As it is discussed in~\citep{osborne1994course}, this property is mathematically convenient but hard to argue for since the sum of coalitional value functions is in general considered to induce an unrelated coalitional game. As a remark, additivity further implies that the benefit allocation of any player responds monotonically to changes in the coalition value $v(\CC)$ (a positive change if $\CC$ contains the area). In the general case, this property cannot be achieved by any benefit allocation chosen from both the core and the least-core, see the discussions and the counter examples in \citep{young1985monotonic}.
	
\section{Definitions and discussions for the nucleolus allocation}\label{app:a2-2}

In this section, we provide the mathematical definition for the nucleolus allocation, compare it with the Shapley value, and discuss additional aspects of its computation and lexicographic minimization property.

Denote the excesses as $\theta(v,\beta,\CC)=v(\CC)-\sum_{a\in\CC}\beta_a$ for any nonempty $\CC\subset\AC$. Let $\theta(v,\beta)\in\mathbb{R}^{2^{|\AC|}-2}$ be the vector whose entries are the excesses but arranged in a nonincreasing order. Given two such ordered vectors $x,y\in\mathbb R^{n_0}$, with $x<_{\text L}y$ we mean that $x$ is lexicographically smaller than $y$, i.e., there exists an index $\nu_0\leq n_0$ such that $x_\nu = y_\nu$ for all $\nu<\nu_0$, and $x_{\nu_0}<y_{\nu_0}$. Let $\X\subset\mathbb{R}^{\AC}$ denote a set of benefit allocations that we are interested in. Then, the nucleolus of the set $\X$ is the benefit allocation that minimizes the excess of all coalitions in a lexicographic manner among all benefit allocations from the set $\X$.
Specifically, the nucleolus of the set $\X$, $\beta^{\text{Nuc}}(v,\X)\in\X$,  is defined by
 \begin{equation*}
 \theta(v,\beta^{\text{Nuc}}(v,\X))<_{\text L}\theta(v,\hat\beta),\ \forall \hat \beta\in\X\ \text{s.t.}\ \hat \beta\neq\beta^{\text{Nuc}}(v,\X).
 \end{equation*}

 In the literature, the nucleolus is generally defined with respect to two sets. Let $\X_{\text{BB}}=\{\beta\in\R^{\AC}|\sum_{a\in\AC}\beta_a=v(\AC)\}$ be the set of all efficient benefit allocations.
 The nucleolus of the set $\X_{\text{BB}}$, $\beta^{\text{Nuc}}(v,\X_{\text{BB}})$, was introduced in~\citet{sobolev1975characterization}, and is also called the prenucleolus.
 This allocation always exists for any coalitional game, and it is unique. Moreover, it lies in the least-core since its definition can be regarded as a stronger version of minimizing the maximum excess among all efficient benefit allocations~\citep{maschler1979geometric}. Hence, it satisfies efficiency, approximate stability and approximate individual rationality.
 
 On the other hand, let $\X_{\text{BB,IR}}=\{\beta\in\R^{\AC}_+|\sum_{a\in\AC}\beta_a=v(\AC)\}$ be the set of all efficient and individually rational benefit allocations.
 The nucleolus of the set $\X_{\text{BB,IR}}$, $\beta^{\text{Nuc}}(v,\X_{\text{BB,IR}})$, was introduced in \citep{schmeidler1969nucleolus}. This allocation is again unique, but it exists if and only if $v(\AC)\geq 0$. It is guaranteed to lie in the least-core only when the coalitional value function~$v$ is nondecreasing, since then it was proven that this allocation coincides with $\beta^{\text{Nuc}}(v,\X_{\text{BB}})$~\citep{aumann1985game}.
 
 In comparison with the Shapley value, both nucleolus allocations above are consistent with the symmetry property but they do not satisfy the additivity property in general, see~\citep{maschler1979geometric}. For the special class of convex graph games, these two nucleolus allocations coincide with the Shapley value~\citep{deng1994complexity}. For our work, the nucleolus will refer to the nucleolus of the set $\X_{\text{BB}}$, since this allocation always exists for a general coalitional game, and it is always a unique allocation from the least-core.
 
 In terms of computational approaches, note that there are also methods to compute the nucleolus by solving a single linear program involving either $4^{|\AC|}$ constraints~\citep{owen1974note} or $2^{|\AC|}!$ constraints~\citep{kohlberg1972nucleolus}. To the best of our knowledge, there are no iterative approaches applicable, and these two methods necessitate the complete evaluation of the coalitional value function. Finally, as discussed in~\citep{maschler1992general}, the lexicographic minimization property attained by the nucleolus allocation may not be relevant to the needs of every application, and it may even be considered hard to grasp in some cases. There could be more intuitively acceptable properties that yield a unique point. For instance, a prominent example is the work by \citet{young1982cost} suggesting to allocate benefits in a water supply project in proportion to the population and the total demand of an area.

\section{Proof of Proposition~\ref{prop:corempt}}\label{app:a3}

	Define $\hat\beta\in\mathbb{R}^{\AC}$ such that $\hat\beta_{a'}=\bar v(\AC)$, and $\hat\beta_{a}=0$ for all $a\neq a'.$ We prove by showing that this allocation lies in the core $K_\text{Core}(\bar v)$. The equality constraint in $K_\text{Core}(\bar v)$ is satisfied by definition. Notice that the condition given in the proposition, $\bar v( \AC\setminus a')=0$, implies that $\bar v(\CC\setminus a')=0$ for all $\CC\subset\AC.$ Using this, inequality constraints are divided into two sets of constraints as $\sum_{a\in\CC}\hat\beta_a\geq \bar v(\CC)$, for all $\CC\ni a'$, and $\sum_{a\in\CC}\hat\beta_a\geq 0$, for all $\CC\not\ni a'$. The first set of inequalities are satisfied, since $\bar v(\AC)\geq \bar v(\CC)$, for all $\CC\ni a'$. The second set of inequalities are satisfied, since $\hat\beta\in\mathbb{R}^{\AC}_+$. This concludes that $\hat\beta\in K_\text{Core}(\bar v)$, and hence the core is nonempty. \hfill\QEDA

According to the proof above, an area that satisfies the condition in Proposition~\ref{prop:corempt} has a right to veto any coalitional deviation that does not include it. The proof constructs the benefit allocation that assigns all the expected cost reduction to this area and then shows that the remaining areas do not have any incentives to form coalitions.
\section{Proof of Proposition~\ref{prop:leastcore_volp}}\label{app:a4}

	For this proof, we extend the method in \citep[Theorem 2.7]{maschler1979geometric} by taking into account that our coalitional value function is defined by $\bar v(\CC) = J(\emptyset) - J(\CC)\geq 0,$ for all $\CC\subset \AC$. We prove by contradiction. Let $\hat \beta \in K_\text{Core}(\bar v,\epsilon^*(\bar v))$ be a benefit allocation with $\hat \beta_{a'}<0$. In this case, we show that there exists $\epsilon <\epsilon^*(\bar v)$ such that $K_\text{Core}(\bar v,\epsilon)\neq \emptyset$. This would contradict the definition of the least-core.
	
	Since $\hat \beta \in K_\text{Core}(\bar v,\epsilon^*(\bar v))$ and $\bar v(a')=J(\emptyset) - J(a')=0$, we have $0>\beta_{a'}\geq \bar v(a') -\epsilon^*(\bar v)=-\epsilon^*(\bar v)$. Notice that, for any $\CC\not\ni a'$, we have $\sum_{a\in\CC}\hat\beta_a\geq \bar v(\CC\cup a')-\epsilon^*(\bar v)-\beta_{a'}$ by adding and subtracting $\beta_{a'}$, and by using the fact that $\epsilon^*(\bar v)>0$ for the case corresponding to $\CC\cup a'=\AC$. Since $\bar v$ is nondecreasing and $\beta_{a'}<0$, we obtain $\sum_{a\in\CC}\hat\beta_a> \bar v(\CC)-\epsilon^*(\bar v)$. 
	
	Notice that we can always find a small positive number $\delta$ such that  $\sum_{a\in\CC}\hat\beta_a - |\CC|\delta> \bar v(\CC)-\epsilon^*(\bar v)+\delta$ holds for any $\CC\not\ni a'$. Next, we show that $K_\text{Core}(\bar v,\epsilon^*(\bar v)-\delta)$ is nonempty for this particular choice. Define $\bar{\beta}$ such that $\bar\beta_a=\hat \beta_a - \delta$ for all $a\neq a'$and $\bar{\beta}_{a'}=\hat{\beta}_{a'}+(|\AC|-1)\delta$. This new allocation $\bar{\beta}$ clearly satisfies the equality constraint in $K_\text{Core}(\bar v,\epsilon^*(\bar v)-\delta)$. For inequality constraints $\CC\not\ni a'$, we have $\sum_{a\in\CC}\bar\beta_a=\sum_{a\in\CC}\hat\beta_a - |\CC|\delta> \bar v(\CC)-\epsilon^*(\bar v)+\delta$, where the strict inequality follows from the definition of $\delta$. For inequality constraints $\CC\ni a'$, we have $\sum_{a\in\CC}\bar\beta_a\geq\sum_{a\in\CC}\hat\beta_a+\delta\geq \bar v(\CC)-\epsilon^*(\bar v)+\delta$. Hence, $\bar \beta \in K_\text{Core}(\bar v,\epsilon^*(\bar v)-\delta)$, in other words, $K_\text{Core}(\bar v,\epsilon^*(\bar v)-\delta)\neq\emptyset.$ This contradicts $K_\text{Core}(\bar v,\epsilon^*(\bar v))$ being the least-core. \hfill\QEDA

\section{MILP reformulation of problem \eqref{eq:ba_findviolation} and the constraint generation algorithm}\label{app:a6}
 Since we have~$\bar v(\CC)=J(\emptyset)-J(\CC)$, optimal solutions to problem \eqref{eq:ba_findviolation} coincide with the ones to the following problem:
\begin{flalign}\label{eq:ba_findviolation_ref}
	&\underset{ \CC\subseteq\AC}{\text{minimize}}\ J(\CC)+\sum_{a\in\CC}\beta_a^k .
\end{flalign}
The problem above is given by the following MILP:
\begingroup
\allowdisplaybreaks
\begin{subequations} \label{mod:B-Pre-k}
\begin{align}
& \bar J(\beta^k)=\,\underset{ \Phi_{\text{PR}}^k }{\text{minimize}}\quad \sum_{i \in \mathcal I}\left(C^{+}_{i} r^{+}_{i} + C^{-}_i r^{-}_{i}\right) + \sum_{i\in \mathcal I}{C}_{i}{p}_{i} \nonumber\\
& \hspace{1.5cm}+ \sum_{s\in\mathcal S} \pi_s \Big[ \sum_{i \in \mathcal I}{C}_{i}\left({p}_{is}^{+}-{p}_{is}^{-}\right)+\sum_{n \in \mathcal N} {C}^{\text{sh}}l^{\text{sh}}_{ns} \Big]+\sum_{a\in\CC}b_a\beta_a^k \label{eq:B-pre-k-obj}
\end{align}
subject to
\begin{align}
 &\hspace{-.25cm} b_a\in\{0,1\},\quad\forall a\in \AC,\label{eq:preempt_bin}\\
		&\hspace{-.25cm}  (1-b_{a_r(e)})\chi_e \le \chi_{e}' \le (1-b_{a_r(e)})\chi_e+b_{a_r(e)}\ \text{and}\nonumber\\
		&\hspace{-.25cm}  (1-b_{a_s(e)})\chi_e \le \chi_{e}' \le (1-b_{a_s(e)})\chi_e+b_{a_s(e)},  \forall e\in\mathcal E,\label{eq:preempt_bin_ta}\\
&\hspace{-.25cm}   \text{Constraints} \;\; \eqref{eq:A-CV-rt-bal-constr} -\eqref{eq:A-CV-rt-flowAC-max}\text{ and }\eqref{eq:A-CV-ref-node},\quad \forall s\in \mathcal S,\label{eq:B-pre-k-CV}\\
	&\hspace{-.25cm} - b_{a_r(\ell)}M  \le f_{\ell s}-{f}_{\ell} \le b_{a_r(\ell)} M \ \text{and}\nonumber\\
	&\hspace{-.25cm} - b_{a_s(\ell)}M  \le f_{\ell s}-{f}_{\ell} \le b_{a_s(\ell)} M ,\nonumber\\  &\hspace{-.25cm}\forall \ell\in\cup_{e\in\mathcal E(\chi)}\Lambda_{e},\quad \forall s\in \mathcal S,\label{eq:preempt_resact}\\
&\hspace{-.25cm}  \text{Constraints} \;\; \eqref{eq:LLR}\  \text{and}\ \eqref{eq:LLD}.
\end{align}
\end{subequations}
\endgroup
where $M$ is a large positive number, $\Phi_{\text{PR}}^k\allowbreak=\allowbreak\{b_a,\allowbreak \forall a \cup\chi_e',\allowbreak \forall e \cup \Phi_{\text R} \cup \Phi_{\text D} \cup  \Phi_{\text B}^{s}, \allowbreak\forall s\}$  is the set of primal optimization variables.
For the sake of brevity, the KKT conditions and the Lagrange multipliers of the lower-level optimization problems are omitted, see~\ref{app:a}.

In problem~\eqref{mod:B-Pre-k}, the parameter ${\beta}_a^k$ can be considered as the activation fee of area~$a$ for participating in the preemptive model. Given such activation fees, problem~\eqref{mod:B-Pre-k} then finds the optimal set of participants for the preemptive model. Notice that this problem involves $|\AC|$ more binary variables than problem~\eqref{mod:B-Pre}. In Section~\ref{sec:IEEE_RTS}, the numerical case studies illustrate that the computation times are still similar for both problems. Let $\{\hat{b}_a, \forall a\}$ denote the optimal binary solution. Finally, we have $\CC^k=\{a\in\AC\,|\,\hat{b}_a=1\}$ and $\bar v(\CC^k)=J(\emptyset)-(\bar J(\beta^k)-\sum_{a\in\CC^k}\beta^k_a)$. 

We can summarize this iterative algorithm as follows.

 \begin{algorithm}[h]
 	\caption{Constraint Generation Algorithm for the Least-Core Selecting Mechanism}
	\begin{algorithmic}[1]\label{alg:const_gen}
 		\renewcommand{\algorithmicrequire}{\textbf{Initialize:} Compute $\bar v(\AC)$,}
 		\REQUIRE 
 		and $\beta^{\text{c}}$, set $k=0$, $\eta_0=1$, $\epsilon_0=0$, initialize $\mathcal{F}^1$ (e.g., $\mathcal{F}^1=\emptyset$). \\
 		\WHILE {$\eta^{k}> \epsilon^k$}
 		\STATE $k \leftarrow k+1$
 		\STATE Obtain $\epsilon^{k}$ and $\beta^{k}$ by solving~\eqref{mod:candidate_k} and~\eqref{mod:candidate_k_unique}, and obtain $\eta^{k}$, $\CC^k$, and $\bar v(\CC^k)$ by solving~\eqref{eq:ba_findviolation}.
		\STATE Set $\mathcal{F}^{k+1}=\mathcal{F}^k\cup\CC^k$.
 		\ENDWHILE\ and \textbf{return} $\beta^k=\hat{\beta}(\bar v,\beta^{\text{c}})$.
 	\end{algorithmic}
 \end{algorithm}

\section{Proof of Proposition~\ref{prop:empty}}\label{app:a7}
	Assume $\beta\in K_\text{Core}(v^s)$. We now prove that this yields a contradiction. Notice that $K_\text{Core}(v^s)\subset\mathbb{R}_+^{\AC}$ since $\beta_{a}\geq v^s(a)=0.$ The scenario-specific core also implies $$\sum_{a\in \AC\setminus \CC}\beta_{a}\leq v^s(\AC)- v^s(\CC)<0.$$ The inequality above follows from combining $\sum_{a\in\AC}\beta_{a}=v^s(\AC)$ and $\sum_{a\in\CC}\beta_{a}\geq v^s(\CC)$. We obtained a contradiction. \hfill\QEDA
	
\section{Market outcomes for three-area example}\label{app:market_outcome}	

The market outcomes for the existing sequential market model and the preemptive model are provided in Table~\ref{tab:three_marketout}. For the preemptive model, the values of $r^-$ and $r^+$ shown in parenthesis indicate the amount of reserves destined to meet the requirements of the neighboring areas. Notice that inflexible generators and wind power generators are not capable of providing any reserves. In the day-ahead stage, the quantities correspond to $p_i$ and $w_j$ for conventional units and wind power generators, respectively. Because of the null production costs, wind power generators are always utilized at their full capacity given by the expected value of the stochastic process. In the balancing stage, the quantities correspond to the changes in the production levels with respect to the ones assigned at the day-ahead stage. Specifically, they are given by $({p}_{is}^{+}-{p}_{is}^{-})$ and $(W_{j s}-{w}_j-w^{\text{spill}}_{j s})$ for conventional units and wind power generators, respectively. It can be verified that there is $31$ MW and $15.6$ MW curtailment in the wind power production in the existing market model for scenarios $s_1$ and $s_2$. This is completely eliminated when the preemptive model is used. 

Observe that the implementation of the preemptive model results in a higher reserve allocation for the lower cost generator at node $8$. Consequently, the day-ahead market takes these new reserve quantities into account while deciding on the day-ahead energy quantities. Notice that there is $9$ MW load shedding in scenario $s_2$ in the existing market model. The preemptive model ensures that we do not resort to any costly load shedding in both scenarios.

\begin{table*}[h]\footnotesize
\centering
\caption{Market outcomes for the existing sequential market and the preemptive model (in MW)}
\begin{tabular}{|c||c|c|c|c|c||c|c|c|c|c|}
\hline
 Models &  \multicolumn{5}{c||}{Existing market model} & \multicolumn{5}{c|}{Preemptive model}\\ \hline
 Trading floors &  \multicolumn{2}{c|}{Reserve} &  Day-ahead & \multicolumn{2}{c||}{Balancing}&  \multicolumn{2}{c|}{Reserve} &  Day-ahead & \multicolumn{2}{c|}{Balancing}\\ \hline
  &  $r^-$&  $r^+$ &   & $s_1$ & $s_2$ &$r^-$&  $r^+$ &   & $s_1$ & $s_2$   \\ \hline
   $i_1$ &$0$  &$0$ &$120$&$0$ &$0$ &$0$ &$0$&$120$ &$0$ &$0$ \\ \hline
   $i_2$ &$8$  &$12$ &$38$&$3$ &$3$ &$8$ &$12$&$38$ &$7.7$ &$7.2$ \\ \hline
   $i_3$ &$0$  &$0$ &$0$  &$0$ &$0$ &$0$ &$0$ &$0$ &$0$ &$0$ \\ \hline
   $i_4$ &$0$  &$0$ &$120$&$0$ &$0$ &$0$ &$0$&$120$ &$0$ &$0$ \\ \hline
   $i_5$ &$9.6$&$6.4$&$33$&$6$ &$6$ &$7.2$ &$4$ &$31.1$ &$-7.2$ &$-7.2$ \\ \hline
   $i_6$ &$0$  &$0$ &$0$  &$0$ &$0$ &$0$ &$0$&$0$ &$0$ &$0$ \\ \hline
   $i_7$ &$0$  &$0$ &$120$&$0$ &$0$ &$0$ &$0$&$120$ &$0$ &$0$ \\ \hline
   $i_8$ &$8$&$12$ &$38$&$12$&$12$&$10.4\,(2.4)$&$14.4\,(2.4)$&$35.6$ &$-10.1$ &$14.4$ \\ \hline
   $i_9$ &$0$  &$0$ &$6.6$&$0$ &$0$ &$0$   &$0$ &$10.9$ &$0$ &$0$ \\ \hline
    $j_3$ &$0$ &$0$ &$42$ & $-3$  &$-12$ &$0$ &$0$ &$42$ &$8$ &$-12$ \\ \hline
   $j_6$ &$0$  &$0$ &$70.4$&$-6.4$&$-6$&$0$ &$0$ &$70.4$ &$-6.4$ &$9.6$ \\ \hline
   $j_9$ &$0$  &$0$ &$42$ & $-12$ &$-12$ &$0$ &$0$ &$42$ &$8$ &$-12$ \\ \hline
\end{tabular}\label{tab:three_marketout}
\end{table*}

In Table~\ref{tab:three_jempty}, we provide the consumer and producer surpluses (CS and PS), and the congestion rents (CR) allocated to each area under both scenarios in the existing sequential market. Using these allocations for all trading floors, we defined a budget balanced cost allocation in the numerical case studies.

\begin{table}[h]\footnotesize
\centering
\caption{Cost allocations in all trading floors for the existing sequential market (in \euro)}
\begin{tabular}{|c||c|c|c|}
\hline
 Areas & Area $1$  & Area $2$ & Area $3$\\ \hline
 CS for~\eqref{mod:A-Res} & $-60.0$  & $-64.0$ & $-70.0$\\ \hline
 PS for~\eqref{mod:A-Res}& $0$  & $0$ & $0$\\ \hline
 CR for~\eqref{mod:A-Res}& $0$  &$0$ & $0$\\ \hline
 CS for~\eqref{mod:A-DA} & $-8{,}448.0$  & $-7{,}239.0$ & $-9{,}900.0$\\ \hline
 PS for~\eqref{mod:A-DA}& $4{,}159.6$  & $3{,}750.2$ & $4{,}392.0$\\ \hline
 CR for~\eqref{mod:A-DA}& $0$  & $99.0$ & $99.0$\\ \hline
  CS for~\eqref{mod:A-CV-rt} in~$s_1$ & $0$  & $0$ & $0$\\ \hline
 PS for~\eqref{mod:A-CV-rt} in~$s_1$& $0$  & $-6{,}400.0$ & $5{,}250.0$\\ \hline
 CR for~\eqref{mod:A-CV-rt} in~$s_1$& $0$  & $0$ & $0$\\ \hline
  CS for~\eqref{mod:A-CV-rt} in~$s_2$ & $0$  & $0$ & $0$\\ \hline
 PS for~\eqref{mod:A-CV-rt} in~$s_2$& $-12{,}000.0$  & $0$ & $2{,}250.0$\\ \hline
 CR for~\eqref{mod:A-CV-rt} in~$s_2$ & $0$  & $0$ & $0$\\ \hline
 $J^{s_1}_a(\emptyset)$& $4{,}348.4$  & $9{,}853.8$ & $229.0$\\ \hline
 $J^{s_2}_a(\emptyset)$& $16{,}348.4$  & $3{,}453.8$ & $3{,}229.0$\\ \hline
\end{tabular}\label{tab:three_jempty}
\end{table}
\section{Impact of the network topology on benefit allocations in the three-area example}\label{app:three_area_netw_top}
\begin{figure}[h]
	\centering
	\begin{tikzpicture}[scale=1, every node/.style={scale=0.5}]
			\node at (0,0) {\includegraphics[width=1.2\textwidth]{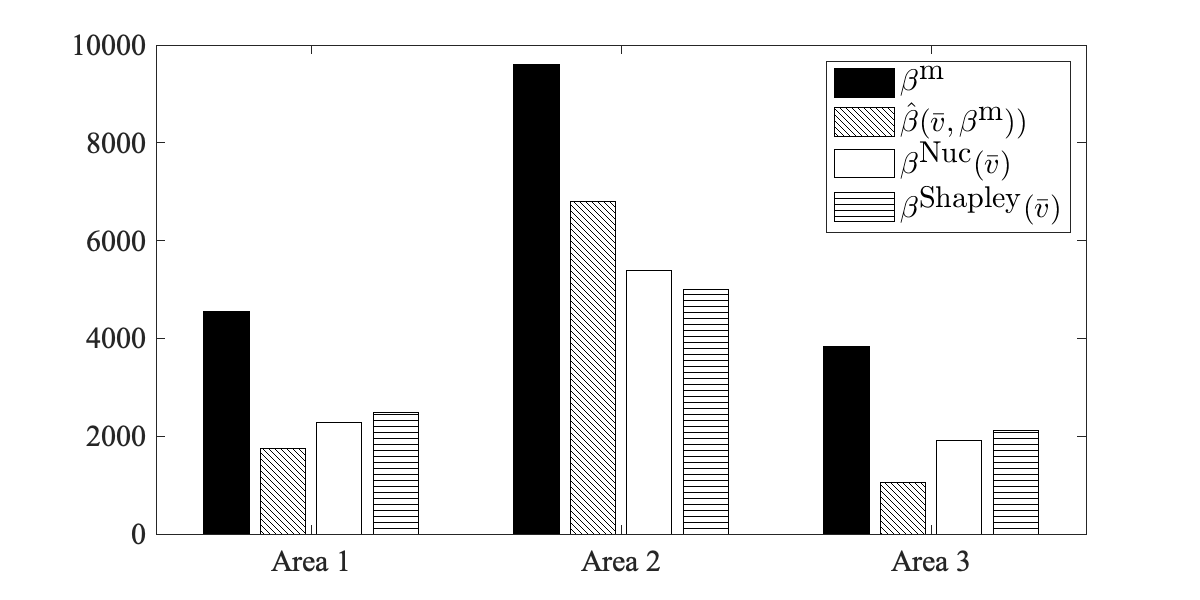}};
	\end{tikzpicture}
	\caption{Benefit allocations after removing the flexibility and the uncertainty from area 2 (in \euro)}\label{fig:three_removeflex}
\end{figure}

In this example, our goal is to illustrate how the network topology and the specific location of each area in the electricity network affects the benefits allocated to this area. To this end, we modify our base model by removing the wind generator from area 2 and by changing all the units in area 2 to be inflexible. Benefit allocations for the expected cost reduction are provided in Figure~\ref{fig:three_removeflex}. We highlight that removing a zero marginal cost wind generator increased the costs globally. In this setup, area 2 receives the highest benefits under every allocation mechanism, even though this area is not capable of directly participating in any reserve exchange. Nevertheless, the role of area 2 in this network arrangement is instrumental, since it enables the coordination between areas 1 and 3 because of its location in the area graph.

As a remark, the core is nonempty, since the condition in Proposition~\ref{prop:corempt} is satisfied. The Shapley value is also in the core, since we verified that the coalition value function is supermodular.

\section{IEEE RTS layout}\label{app:IEEERTS96_area}

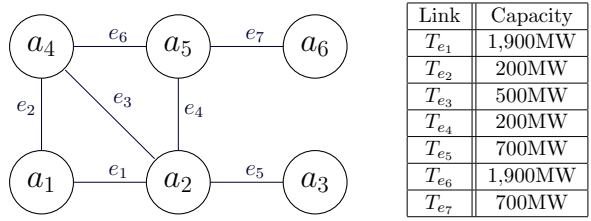
\begin{figure}[h]
	\centering
	\begin{tikzpicture}[scale=1.2, every node/.style={scale=0.72}]
    \draw[-,black!80!blue,line width=.1mm] (0.35,0) -- (1.15,0) node[above=.2cm, left=.2cm] {\large $e_1$};
    \draw[-,black!80!blue,line width=.1mm] (0,0.35) -- (0,1.15) node[above=-.4cm, left=-1.8cm] {\large $e_3$};
    \draw[-,black!80!blue,line width=.1mm] (1.5,0.35) -- (1.5,1.15) node[above=-.6cm, left=-.6cm] {\large $e_4$};
    \draw[-,black!80!blue,line width=.1mm] (1.24,0.26) -- (0.26,1.24) node[above=-.7cm, left=.4cm] {\large $e_2$};
    \draw[-,black!80!blue,line width=.1mm] (1.85,0) -- (2.65,0) node[above=.2cm, left=.2cm] {\large $e_5$};
    \draw[-,black!80!blue,line width=.1mm] (1.85,1.5) -- (2.65,1.5) node[above=.2cm, left=.2cm] {\large $e_7$};
    \draw[-,black!80!blue,line width=.1mm] (0.35,1.5) -- (1.15,1.5) node[above=.2cm, left=.2cm] {\large ${e_6}$};
    %\draw[-,black!80!blue,line width=.1mm] (1.7,1.7) -- (2,1);
    %\draw[-,black!80!blue,line width=.1mm] (1.7,1.7) -- (1,0);
    %\draw[-,black!80!blue,line width=.1mm] (0.35,1.7) -- (1,0);
    \draw (0,0) circle (.35cm) node {\LARGE $a_1$};
    \draw (1.5,0) circle (.35cm) node {\LARGE $a_2$};
    \draw (3,0) circle (.35cm) node {\LARGE $a_3$};
    \draw (0,1.5) circle (.35cm) node {\LARGE $a_4$};
    \draw (1.5,1.5) circle (.35cm) node {\LARGE $a_5$};
    \draw (3,1.5) circle (.35cm) node {\LARGE $a_6$};
    \node at (5,.8) {%
  \begin{tabular}{|c||c|}
\hline
 Link & Capacity \\ \hline
 $T_{e_1}$& $1{,}900$MW \\ \hline
 $T_{e_2}$& $200$MW \\ \hline
 $T_{e_3}$& $500$MW \\ \hline
 $T_{e_4}$& $200$MW \\ \hline
 $T_{e_5}$& $700$MW \\ \hline
 $T_{e_6}$& $1{,}900$MW \\ \hline
 $T_{e_7}$& $700$MW \\ \hline
\end{tabular}};
	\end{tikzpicture}
	\caption{The area graph for the IEEE RTS}\label{fig:IEEERTS96_area_graph}
\end{figure}

The layout can be found in Figure~\ref{fig:IEEERTS96_area}. The area graph is illustrated in Figure~\ref{fig:IEEERTS96_area_graph}.
\input{Fig-IEEE-24.tex}

\section{Modifications of IEEE RTS case study}\label{app:modif_rts}
\subsection{Impact of the wind power penetration levels on benefit allocations}

\begin{figure}[h]
	\centering
	\begin{tikzpicture}[scale=1, every node/.style={scale=0.5}]
			\node at (0,0) {\includegraphics[width=1.2\textwidth]{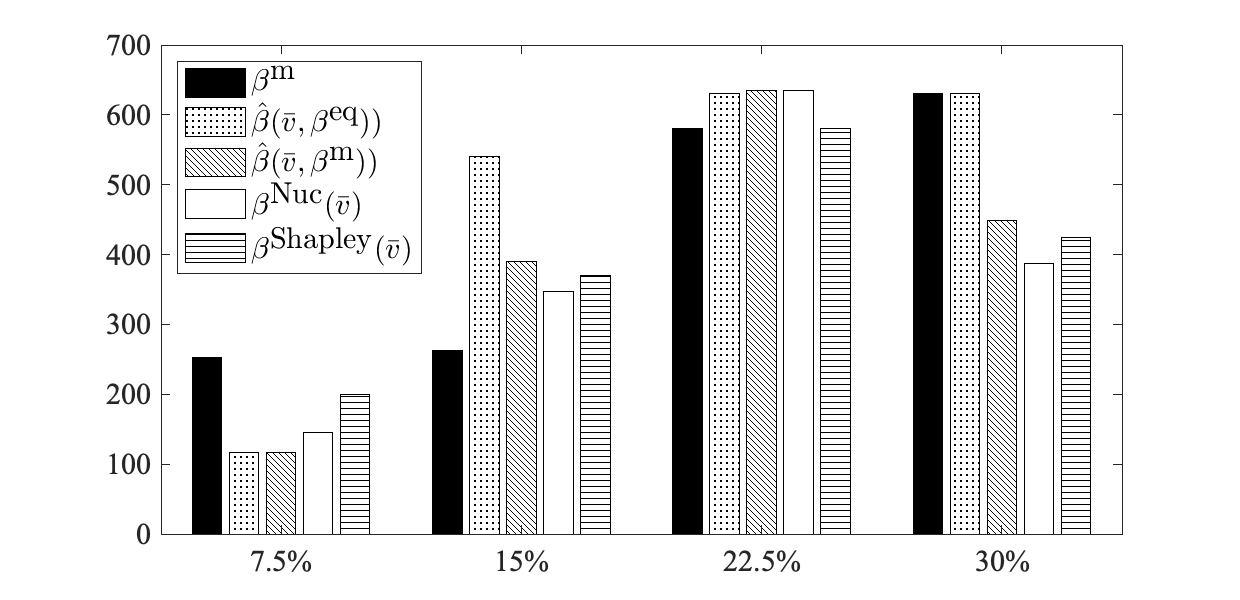}};
	\end{tikzpicture}
	\caption{Benefit allocations of area 1 for the different levels of wind power penetration in area 1 (in \euro)}\label{fig:rts_benefit_2}
\end{figure}

The wind power penetration level of an area is defined as the ratio between the expected wind power production and the total demand of that area. In this example, we change the level of wind power penetration for area 1. The default value is given by $30\%$ from the previous section. The resulting benefits allocations for area~1 are provided in Figure~\ref{fig:rts_benefit_2}. For all efficient benefit allocation mechanisms, observe that the benefits initially increase and then decrease with the wind power penetration level. 
It is generally hard to anticipate such changes in the benefits, since the wind power generation has two impacts acting in opposite directions. {On the one hand}, it has null production cost bringing in low cost energy to the coalition. {On the other hand}, it increases the need for the reserve and balancing services. As a remark, in this study, the core is empty for the levels $15\%$ and $22.5\%$.

\subsection{Impact of the available flexibility on benefit allocations}
\begin{figure}[h]
	\centering
	\begin{tikzpicture}[scale=1, every node/.style={scale=0.5}]
			\node at (0,0) {\includegraphics[width=1.3\textwidth]{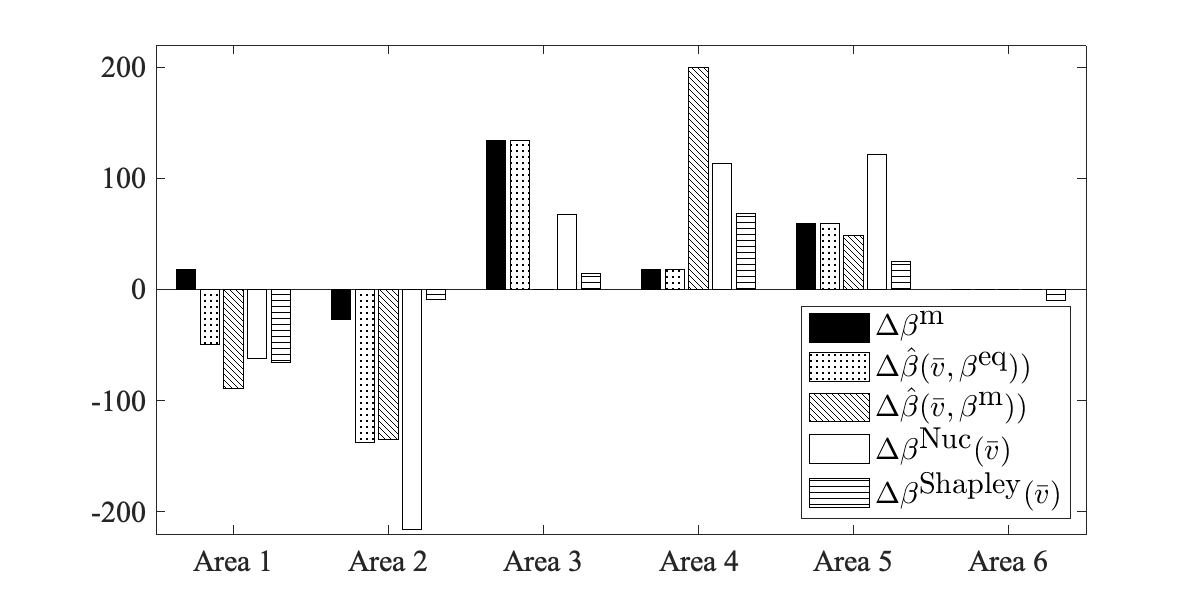}};
	\end{tikzpicture}
	\caption{The changes in benefit allocations after increasing the flexibility in area 4 (in \euro)}\label{fig:rts_benefit_3}
\end{figure}

In order to assess the impact of {available flexibility}, we {double} the capacities of all flexible generators in area 4. The changes in the benefits allocations are provided in Figure~\ref{fig:rts_benefit_3}, denoted by $\Delta \beta$. The figure shows that the benefits allocated to area~4 are increased under every allocation mechanism compared to Figure~\ref{fig:rts_benefit}. All mechanisms account for the flexibility offered by the generators in area~4 both in the reserve capacity and the balancing markets. Note that an increase in the capacity of flexible generators reduces also the day-ahead system cost both with or without the implementation of the preemptive model. Finally, predicting the changes in the benefits of the other areas is quite difficult. We observe that this additional flexibility replaced the flexibility in areas 1 and 2, while increasing the contributions of areas 3 and 5 in subcoalitions.

\end{document}